\documentclass[11pt]{article}
\usepackage{amsmath, latexsym, epsf, amssymb, framed}
\usepackage{enumerate}
\usepackage[margin =1in]{geometry}
\usepackage{amssymb,amsthm}
\usepackage{amsmath}
\usepackage{bbm}
\usepackage{enumerate}
\usepackage{float}
\usepackage{color}
\usepackage{graphicx}
\usepackage{wrapfig}
\usepackage{pgf,tikz,pgfplots}
\usepackage{cancel}
\usepackage{tikzsymbols}
\usepackage[bb=boondox]{mathalfa}
\pgfplotsset{compat=1.15}
\usepackage{mathrsfs}
\usetikzlibrary{arrows}
\usepackage{hyperref}
\hypersetup{hidelinks}
\usepackage{mathtools}

\usepackage{subcaption}
\usepackage{bm}
\usepackage{authblk}

\newcommand{\inclu}[0] {\ar@{^{(}->}}

\graphicspath{ {./figures/} }

\newtheorem*{theorem*}{Theorem}

\usepackage{algorithm}
\usepackage{algpseudocode}
\usepackage{tabularx}
\usepackage{hyperref}
\usepackage{xurl}
\usepackage{caption}
\usepackage{cleveref}

\newcommand{\objcref}[1]{objective~(\ref{#1})}
\newcommand{\Objcref}[1]{Objective~(\ref{#1})}

\newcommand{\constrcref}[1]{constraint~(\ref{#1})}
\newcommand{\Constrcref}[1]{Constraint~(\ref{#1})}

\newcommand{\constrscref}[1]{constraints~(\ref{#1})}
\newcommand{\Constrscref}[1]{Constraints~(\ref{#1})}

\newcommand{\nodeU}{u}
\newcommand{\nodeV}{v}

\newcommand{\graphT}{G_T}
\newcommand{\nodesT}{V_T}
\newcommand{\edgesT}{E_T}
\newcommand{\graphB}{G_B}
\newcommand{\nodesB}{V_B}
\newcommand{\edgesB}{E_B}

\newcommand{\graphM}{G_M}
\newcommand{\nodesM}{V_M}
\newcommand{\edgesM}{E_M}
\newcommand{\edgesFirst}{E^F}
\newcommand{\edgesLast}{E^L}
\newcommand{\nodesOri}{V^O}
\newcommand{\nodesDest}{V^D}

\newcommand{\lineSet}{\mathcal{L}}
\newcommand{\oneline}{l}
\newcommand{\graphline}{P_l}
\newcommand{\nodesline}{V_l}
\newcommand{\edgesline}{E_l}

\newcommand{\requestSet}{\mathcal{D}}

\newcommand{\ori}{s}

\newcommand{\dest}{t}
\newcommand{\requestST}{d_{s,t}}

\newcommand{\costUV}{c_{u,v}}
\newcommand{\costL}{c_l}
\newcommand{\capL}{\kappa_l}
\newcommand{\capacity}{\kappa}
\newcommand{\freqL}{\mathcal{F}_l}
\newcommand{\lbL}{\tau_l}

\newcommand{\scaleBus}{\gamma}
\newcommand{\scaleTaxi}{\alpha}
\newcommand{\budget}{B}

\newcommand{\subEdge}{h}
\newcommand{\source}{S}
\newcommand{\sink}{T}

\newcommand{\detourConst}{\theta}
\newcommand{\maxLength}{\Omega}

\newcommand{\numTaxi}{x}
\newcommand{\numBus}{Y}
\newcommand{\directTrip}{q}
\newcommand{\getOn}{y}
\newcommand{\getOff}{w}
\newcommand{\toDest}{z}
\newcommand{\busFlow}{f}

\newcommand{\busDual}{p}
\newcommand{\taxiDual}{k}
\newcommand{\budgetDual}{\beta}
\newcommand{\capDual}{r}

\newcommand{\taxi}{on-demand}
\newcommand{\bus}{bus}
\newcommand{\Taxi}{On-demand}

\newcommand{\Bus}{Bus}

\newcommand{\mip}{\mathcal{M}}
\newcommand{\mipPrime}{\mathcal{M}^{\prime}}
\newcommand{\pricing}{\mathcal{P}}
\newcommand{\pricingPrime}{\mathcal{P}^{\prime}}

\newtheorem{Assumption}{Assumption}
\newtheorem*{Example*}{Example}

\newcommand{\ttilde}[1]{\smash{\tilde{#1}}}

\allowdisplaybreaks

\title{Integrated Multi-Modal Transit Network Design \\ via Dynamic Line Generation}

\author{Ning Duan\thanks{School of Operations Research and Information Engineering, Cornell University} \hspace{1cm} 
Oktay G\"unl\"uk\thanks{School of Industrial and Systems Engineering, Georgia Institute of Technology} \hspace{1cm} 
Samitha Samaranayake\thanks{School of Civil and Environmental Engineering, Cornell University}}
\date{}

\begin{document}
\maketitle

\abstract{The integration of fixed-route public transit and on-demand mobility services presents both a modeling challenge and a computational opportunity for large-scale network design. We propose a flow-based mixed-integer programming formulation that jointly optimizes transit line planning and service frequencies while explicitly capturing first- and last-mile connectivity via on-demand services, under a fixed operating budget. To achieve tractability at urban scale, we develop a novel column generation heuristic scheme with tailored pricing subproblems. Applied to networks and demand in Boston and Chicago, the framework yields operationally feasible designs that substantially increase demand served. Relative to transit-only and on-demand-only baselines, ridership increases by up to 20.99\% and 93.58\% in Boston, and by up to 10.63\% and 149.84\% in Chicago. Compared to a multi-modal benchmark, our approach improves ridership by 5.42\% and 5.80\% in Boston and Chicago, respectively. These results demonstrate that (i) joint co-design of transit routes, frequencies, and on-demand legs within a unified optimization framework yields substantially greater ridership than single-mode or decoupled approaches under equivalent budget constraints; and (ii) the proposed formulation and column generation pricing scheme admit tractable, operationally feasible, high-performing solutions relative to tested baselines.}

\section{Introduction}

Efficient, affordable and ubiquitous access to personal mobility is a fundamental societal need, with broad implications for personal well-being, economics, education, and public health~\cite{Weber2022}. However, current transportation systems, dominated  by personal automobiles, fail to deliver on these fronts. The reliance on single-occupancy personal vehicles has led to severe gridlock across major metropolitan areas, with congestion costs in the United States alone reaching \$224 billion in 2022~\cite{congestionDollar}, approximately 0.9\% of GDP, including 8.5 billion hours lost to traffic and 3.3 billion gallons of wasted fuel. These figures do not account for other significant negative externalities such as greenhouse gas emissions and particulate matter pollution~\cite{congestionEmission}, travel-time uncertainty~\cite{congestionUncertainty}, and increased accident risk~\cite{Berhanu2023}.  

Over the past decade, new mobility technologies have emerged with the promise of more efficient and sustainable urban mobility. The large-scale adoption of ride-hailing has expanded access to mobility,
particularly for people who cannot or prefer not to drive, and reduced
drunk driving \cite{Misra2022,Dills2018}. However, ride-hailing services largely inherit---and in some cases exacerbate---the inefficiencies of private vehicle usage, particularly through deadheading, i.e., empty repositioning trips between passenger rides~\cite{sejourne2018price, balding2019estimated}. While ride-pooling offers a partial mitigation by increasing vehicle occupancy, its effectiveness is limited in low-density environments and remains constrained by operational costs. Conventional mass transit systems, especially \bus~networks, offer significantly higher capacity and lower per-passenger externalities~\cite{Higashide2019Better}. Nevertheless, deploying high-quality \bus~services faces many challenges---efficient operations require high-demand corridors, yet the public nature of \bus~service demands broad accessibility, especially for low-income populations who have been pushed to the outskirts of urban areas due to the suburbanization of poverty~\cite{Wang2017Relationship}. This tension is further exacerbated by the spatial structure of many U.S. cities. 

Recently, there have been attempts to enhance transit systems with the multi-modal integration of ride-hailing and micro-transit services as first-mile and last-mile connectors to extend their reach and efficiency~\cite{Brown2021Can}. While these approaches can extend the reach of transit, they are limited by high costs per rider (due to high labor costs) and are constrained by the existing mass-transit network, which was not designed to be part of a multi-modal system.

Motivated by existing limitations, we propose end-to-end design of multi-modal transit systems. In our setting, we use “\bus” to denote the transit system, although the proposed methodology readily generalizes to other types of transit. We use “\taxi” to denote the flexible transportation mode in our system, referring to single-seat vehicle services such as taxis. In this context, we exploit both the efficiency of high-occupancy, high-frequency \bus~corridors and the flexibility of \taxi~services. We consider a transportation planning problem where a transit system planner aims to design a service that integrates conventional fixed-route \bus~lines and  \taxi~components. Each traveler can be transported from their origin to their destination either by \bus,~\taxi~or a combination of the two modes (e.g., \taxi ~to \bus  ~to \taxi~or a single \taxi~leg). Given a road network, origin-destination travel demand, \bus~capacity, \taxi~and \bus~operating costs, and an operating budget, we develop a new flow-based mixed-integer programming (MIP) model for the problem and design novel column generation pricing subproblems that price candidate lines using the dual values of the master problem's LP relaxation. We jointly optimize service routes and frequencies to maximize ridership. Our main contribution advances the field by demonstrating a framework that can solve the multi-modal \bus~line planning problem at urban scale, utilizing travel demand and road network data from major U.S. metropolitan areas.

\subsection{Contributions}
The multi-modal transit network planning problem poses significant modeling and algorithmic challenges. We summarize our main contributions below.

\subsubsection{Realistic Modeling}
We propose a novel flow-based mixed-integer programming formulation for the joint design of \bus~lines and operating frequencies in a multi-modal transit system. We introduce behavioral and service-quality constraints that are not captured by standard formulations. We assume that the transit system must satisfy a number of constraints related to passenger acceptance and quality of service. First, the resulting routes should not include  many transfers, as passengers associate a high level of dissatisfaction with transfers \cite{Grise2019Transferring}. This quality-of-service constraint precludes the use of the standard multi-commodity flow models utilized in the telecommunication network literature \cite{Gunluk1999A, Bienstock1996Capacitated}, where data packets are not assigned to fixed routes and essentially can take any path in the network, implying unlimited transfers. It also eliminates the use of multi-commodity flow models in \bus-only design for the same reason \cite{Borndorfer2007}. Second, passengers prefer shorter trips and are less willing to use transit if the routes are much longer than a direct personal vehicle or \taxi~trip. This includes the waiting time for a \bus~to arrive, which necessitates \bus~routes to operate at a high frequency or short headways (e.g., every 10--15 minutes)---studies show that passengers have a greater disutility for time spent waiting for a \bus~compared to time spent on a \bus~\cite{Litman2008Valuing, Fan2016Waiting}. Consequently, each operational \bus~line  needs a minimum number of operating vehicles that depends on the length of the line, to ensure the required operating frequency. Third, we ensure these \bus~lines are efficient and acceptable to riders by requiring that their end-to-end distance is no more than double the shortest-path distance between their endpoints, while also adhering to maximum length constraints.

\subsubsection{Tailored Column Generation}
A challenge in constructing transit lines is the design of column generation pricing subproblems in the absence of explicit dual variables for line-specific constraints. We address this difficulty by introducing  subproblems that approximate the reduced cost of candidate lines, enabling effective column generation despite the lack of corresponding dual information.

\subsubsection{Computational Experiments on Urban-Scale Instances}
To our knowledge, this work is the first to jointly generate end-to-end lines, frequencies, and \taxi~legs from scratch at urban scale. To ensure scalability at the urban level, we also develop a suite of LP-based heuristics integrated within both the pricing subproblems and the master MIP. In each column generation iteration, we restrict the search to a reduced network subgraph prioritized by dual information from the LP relaxation. After generating candidate lines, we apply a pre-selection heuristic to identify the most promising subset for the master problem, significantly reducing the computational burden arising from the interaction of continuous, integer, and semi-continuous decision variables.

\subsection{Related Literature}

The line planning problem in transit networks is a classical network design problem and is computationally challenging to solve. The problem determines a set of transit lines and their service frequencies over a given transportation network, typically balancing passenger service quality and operating cost.~\cite{schobel2012line} provides a comprehensive review of models and algorithms for the problem. In the passenger-oriented setting considered in this paper, a transit planner with a fixed operating budget and known origin-destination demand aims to find a set of bus routes and their corresponding frequencies such that as many passengers as possible can be served. Existing mathematical-programming-based approaches to line planning differ in how they model passenger service and the feasible set of lines. Many formulations assume the existence of a candidate set of lines~\cite{gattermann2017line}, while others generate lines dynamically or use heuristic procedures to construct promising routes. These approaches also differ in whether they explicitly model passenger routing, capacity constraints, and operational restrictions such as transfers and route length limits.

A central computational challenge in line planning is the existence of an exponentially large set of possible routes. Optimizing over this route space requires solving difficult route-generation or pricing subproblems. Column generation and related decomposition methods are natural tools in such settings, since they allow candidate routes to be generated dynamically rather than enumerated in advance~\cite{Barnhart1998Branch}. However, directly applying column generation in general line planning settings is challenging because the pricing problem itself can be computationally difficult; tractable cases often require restrictive assumptions on route structure, such as short lines whose length grows logarithmically with network size~\cite{Borndorfer2007}, and related pricing problems admit hardness-of-approximation results~\cite{Perivier2021}. As a result, scalable approaches often rely on predefined candidate sets, problem-specific relaxations, heuristic pricing procedures, or other modeling simplifications, such as relaxing capacity constraints~\cite{Bertsimas2021}. The single-mode version of our model is closest to passenger-oriented line planning, but differs from prior formulations in its combination of route generation, capacity-constrained passenger assignment, and urban-scale computation. Our formulation maintains explicit passenger flow and capacity constraints while using a column-generation-based approach with heuristic relaxations of the pricing problem to compute line plans at scale. Moreover, our model captures several operationally critical aspects of transit systems, including explicit control over the number of passenger transfers.

Beyond the single-mode line planning problem, there is a growing body of work on integrating multiple transportation services, often with a fixed-route transit component. A common theme in this literature is to use flexible or on-demand mobility as a first- and last-mile complement to public transit, including ridesharing-based integration with existing transit services~\cite{Stiglic2018}. These studies highlight the potential of flexible mobility services to expand the reach of public transit and improve accessibility, especially in settings where fixed-route services alone are inefficient or difficult to operate at high frequency. However, much of this work either assumes that the fixed-route transit network is given or restricts the design space to predefined routes, hubs, road segments, service zones, or aggregate service levels~\cite{Auad-Perez2022, Ng2024}. Existing studies also often emphasize operational or tactical decisions, such as fleet sizing, dispatching, and rebalancing~\cite{Pinto2020}, frequently relying on simulation-based evaluation or heuristic solution methods. In contrast, our model jointly designs end-to-end transit lines, service frequencies, and on-demand legs under a common operating budget, while explicitly accounting for passenger routing and operational constraints.

\subsection{Outline}
The remainder of the paper is organized as follows. \Cref{sec:statement} defines the problem formally and presents our models. \Cref{sec:heuristics} presents our computational algorithms for solving the models for large-scale instances. \Cref{sec:computatation} evaluates our approach on urban-scale instances.

\section{Model}\label{sec:statement}
In this section, we present the problem definition, our mixed-integer programming model and column generation pricing subproblems.  

Let $\graphT=(\nodesT,\edgesT)$ be a directed complete graph representing the transportation network on which \taxi~vehicles may travel. Let $\scaleTaxi\costUV$ denote the operating cost of an \taxi~vehicle traveling on edge $(\nodeU, \nodeV)\in \edgesT$, where $\scaleTaxi$ is a scaling parameter and $\costUV$ denotes the shortest-path distance between $\nodeU$ and $\nodeV$ in the road network. The static passenger demand $\requestST$ denotes the number of passengers seeking to travel from origin $\ori\in \nodesT$ to destination $\dest\in \nodesT$ over a given time span. Let $\requestSet = \left\{(\ori, \dest)\in\nodesT\times\nodesT\mid \requestST > 0\right\}$ denote the set of origin-destination (OD) pairs with positive passenger demand. Let $\graphB = (\nodesB, \edgesB)$, a directed graph, be the \bus \ network where $\nodesB \subset \nodesT$ and $\edgesB \subset \edgesT$.  The route of a bus line $\oneline$ is represented by $\graphline = (\nodesline, \edgesline)$, consisting of copies of nodes and edges of $\graphB$. $\graphline$ is a walk that is the concatenation of a simple path and its reverse. To ease the terminology, we call such a walk a ``bidirectional simple path". \Bus es operate on $\graphline$ by circulating along the loop. Trip options in our multi-modal transportation system are direct \taxi \ trips and trips with a \bus \ leg (with or without complementary \taxi~services).

We assume that the transit authority guarantees passengers a maximum waiting time at each \bus~stop on each line. This  requirement effectively enforces a service headway of at most the maximum waiting time between consecutive~\bus es on the same line. If the time required for a \bus~to complete a full loop exceeds this maximum allowable waiting time, multiple vehicles must operate concurrently on that line during the time span of the waiting time. In our model, demand within a representative headway interval is satisfied by determining the number of vehicles assigned to each line. Consequently, both operating costs and vehicle capacities are normalized to the service provided by each vehicle within a single headway interval rather than over a complete route loop. We define the operating cost of each vehicle to be proportional to the distance traveled. Assuming that each \bus~vehicle travels at the same speed and operates continuously throughout the headway interval, the normalized operating cost is identical across vehicles. We further assume that all vehicles have identical physical capacities. However, during a single headway interval, a vehicle may traverse only a portion of the line if completing a full loop requires more than one headway period. Therefore, we normalize the effective per-vehicle capacity over the entire line by scaling the full vehicle capacity by the fraction of the loop traversed during the headway interval.

Specifically, suppose that the operator promises that the maximum passenger waiting time at a \bus\ stop is $T$ minutes. We assume each \bus\ travels at the same speed, and is able to travel $R$ meters in $T$ minutes. We also assume that each \bus\ vehicle has $\capacity$ physical capacity (for example, 50). For any line $\oneline$, suppose it has total length $M_l$ obtained by summing over the edges of the line. To fulfill the maximum waiting time promise, at least $\lbL = \left\lceil\frac{M_l}{R}\right\rceil$ \bus es are required to ensure that passengers wait no longer than $T$ minutes at the bus stops. We need $\left\lceil\frac{M_l}{R}\right\rceil$ vehicles to collaboratively cover the entire line, as each covers $R$ meters during the time window.
Hence, the number of \bus es operating on line $\oneline$ must be one of the following: 
\[
\freqL = \{0, \lbL, \lbL+1, \lbL+2, \ldots\}.
\]
In a time window of $T$ minutes, the distance each vehicle travels on line $\oneline$ is given by:
\[
\costL = R.
\]
It follows that the per-vehicle cost $\scaleBus\costL$ on line $\oneline$ in this case is $\scaleBus R$. Since each \bus\ has capacity $\capacity$, it provides $\capacity$ seats to passengers along the $R$ meters it travels. We average this capacity over the entire line. For line $\oneline$ with distance $M_l$, each vehicle provides seats
\[
\capacity_{\oneline} = \capacity \cdot \frac{\costL}{M_l} = \capacity \cdot \frac{R}{M_l}
\]
over the line as an averaging approximation that ignores position-dependent loading.

In our setting, we are given a total budget $\budget$ for operating \taxi~and \bus~vehicles, the transportation network $\graphT$, the \bus~network $\graphB$ and the total passenger demand $\requestST$ for OD pairs. We aim to maximize the passenger demand satisfied by our multi-modal transportation system (i.e., maximize ridership) given our available trip options. To achieve this, we need to determine what \bus \ lines to operate and decide on the number of \bus es assigned to each \bus\ line. We begin by formulating the problem as a mixed-integer program over a prespecified candidate set of \bus~lines, as described in \Cref{sec:masterproblem}. This formulation treats the candidate line set as given. To endogenize this set, we subsequently develop a column generation approach that constructs promising \bus~lines from scratch via pricing subproblems, detailed in \Cref{sec:subproblem}.

\subsection{Master Problem: Optimization of Bus Lines and Operating Frequencies}\label{sec:masterproblem}

\begin{figure}[!htb]
\centering
\includegraphics[scale=0.3]{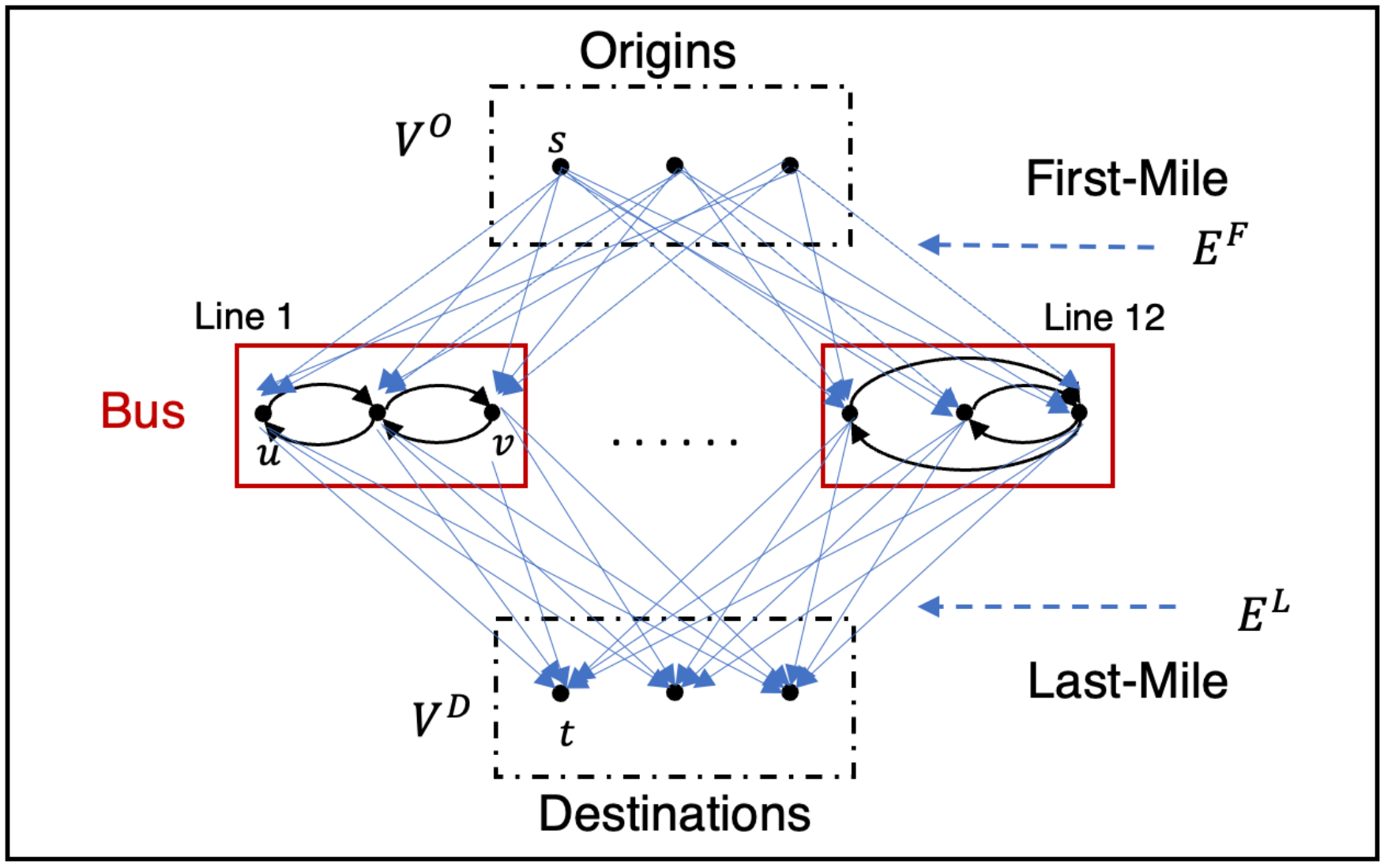}
\caption{An example of $\graphM$ with 12 \bus~lines, modeling passenger flow with a \bus~leg (with or without complementary \taxi~services). $\nodesOri$ denotes the origin nodes. $\nodesDest$ denotes the destination nodes. $\edgesFirst$ is the set of edges of passenger flow using first-mile \taxi~services. $\edgesLast$ is the set of edges of passenger flow using last-mile \taxi~services. In the middle layer (red boxed), we have $\graphline$ for each individual line modeling passenger flow on each line.}
\label{fig:3-layer}
\end{figure}

Given a candidate set of \bus~lines, our MIP formulation models passenger flows from origins to destinations under the available multi-modal trip options. Before presenting the formulation, we first construct the underlying graphs used to represent passenger movements. We begin by describing how flow is tracked for trips that include a \bus~leg. Let $\nodesOri = \left\{\ori\in\nodesT \mid \sum_{\dest} \requestST > 0\right\}$ and $\nodesDest=\left\{\dest \in\nodesT \mid \sum_{\ori} \requestST > 0\right\}$ be the sets of copies of nodes corresponding to origins and destinations. For each line $\oneline \in \lineSet$, let $\graphline = ( \nodesline, \edgesline)$ be the graph consisting of copies of the nodes in $\nodesB$ and edges in $\edgesB$ that belong to the line $\oneline$. We define the edge set $\edgesFirst = \{(\nodeU, \nodeV)\mid \nodeU\in \nodesOri, \nodeV\in \nodesline, \oneline\in \lineSet\}$ to model first-mile \taxi~edge flow connecting origins to \bus \ stops and $\edgesLast = \{(\nodeU, \nodeV)\mid \nodeV\in \nodesDest, \nodeU\in \nodesline, \oneline\in \lineSet\}$ to model last-mile \taxi~edge flow connecting \bus \ stops to destinations.  
We define $\graphM=(\nodesM, \edgesM)$ where $\nodesM=\nodesOri\cup \nodesDest\cup \left(\bigcup_{\oneline\in \lineSet}  \nodesline\right)$ and $\edgesM=\edgesFirst \cup \edgesLast \cup \left(\bigcup_{\oneline\in \lineSet} \edgesline\right)$. The graph models flow with a \bus \ leg with or without complementary first-mile and last-mile \taxi~services. \Cref{fig:3-layer} illustrates an example of $\graphM$ where the middle layer consists of $\graphline$ for each \bus~line $\oneline$. To model passenger flow traveling directly from origins to destinations using only \taxi~services, we use the directed complete graph $\graphT$ where each edge represents such a flow. For vehicles running to serve the passenger flow in $\graphM$, an \taxi~vehicle  traveling between $\nodeU$ and $\nodeV$ for the first- or last-mile segment incurs a cost of $\scaleTaxi\costUV$ when $\nodeU\neq \nodeV$ and zero if $\nodeU = \nodeV$. A \bus~vehicle traveling on $\graphline$ in $\graphM$ to serve the \bus~passenger flow incurs a cost of $\scaleBus\costL$. An \taxi~vehicle traveling on edge $(\nodeU, \nodeV)\in \edgesT$ to serve passenger flow traveling directly from origin $\ori$ to destination $\dest$ incurs a cost of $\scaleTaxi\costUV$.

 In the illustrative example of \Cref{fig:3-layer}, a passenger traveling from origin node $\ori$ to destination node $\dest$ who takes \bus \ line $\oneline_1$ from $\nodeU$ to $\nodeV$ travels edge $(\ori, \nodeU)$, the path from $\nodeU$ to $\nodeV$ in $ E_{\oneline_1}$, and edge $(\nodeV, \dest)$. Note that our model ensures that passengers do not transfer between \bus \ lines, which is desirable for passenger convenience.

We formulate our mixed-integer program (MIP) based on $\graphM$ and $\graphT$ to optimize passenger flow and vehicle assignment. The decision variables are the following. For each $\nodeU\in \nodesT$ and $\nodeV\in \nodesT$, let $\numTaxi_{\nodeU, \nodeV}$ be the total number of single-capacity \taxi~vehicles traveling from $u$ to $v$.  For each $\oneline \in \lineSet$, let $\numBus_{\oneline}$ be the number of \bus \ vehicles traveling on line $\oneline$, which takes value in the feasible set $\freqL = \{0, \lbL, \lbL+1, \lbL+2, \cdots\}$. For each $(\ori, \dest) \in \requestSet$, let $\directTrip_{\ori, \dest}$ be the amount of passenger flow traveling directly from origin $\ori$ to destination $\dest$, which is served by \taxi~vehicles only. For convenience, let $\directTrip_{\nodeU, \nodeV}=0$ for $(\nodeU, \nodeV) \in (\nodesT\times \nodesT)\setminus \requestSet$. For each origin $\ori\in \nodesOri$, each line $\oneline\in \lineSet$ and each line node $\nodeV\in  \nodesline$, let $\getOn^{\ori}_{\oneline, \nodeV}$ be the amount of passenger flow in $\graphM$ originating from $\ori\in \nodesOri$ using first-mile \taxi~edge $(\ori, \nodeV)\in \edgesFirst$ to get on bus line $\oneline$ at $\nodeV$.  For each origin $\ori\in \nodesOri$, each \bus~line $\oneline\in \lineSet$ and each line edge $(\nodeU, \nodeV)\in \edgesline$, let $\busFlow^{\ori}_{\oneline, \nodeU, \nodeV}$ be the amount of passenger flow in $\graphM$ originating from $\ori\in \nodesOri$ using edge $(\nodeU, \nodeV)$ in $ \edgesline$. For each origin $\ori\in \nodesOri$, each line $\oneline\in \lineSet$ and each line node $\nodeV\in \nodesline$, let $\getOff^{\ori}_{\oneline, \nodeV}$ be the  amount of passenger flow in $\graphM$ originating from $\ori\in \nodesOri$ alighting from \bus~line $\oneline$ at $\nodeV$. For each OD pair $(\ori, \dest)\in \requestSet$ and each bus node $\nodeV\in \nodesB$, let $\toDest^{\ori}_{\nodeV, \dest}$ be the amount of passenger flow in $\graphM$ originating from $\ori\in \nodesOri$ alighting from any line at location $\nodeV$ and continuing to destination $\dest$. All decision variables are continuous and non-negative, except for the number of \bus~vehicles on each line, which is a semi-continuous integer decision variable. This modeling choice allows \bus~line decisions to remain discrete while representing passenger flows continuously at the planning stage. The main notations are summarized in \Cref{tab:notations}; additional notation is introduced and defined as needed throughout the paper.

\begin{table}[!h]
\centering
\caption{Main Notations}
\footnotesize
\begin{tabularx}{0.85\textwidth}{cX}
    \hline
    \multicolumn{2}{c}{Network and Demand}\\
    \hline 
    $\graphT = (\nodesT, \edgesT)$ & transportation (\taxi) network\\
    $\graphB = (\nodesB, \edgesB)$ & \bus \ network \\
    $\lineSet$ & set of \bus \ lines\\
    $\graphline = (\nodesline, \edgesline)$ & nodes and edges for line $\oneline$ \\ 
    $\graphM = (\nodesM, \edgesM)$ & network modeling flow using \bus es with or without \taxi~services \\
    $\requestST$ & passenger demand seeking to travel from origin $\ori$ to destination $\dest$\\
    $\requestSet$ & set of OD pairs with positive passenger demand\\
    \hline
    \multicolumn{2}{c}{Input}\\
    \hline
    $\costUV$ & distance of edge $(\nodeU, \nodeV)$ in $\graphT$\\
    $\costL$ & distance traveled by a \bus~vehicle on line $\oneline$\\
    $\capL$ & seats provided by a \bus~vehicle on line $\oneline$\\
    $\freqL$ & set of feasible number of vehicles operating on line $\oneline$\\
    $\scaleTaxi$ & \taxi \ cost scaling factor\\
    $\scaleBus$ & \bus \ cost scaling factor\\
    $\budget$ & total budget \\
    $\detourConst$ & subproblem detour constraint parameter\\
    $\maxLength$ & subproblem maximum length constraint parameter\\
    \hline
    \multicolumn{2}{c}{Decision Variables}\\
    \hline
    $\numTaxi_{\nodeU,\nodeV}$ & number of \taxi \ vehicles traveling from $\nodeU$ to $\nodeV$ \\
    $\numBus_{\oneline}$ & number of \bus \ vehicles operating on line $\oneline$\\
    $\directTrip_{\nodeU,\nodeV}$ & flow from $\nodeU$ to $\nodeV$ served by \taxi~services only \\ 
    $\getOn^{\ori}_{\oneline,\nodeV}$ & flow with origin $\ori$ boarding line $\oneline$ at node $\nodeV$\\
    $\busFlow^{\ori}_{\oneline,\nodeU,\nodeV}$ & flow with origin $\ori$ traveling on line edge $(\nodeU,\nodeV)\in \edgesline$\\
    $\getOff^{\ori}_{\oneline,\nodeV}$ & flow with origin $\ori$ alighting from line $\oneline$ at node $\nodeV$\\
    $\toDest^{\ori}_{\nodeV,\dest}$ & flow with origin $\ori$ and destination $\dest$ alighting from a line at location $\nodeV$
\end{tabularx}
\label{tab:notations}
\end{table}

Given input demand $\requestST$ for OD pairs and input graphs $\graphT$, $\graphB$, $\graphM$ constructed from a moderate-sized input candidate set of \bus~lines $\lineSet$, the MIP formulation $\mip$ is the following. Our formulation models the passenger flow separately for each origin $\ori\in \nodesOri$, while correctly accounts for the number of \taxi~vehicles and \bus~vehicles needed for the whole system.

\begin{align}
\left(\mip\right) &&\max  \sum_{(\ori,\dest)\in \requestSet} \directTrip_{\ori,\dest} + \sum_{(\ori,\dest)\in \requestSet} \sum_{\nodeV\in \nodesB} &\toDest^{\ori}_{\nodeV,\dest}
\label{eqn:Ridership}\\
&&\text{s.t.}\hskip3cm&&\nonumber\\
&&
\directTrip_{\ori,\dest} + \sum_{\nodeV\in \nodesB} \toDest^{\ori}_{\nodeV,\dest} &\leq \requestST  &&\forall (\ori,\dest)\in \requestSet\label{eqn:Request}\\
&& \getOff^{\ori}_{\oneline,\nodeU} + \sum_{(\nodeU,\nodeV)\in  \edgesline} \busFlow^{\ori}_{\oneline,\nodeU,\nodeV} -\getOn^{\ori}_{\oneline,\nodeU} - \sum_{(\nodeV,\nodeU)\in  \edgesline} \busFlow^{\ori}_{\oneline,\nodeV,\nodeU} &= 0    && \forall \ori\in \nodesOri, ~\forall \oneline\in \lineSet, ~\forall \nodeU\in \nodesline \label{eqn:Flow}\\
\textcolor{blue}{[\tilde \busDual_{\ori,\nodeV}]} \ &&
\sum_{\dest\in \nodesDest: (\ori,\dest)\in \requestSet} \toDest^{\ori}_{\nodeV,\dest} - \sum_{\oneline\in \lineSet: \nodeV\in  \nodesline} \getOff^{\ori}_{\oneline,\nodeV} &= 0    &&\forall \ori\in \nodesOri, ~ \forall \nodeV\in \nodesB  \label{eqn:Agg}\\
&&
-  \capL \numBus_{\oneline} + \sum_{\ori\in \nodesOri} \busFlow^{\ori}_{\oneline,\nodeU,\nodeV}    &\leq 0  &&\forall \oneline\in \lineSet, ~ \forall (\nodeU, \nodeV)\in  \edgesline \label{eqn:Capa}
\end{align}\vskip-.5cm\begin{align}
\textcolor{blue}{[\tilde \taxiDual_{\nodeU,\nodeV}]} \ &&
-\numTaxi_{\nodeU,\nodeV} + \directTrip_{\nodeU,\nodeV} + \mathbbm{1}_{\{\nodeU\in \nodesB\}} \sum_{\ori\in \nodesOri: (\ori, \nodeV)\in \requestSet} \toDest^{\ori}_{\nodeU, \nodeV} + \mathbbm{1}_{\{\nodeU\in \nodesOri\}}\sum_{\oneline\in \lineSet: (\nodeU, \nodeV)\in  \edgesFirst} \getOn^{\nodeU}_{\oneline,\nodeV} &\leq 0  && \forall \nodeU\in \nodesT,~\nodeV\in \nodesT  \label{eqn:Taxi}
\end{align}\vskip-.5cm\begin{align}
&&
\sum_{\nodeV\in \nodesT\setminus\{\nodeU\}} \numTaxi_{\nodeV,\nodeU} - \sum_{\nodeV\in \nodesT\setminus\{\nodeU\}} \numTaxi_{\nodeU,\nodeV} &= 0  &&  \forall \nodeU\in \nodesT\label{eqn:TaxBal}\\
\textcolor{blue}{[\tilde \budgetDual]} \ &&\sum_{(\nodeU,\nodeV)\in \edgesT} \scaleTaxi\costUV \numTaxi_{\nodeU,\nodeV} + \sum_{\oneline\in \lineSet} \scaleBus \costL \numBus_{\oneline} &\leq \budget  \label{eqn:Budget}
\end{align}\vskip-.5cm\begin{align}
&& \numTaxi_{\nodeU,\nodeV} \geq 0  \quad &\forall \nodeU\in \nodesT, ~\forall \nodeV\in \nodesT\\
&& \numBus_{\oneline} \in \freqL \quad  &\forall \oneline\in \lineSet\\
&& \directTrip_{\ori,\dest} \geq 0 \quad &\forall (\ori,\dest)\in \requestSet\\
&& \getOn^{\ori}_{\oneline,\nodeV} \geq 0 \quad &\forall \ori\in \nodesOri, ~\forall \oneline\in \lineSet, ~\forall \nodeV\in  \nodesline\\
&& \busFlow^{\ori}_{\oneline,\nodeU,\nodeV} \geq 0   \quad &\forall \ori\in \nodesOri, ~\forall \oneline\in \lineSet, ~\forall (\nodeU,\nodeV)\in  \edgesline\\
&& \getOff^{\ori}_{\oneline,\nodeV} \geq 0 \quad &\forall \ori\in \nodesOri, ~\forall \oneline\in \lineSet, ~\forall \nodeV\in  \nodesline\\
&& \toDest^{\ori}_{\nodeV,\dest}  \geq 0 \quad &\forall (\ori,\dest)\in \requestSet, ~\forall \nodeV\in \nodesB
\end{align}


\Objcref{eqn:Ridership} maximizes total satisfied demand by aggregating passenger flows across \taxi-only trips $\sum_{(\ori,\dest)\in \requestSet} \directTrip_{\ori,\dest}$ and all other available travel options $\sum_{(\ori,\dest)\in \requestSet} \sum_{\nodeV\in  \nodesB} \toDest^{\ori}_{\nodeV,\dest}$. For each OD pair $(\ori, \dest)\in \requestSet$, \constrscref{eqn:Request} require  that the combined demand satisfied by \taxi-only services $\directTrip_{\ori,\dest}$ and all other travel options $\sum_{\nodeV\in  \nodesB} \toDest^{\ori}_{\nodeV,\dest}$ not exceed the total demand. For each origin $\ori\in \nodesOri$, each line $\oneline\in \lineSet$ and each line node $\nodeU\in  \nodesline$, \constrscref{eqn:Flow} enforce flow conservation for passengers originating at each origin $\ori\in \nodesOri$ at every node $\nodeU\in  \nodesline$ in $\graphM$. For each origin $\ori\in \nodesOri$ and each bus node $\nodeV\in \nodesB$, \constrscref{eqn:Agg} require that the passenger flow from origin $\ori\in \nodesOri$ alighting at location $\nodeV$ toward all destinations $\sum_{\dest\in \nodesDest: (\ori,\dest)\in \requestSet} \toDest^{\ori}_{\nodeV,\dest}$  equals the flow from $\ori\in \nodesOri$ arriving at that location $\sum_{\oneline\in \lineSet: \nodeV\in  \nodesline} \getOff^{\ori}_{\oneline,\nodeV}$. For each line $\oneline \in \lineSet$ and each line edge $(\nodeU, \nodeV)\in  \edgesline$,  \constrscref{eqn:Capa} require that the total passenger flow on each edge in $\sum_{\ori\in \nodesOri} \busFlow^{\ori}_{\oneline,\nodeU,\nodeV}$ not exceed the seating capacity provided by \bus~vehicles operating on the line $\capL \numBus_{\oneline}$. For each edge $(\nodeU, \nodeV)\in \edgesT$, \constrscref{eqn:Taxi} require that the combined passenger flow from \taxi-only services $\directTrip_{\nodeU,\nodeV}$, last-mile trips $\mathbbm{1}_{\{\nodeU\in \nodesB\}} \sum_{\ori\in \nodesOri: (\ori, \nodeV)\in \requestSet} \toDest^{\ori}_{\nodeU, \nodeV}$, and first-mile trips $\mathbbm{1}_{\{\nodeU\in \nodesOri\}}\sum_{\oneline\in \lineSet: (\nodeU, \nodeV)\in  \edgesFirst} \getOn^{\nodeU}_{\oneline,\nodeV}$ not exceed the seating capacity of \taxi~vehicles. For each node $\nodeU\in \nodesT$,  \constrcref{eqn:TaxBal} requires \taxi~vehicles to rebalance. Finally, \constrcref{eqn:Budget} requires that the combined operating cost of \taxi~services $\sum_{(\nodeU,\nodeV)\in \edgesT} \scaleTaxi\costUV \numTaxi_{\nodeU,\nodeV}$ and \bus~services $\sum_{\oneline\in \lineSet} \scaleBus \costL \numBus_{\oneline}$ not exceed the available budget $\budget$. We denote by $\tilde \busDual$, $\tilde \taxiDual$, and $\tilde \budgetDual$ the dual variables of the linear programming (LP) relaxation of $\mip$ corresponding to \constrscref{eqn:Agg}, \constrscref{eqn:Taxi}, and \constrcref{eqn:Budget},  respectively. We formulate column generation subproblems for constructing the candidate set $\lineSet$ using these dual variables (\Cref{sec:subproblem}). 

Our MIP formulation $(\mip)$ ensures no \bus-to-\bus~transfers, accounting for passengers' convenience. Passenger flow originating at $\ori\in \nodesOri$ is tracked by the origin-indexed line-specific variables $\busFlow^{\ori}_{\oneline, \nodeU, \nodeV}$, $\getOn^{\ori}_{\oneline, \nodeV}$ and $\getOff^{\ori}_{\oneline,\nodeV}$. A passenger boards line $\oneline$ at stop $\nodeU$ and alights at stop $\nodeV$, exiting the \bus~layer of $\graphM$ entirely via the last-mile edges $\edgesLast$ into destinations $\nodesDest$, as specified by \constrscref{eqn:Agg}. Since no edge in $\edgesM$ connects an alighting node of one line to a boarding node of another, the only path between two \bus~lines passes through $\nodesDest$ and back through $\nodesOri$, which would constitute an entirely separate trip. Therefore, every passenger uses at most one \bus~leg. This structural no-transfer guarantee distinguishes $(\mip)$ from standard multi-commodity flow formulations, where flow may traverse arbitrarily many intermediate nodes.

\Constrcref{eqn:TaxBal} models \taxi~vehicle rebalancing to serve as a proxy for the real-time operational behavior of \taxi~services. \Taxi~operations are highly dependent on and sensitive to demand patterns, which are difficult to model exactly at the planning stage. Note that the formulation remains valid when $\bigcup_{\oneline\in\lineSet} \nodesline \subsetneq \nodesB$. \Constrcref{eqn:Agg} correctly computes the flow from any in-use bus node to the respective destinations, leaving the flow from any unused bus node to zero, i.e. $\toDest^{\ori}_{\nodeV, \dest}=0, ~ \forall \nodeV\in \nodesB\setminus \bigcup_{\oneline\in\lineSet} \nodesline, \forall (\ori, \dest)\in \requestSet$. This modeling choice supports the column generation framework used for constructing the candidate set of lines. Additionally, note that $\graphM$ may contain phantom direct \taxi~passenger flow from origin $\ori$ to destination $\dest$ that is not explicitly captured by the decision variable $\directTrip_{\ori, \dest}$; this flow can be computed as $\sum_{\ori\in \nodesOri}\sum_{\oneline\in \lineSet}\sum_{\nodeV\in  \nodesline} \min \left\{\getOn^{\ori}_{\oneline,\nodeV}, \getOff^{\ori}_{\oneline,\nodeV}\right\}$. It represents the passenger flow in $\graphM$ that originates from the same origin, boarding and alighting at the same \bus~stop. Therefore, the exact total demand satisfied by \taxi~services is given by $\sum_{(\ori, \dest)\in \requestSet}\directTrip_{\ori, \dest} + \sum_{\ori\in \nodesOri}\sum_{\oneline\in \lineSet}\sum_{\nodeV\in  \nodesline} \min \left\{\getOn^{\ori}_{\oneline,\nodeV}, \getOff^{\ori}_{\oneline,\nodeV}\right\}$ as explained. Note that the phantom routing is weakly dominated by the direct routing due to the triangle inequality of the cost structure. Therefore, the optimal LP and MIP objective values of $(\mip)$ are not affected by possible phantom routing. Phantom trips can only appear in an optimal solution of $(\mip)$ when the same boarding and alighting stop is on a shortest path.

\subsection{Pricing Problem: Construction of Candidate Bus Lines}\label{sec:subproblem}
The MIP maximizes ridership by optimizing over a candidate set of \bus~lines. However, explicitly enumerating all feasible routes is computationally intractable due to their exponential number. To address this, we employ a heuristic column generation scheme~\cite{Barnhart1998Branch} that dynamically constructs promising lines via the pricing problem as needed. 

In our formulation, a \bus~line is considered feasible if it follows a realistic path, defined as a directed walk formed by the concatenation of a simple path and its reverse. 
We ensure these lines are operationally efficient and acceptable to riders by requiring that the end-to-end length in either direction not exceed a $\detourConst$ multiple of the shortest-path distance between the two endpoints in $\graphT$, where edge distances in $\graphT$ are defined as shortest-path distances in the underlying road network (as specified in the experimental design). We additionally enforce a maximum length constraint, requiring that the total length of the bidirectional path not exceed $\maxLength$. The subproblem is formulated as a mixed-integer program that selects the sequence of edges to construct such high-value routes, where $\detourConst$ and $\maxLength$ are input parameters.

We present two novel pricing problem formulations. The first pricing problem, detailed below, approximates the exact reduced-cost calculation for the generated line. We also design a second pricing problem in \Cref{subsec:CG2} with reduced complexity to complement the first subproblem, as it may offer faster solving times. Together, we use both subproblems in our algorithm to balance computational complexity and performance. Recall that in \Cref{sec:statement} we introduced the \bus~network $\graphB=(\nodesB,\edgesB)$ as a directed subgraph of $\graphT$. We generate candidate lines from $\graphB$ under the following assumption.

\begin{Assumption}\label{assump:connectivity}
$\graphB=(\nodesB,\edgesB)$ is strongly connected. Moreover, for each $(u,v)\in \edgesB$, we also have $(v,u) \in \edgesB$. 
\end{Assumption}

Given $\graphB=(\nodesB,\edgesB)$, we define an augmented graph $\graphB'=(\nodesB',\edgesB')$ for the generation of candidate lines. Specifically, let $\nodesB' = \nodesB \cup \{\source, \sink\}$ and $\edgesB' = \edgesB \cup \{(\source,\nodeU), (\nodeU, \sink)\mid u\in \nodesB\}$. That is, $\graphB'$ is obtained by adding a source node $\source$ and a sink node $\sink$ to $\graphB$, with edges connecting the source to all nodes and all nodes to the sink. $\graphB'$ serves as a helper graph that enables the subproblem to generate a simple path from $\source$ to $\sink$, which corresponds to selecting a simple path in $\graphB$.

\subsubsection{Pricing Problem Formulation I}\label{subsec:CG1}
The formulation uses the dual variables associated with \constrscref{eqn:Agg}, \constrscref{eqn:Taxi}, and \constrcref{eqn:Budget} in the LP relaxation of $\mip$, denoted by $\tilde \busDual$, $\tilde \taxiDual$, and $\tilde \budgetDual$, respectively. The decision variables are defined as follows. For each edge $(\nodeU, \nodeV)\in \edgesB'$, let $\subEdge_{\nodeU, \nodeV}$ be a binary variable indicating whether edge $(\nodeU, \nodeV)$ lies on the directed simple path selected by the subproblem. The line generated is the bidirectional simple path formed by this simple path together with its reverse. Accordingly, each edge pair of the line is indicated exactly once, and the terms $(\costUV + c_{\nodeV, \nodeU})\subEdge_{\nodeU, \nodeV}$ appearing in the objective and in the length constraints account for both directions of travel along the line. For each origin $\ori \in \nodesOri$ and node $\nodeU \in \nodesB$, let $\getOn^{\ori}_{\nodeU}$ denote the flow originating from $\ori$ that boards the new line at $\nodeU$, and let $\getOff^{\ori}_{\nodeU}$ denote the flow originating from $\ori$ that alights from the new line at $\nodeU$.
For each origin $\ori \in \nodesOri$ and edge $(\nodeU, \nodeV) \in \edgesB$, let $\busFlow^{\ori}_{\nodeU,\nodeV}$ denote the flow originating from $\ori$ that traverses edge $(\nodeU, \nodeV)$ on the new line.

\begin{equation}
\left(\pricing\right) \qquad  \max -\frac{\scaleBus}{\capacity}\sum_{(\nodeU,\nodeV)\in \edgesB} \tilde\budgetDual\left(\costUV+ c_{v,u}\right)\subEdge_{\nodeU,\nodeV} -\sum_{\ori\in \nodesOri} \sum_{\nodeU\in \nodesB} \tilde \taxiDual_{\ori,\nodeU}\getOn^{\ori}_{\nodeU} + \sum_{\ori\in \nodesOri}\sum_{\nodeU\in \nodesB}  \tilde \busDual_{\ori,\nodeU}\getOff^{\ori}_{\nodeU}\label{eqn:Sub1-maxride} 
\end{equation}
\begin{align}
&&s.t.\hskip3cm&&\\[.3cm]
&& \sum_{(\source,\nodeU)\in \edgesB'} \subEdge_{\source,\nodeU} & = 1 \label{eqn:Sub1-source}\\
&& \sum_{(\nodeU,\sink)\in \edgesB'} \subEdge_{\nodeU,\sink} &= 1 \label{eqn:Sub1-sink}\\
&& \sum_{(\nodeU,\nodeV)\in \edgesB'}\subEdge_{\nodeU,\nodeV} - \sum_{(\nodeV,\nodeU)\in \edgesB'}\subEdge_{\nodeV,\nodeU} & = 0 && \quad \forall \nodeU\in \nodesB\label{eqn:Sub1-conservation}\\
&& \sum_{(\nodeU,\nodeV)\in \edgesB: \nodeU\in F, \nodeV\in F} \subEdge_{\nodeU,\nodeV} &\leq |F| - 1 && \quad \forall F\subsetneq \nodesB, ~|F|\geq 2\label{eqn:Sub1-subtour} \\
&& \sum_{(\nodeU, \nodeV)\in \edgesB}\subEdge_{\nodeU, \nodeV}\costUV &\leq \detourConst\sum_{\nodeU, \nodeV\in \nodesB}\costUV\subEdge_{\source, \nodeU}\subEdge_{\nodeV, \sink}  \label{eqn:Sub1-detour1}\\
&& \sum_{(\nodeU, \nodeV)\in \edgesB}\subEdge_{\nodeU, \nodeV}c_{\nodeV, \nodeU} &\leq \detourConst\sum_{\nodeU, \nodeV\in \nodesB}c_{\nodeV, \nodeU}\subEdge_{\source, \nodeU}\subEdge_{\nodeV, \sink}  \label{eqn:Sub1-detour2}\\
&& \sum_{(\nodeU, \nodeV)\in \edgesB} (\costUV+c_{\nodeV, \nodeU})\subEdge_{\nodeU, \nodeV} &\leq \maxLength \label{eqn:Sub1-limit} \\
&& \sum_{\ori\in \nodesOri}\busFlow^{\ori}_{\nodeU,\nodeV} &\leq \subEdge_{\nodeU,\nodeV} + \subEdge_{\nodeV,\nodeU} && \quad  \forall (\nodeU,\nodeV)\in \edgesB \label{eqn:Sub1-cap}\\
&& \getOn^{\ori}_{\nodeU} + \sum_{(\nodeV, \nodeU)\in \edgesB}\busFlow^{\ori}_{\nodeV,\nodeU} &= \getOff^{\ori}_{\nodeU}+\sum_{(\nodeU, \nodeV)\in\edgesB}\busFlow^{\ori}_{\nodeU,\nodeV} &&\quad  \forall \ori \in \nodesOri, ~\forall \nodeU\in \nodesB\label{eqn:Sub1-flow}\\
&& \subEdge_{\nodeU,\nodeV} & \in \{0,1\}  &&\quad  \forall (\nodeU,\nodeV)\in \edgesB' \\
&& \getOn^{\ori}_{\nodeU} & \geq 0 && \quad  \forall \ori\in \nodesOri, ~\forall \nodeU\in \nodesB\\
&& \getOff^{\ori}_{\nodeU} &\geq 0 && \quad  \forall \ori\in \nodesOri, ~\forall \nodeU\in \nodesB\\
&& \busFlow^{\ori}_{\nodeU,\nodeV} & \geq 0 && \quad  \forall \ori\in \nodesOri, ~\forall (\nodeU,\nodeV)\in \edgesB
\end{align}

\Objcref{eqn:Sub1-maxride} maximizes the estimated per-capacity reduced cost of the line generated by the subproblem to be added to the candidate set. The estimated per-capacity reduced cost is computed as the marginal benefit of creating the bidirectional simple path minus its construction cost, by using the corresponding dual variables of the LP relaxation of $(\mip)$. The computation credits the line's boarding and alighting flow against their respective master duals, with opposite signs as detailed below. For any variable with a zero coefficient in the master problem's original objective, its contribution to the reduced-cost objective is $-\sum_i a_{ij}\tilde\pi_i$, summed over the constraints $i$ in which it appears, where $a_{ij}$ is its coefficient in constraint $i$ and $\tilde\pi_i$ is the corresponding dual value. Since $\getOn^{\ori}_{\nodeU}$ enters \eqref{eqn:Taxi} with coefficient $+1$, it contributes $-\tilde\taxiDual\getOn$; since $\getOff^{\ori}_{\nodeU}$ enters \eqref{eqn:Agg} with coefficient $-1$, it contributes $+\tilde\busDual\getOff$.  We emphasize what is and is not exact in this estimate. The exact reduced cost of a candidate line variable $\numBus_{\oneline}$ depends on the duals of the line-specific capacity \constrscref{eqn:Capa}, which do not exist until the line has been added to the restricted master problem. \Objcref{eqn:Sub1-maxride} therefore prices the line using only duals that are available in advance---the alighting-balance duals $\tilde \busDual$, the \taxi~capacity duals $\tilde \taxiDual$, and the budget dual $\tilde \budgetDual$---and credits the boarding and alighting flow that the line would carry. A positive objective value of $(\pricing)$, even when optimized over the pricing formulation, therefore indicates a promising rather than a certified improving column, which is why we describe the overall scheme as a heuristic column generation procedure. Note that the dual variable $\tilde \taxiDual_{\nodeU,\nodeV}$ is defined for each $(\nodeU,\nodeV)\in \edgesT$. Since we define $\edgesT = \{(\nodeU, \nodeV): \nodeU, \nodeV\in \nodesT, \nodeU\neq \nodeV\}$ and the origin nodes $\nodesOri$ and \bus~nodes $\nodesB$ are subsets of $\nodesT$, the term $\tilde \taxiDual_{\ori,\nodeU}$ appearing in \objcref{eqn:Sub1-maxride} is therefore well-defined whenever $\ori\neq \nodeU$; for convenience, let $\tilde \taxiDual_{\nodeU,\nodeU} = 0$ for $\nodeU\in \nodesT$.  \Constrscref{eqn:Sub1-source}, \constrscref{eqn:Sub1-sink} and \constrscref{eqn:Sub1-conservation} form a path from the source $\source$ to the sink $\sink$ in $\graphB'$. \Constrscref{eqn:Sub1-subtour} eliminate cycles in either direction of the new line. \Constrcref{eqn:Sub1-detour1} and \constrcref{eqn:Sub1-detour2} ensure that the length of either direction of the new line does not exceed $\detourConst$ times the distance between its endpoints in $\graphT$. \Constrcref{eqn:Sub1-limit} restricts the total length of the bidirectional simple path to be at most $\maxLength$.
\Constrscref{eqn:Sub1-cap} ensure that passenger flow occurs only on edges included in the new line, and \constrscref{eqn:Sub1-flow} enforce flow conservation at each node along the line. We note that \constrscref{eqn:Sub1-source}--\eqref{eqn:Sub1-conservation} alone admit a degenerate solution in which only the connector edges $(\source,u)$ and $(u,\sink)$ are selected for a single node $u$, with no edge of $\edgesB$ chosen, corresponding to a single-node ``line'' with $M_l=0$. Our implementation excludes this case by requiring at least one edge of $\edgesB$ to be selected, in $(\pricing)$ here and $(\pricingPrime)$ below as well as in their cost-minimization counterparts introduced in \Cref{subsec:experiments-budget}.

Note that \constrscref{eqn:Sub1-subtour} are difficult to impose explicitly, as their number grows exponentially with the number of nodes. In addition, \constrscref{eqn:Sub1-detour1} and \constrscref{eqn:Sub1-detour2} are endpoint-dependent bilinear constraints. We therefore enforce these constraints dynamically using callback functions provided by optimization solvers (such as Gurobi). Specifically, we solve $(\pricing)$ initially without constraints~\labelcref{eqn:Sub1-subtour,eqn:Sub1-detour1,eqn:Sub1-detour2}. At each integer solution encountered during the branch-and-bound process, the callback extracts the selected edges and checks whether they form a single simple path or contain disconnected cycles. If a subtour violation is detected, we identify the corresponding node subset $F\subsetneq \nodesB$ and add the violated subtour-elimination constraint. For the detour constraints, once the selected edges define a bidirectional simple path, the callback identifies the endpoints of the path and checks whether either direction violates the prescribed detour bound. If a violation is detected, we add the corresponding endpoint-specific linearized detour constraint. This avoids explicitly imposing the bilinear  constraints over all endpoint pairs while ensuring that any integer solution returned by the solver satisfies both the subtour-elimination and detour requirements. In the following, constraints~\labelcref{eqn:Sub2-subtour,eqn:Sub2-detour1,eqn:Sub2-detour2} are handled in the same manner, and we omit the details for brevity.

\subsubsection{Pricing Problem Formulation II}\label{subsec:CG2}
Our second column generation formulation is based on the LP relaxation of an MIP modified from ($\mip$). The modification replaces the capacity \constrscref{eqn:Capa} in  ($\mip$) with the aggregated \constrscref{eqn:AggCapa} below, which remove line-specific dependence, thereby relaxing the original formulation ($\mip$). The dual variables to \constrscref{eqn:AggCapa} are denoted by $\tilde \capDual_{\nodeU,\nodeV}$. Keeping all other parts the same as ($\mip$), let the new MIP obtained by substituting \constrscref{eqn:Capa} with \constrscref{eqn:AggCapa} be denoted by ($\mipPrime$).

\begin{align}
\text{[Bus Capacity]:} -\capL \numBus_{\oneline} + \sum_{\ori\in \nodesOri} \busFlow^{\ori}_{\oneline,\nodeU,\nodeV}  &\leq 0  \quad \forall \oneline\in \lineSet, ~ \forall (\nodeU, \nodeV)\in  \edgesline\nonumber\\
\quad \quad &\Downarrow \nonumber\\
\textcolor{blue}{[\tilde \capDual_{\nodeU,\nodeV}]}\text{[Aggregated Bus Capacity]:}
- \sum_{\oneline\in \lineSet: (\nodeU, \nodeV)\in \edgesline}\capL \numBus_{\oneline} &+ \sum_{\oneline\in \lineSet: (\nodeU, \nodeV)\in \edgesline} \sum_{\ori\in \nodesOri} \busFlow^{\ori}_{\oneline,\nodeU,\nodeV}    \leq 0  \quad \forall (\nodeU,\nodeV)\in \edgesB \label{eqn:AggCapa}
\end{align}

Our second column generation formulation ($\pricingPrime$) uses the dual variables associated with \constrcref{eqn:Budget} and \constrscref{eqn:AggCapa} in the LP relaxation of $(\mipPrime)$, denoted by $\tilde \budgetDual$ and $\tilde \capDual_{\nodeU,\nodeV}$, respectively. $(\pricingPrime)$ prices lines against the aggregated capacity duals $\tilde \capDual$ of $(\mipPrime)$, which are defined on $\edgesB$ independently of the current line set. The objective is thus a close approximation of the per-capacity reduced cost of the line generated by the subproblem for the LP relaxation of $(\mipPrime)$. As in $(\pricing)$, for each edge $(\nodeU, \nodeV)\in \edgesB'$, let $\subEdge_{\nodeU, \nodeV}$ be a binary variable indicating whether edge $(\nodeU, \nodeV)$ lies on the directed simple path selected by the subproblem, with the generated line given by this path together with its reverse.

\begin{equation}
(\pricingPrime)  \qquad \max -\frac{\scaleBus}{\capacity}\sum_{(\nodeU,\nodeV)\in \edgesB}\tilde \budgetDual\left( \costUV+  c_{\nodeV, \nodeU}\right)\subEdge_{\nodeU, \nodeV} + \sum_{(\nodeU,\nodeV)\in \edgesB} (\tilde \capDual_{\nodeU,\nodeV} +\tilde \capDual_{\nodeV,\nodeU})\subEdge_{\nodeU, \nodeV}\label{eqn:Sub2-maxride}
\end{equation}
\begin{align}
&&s.t.\hskip3cm&&\\[.3cm]
&& \sum_{(\source,\nodeU)\in \edgesB'} \subEdge_{\source,\nodeU} & = 1 \label{eqn:Sub2-source}\\
&& \sum_{(\nodeU,\sink)\in \edgesB'} \subEdge_{\nodeU,\sink} &= 1 \label{eqn:Sub2-sink}\\
&& \sum_{(\nodeU,\nodeV)\in \edgesB'}\subEdge_{\nodeU,\nodeV} - \sum_{(\nodeV,\nodeU)\in \edgesB'}\subEdge_{\nodeV,\nodeU} & = 0 && \quad \forall \nodeU\in \nodesB\label{eqn:Sub2-conservation}\\
&& \sum_{(\nodeU,\nodeV)\in \edgesB: \nodeU\in F, \nodeV\in F} \subEdge_{\nodeU,\nodeV} &\leq |F| - 1  && \quad \forall F\subsetneq \nodesB, ~|F|\geq 2\label{eqn:Sub2-subtour}\\
&& \sum_{(\nodeU, \nodeV)\in \edgesB}\subEdge_{\nodeU, \nodeV}\costUV &\leq \detourConst\sum_{\nodeU, \nodeV\in \nodesB}\costUV\subEdge_{\source, \nodeU}\subEdge_{\nodeV, \sink}  \label{eqn:Sub2-detour1}\\
&& \sum_{(\nodeU, \nodeV)\in \edgesB}\subEdge_{\nodeU, \nodeV}c_{\nodeV, \nodeU} &\leq \detourConst\sum_{\nodeU, \nodeV\in \nodesB}c_{\nodeV, \nodeU}\subEdge_{\source, \nodeU}\subEdge_{\nodeV, \sink}  \label{eqn:Sub2-detour2}\\
&& \sum_{(\nodeU, \nodeV)\in \edgesB} (\costUV+c_{\nodeV, \nodeU})\subEdge_{\nodeU, \nodeV} &\leq \maxLength  \label{eqn:Sub2-limit}\\
&&\subEdge_{\nodeU,\nodeV} &\in \{0,1\} && \quad  \forall (\nodeU,\nodeV)\in \edgesB'
\end{align}

\Objcref{eqn:Sub2-maxride} maximizes the marginal benefit of constructing the bidirectional simple path, as captured by the relaxed capacity \constrscref{eqn:AggCapa}, minus its construction cost.
\Constrscref{eqn:Sub2-source}, \constrscref{eqn:Sub2-sink} and \constrscref{eqn:Sub2-conservation} form a path from the source $\source$ to the sink $\sink$ in $\graphB'$. \Constrscref{eqn:Sub2-subtour} eliminate cycles in either direction of the new line. \Constrscref{eqn:Sub2-detour1} and \constrscref{eqn:Sub2-detour2} ensure that the length of either direction of the new line does not exceed $\detourConst$ times the distance between its endpoints in $\graphT$. \Constrcref{eqn:Sub2-limit} restricts the total length of the bidirectional simple path to be at most $\maxLength$.

\textbf{Restricted Master Problem Update} In either subproblem formulation, the solution returned by the program identifies a new \bus~line $\oneline$. This line is added to the candidate set $\lineSet$ in the restricted master linear program, along with the associated variables and constraints. In addition to introducing the primary variable $\numBus_{\oneline}$, we concurrently create the following line-dependent primal variables: $\getOn^{\ori}_{\oneline,\nodeV}$, $\busFlow^{\ori}_{\oneline,\nodeU,\nodeV}$, and $\getOff^{\ori}_{\oneline,\nodeV}$.
We then add new constraints or update existing line-dependent constraints in $(\mip)$, including \constrscref{eqn:Flow}, \constrscref{eqn:Agg}, \constrscref{eqn:Capa}, \constrscref{eqn:Taxi}, and \constrcref{eqn:Budget}. For $(\mipPrime)$, we add or update \constrscref{eqn:Flow}, \constrscref{eqn:Agg}, \constrscref{eqn:AggCapa}, \constrscref{eqn:Taxi}, and \constrcref{eqn:Budget}.
After incorporating the new line, we re-solve the restricted master linear program $(\mip)$ or $(\mipPrime)$ to obtain updated dual values for the next iteration of the column generation subproblems.
Building on this iterative framework, we next address the computational challenges that arise in large-scale settings.

\section{Algorithms and Heuristics}\label{sec:heuristics}
Large-scale instances involving integer decision variables in ($\mip$), as well as in ($\pricing$) and ($\pricingPrime$), make these models difficult to solve. Therefore, we design linear programming (LP)-based heuristics to accelerate the computation of high-quality solutions. \Cref{sec:master-heur} presents heuristic methods for solving the MIP, and \Cref{sec:subproblem-heur} describes heuristics for the column generation subproblems. The overall algorithm is presented in \Cref{alg:full-alg}.

\subsection{Heuristic Solution of the Master MIP} \label{sec:master-heur}
Given the set of OD pairs $\requestSet$, the candidate set of lines $\lineSet$, and the network $\graphT = (\nodesT, \edgesT)$, the number of variables in $(\mip)$ is $O\left(|\nodesT|^3|\lineSet| + |\requestSet||\nodesT|\right)$, and the number of constraints is $O\left(|\nodesT|^2|\lineSet| + |\requestSet|\right)$. The size of $\lineSet$ can therefore substantially impact the computational effort required to solve $(\mip)$. The computational barrier is further driven by the presence of semi-continuous integer variables $\numBus_{\oneline}$, each with lower bound $\lbL$ for $\oneline \in \lineSet$. These variables constitute a major source of difficulty. To address these challenges, \Cref{alg:mip-heur} leverages the LP relaxation to first identify a subset of promising lines and then solves $(\mip)$ restricted to this selected set. The selected set prioritizes lines with positive $\numBus_{\oneline}$ values in the LP based on the ratio $\numBus_{\oneline}/\lbL$, with larger ratios receiving higher priority.

\begin{center}
\begin{minipage}{0.85\linewidth}
\begin{algorithm}[H]
\caption{Heuristic for Solving ($\mip$) in Large-Scale Instances}\label{alg:mip-heur}
\begin{algorithmic}[1]
\Require Candidate set of lines $\lineSet$, budget $\budget$, parameter $m$
\State Solve the LP relaxation of $(\mip)$ with candidate set $\lineSet$ and budget $\budget$
\State Obtain the optimal solution $\tilde \numBus_{\oneline}$ and reduced cost $r(\tilde \numBus_{\oneline})$ for all $\oneline \in \lineSet$
\State Sort $\tilde \numBus_{\oneline}/\lbL$ for all $\oneline \in \lineSet$ with $\tilde \numBus_{\oneline} > 0$
\State Sort $r(\tilde \numBus_{\oneline}) \cdot \lbL$ for all $\oneline \in \lineSet$ with $\tilde \numBus_{\oneline} = 0$
\State Select a total of $m$ lines with the largest values of $\tilde \numBus_{\oneline}/\lbL$, followed by those with the largest values of $r(\tilde \numBus_{\oneline}) \cdot \lbL$
\State Solve $(\mip)$ using the $m$ pre-selected candidate lines and budget $\budget$\\
\Return Solution of $(\mip)$ with the $m$ pre-selected candidate lines and budget $\budget$
\end{algorithmic}
\end{algorithm} 
\end{minipage}
\end{center}

\subsection{Heuristic Solution of the Pricing Problem}\label{sec:subproblem-heur}

Given the \taxi~network $\graphT = (\nodesT, \edgesT)$ and the \bus~network $\graphB = (\nodesB, \edgesB)$, the column generation subproblem $(\pricing)$ has $O(|\nodesT||\nodesB|^2)$ variables and $O(|\nodesB|^2 + |\nodesT||\nodesB|)$ constraints, while $(\pricingPrime)$ has $O(|\nodesB|^2)$ variables and $O(|\nodesB|^2)$ constraints. To address the computational burden for large graphs $\graphB$, we design \Cref{alg:subprob-heur} to restrict each column generation iteration to a reduced subgraph, where \bus~edges are selected based on dual information from the LP relaxation. The algorithm is described for $(\pricing)$ and applies analogously to $(\pricingPrime)$. 

\begin{center}
\begin{minipage}{0.85\linewidth}
\begin{algorithm}[H]
\caption{Heuristic for Solving ($\pricing$)}\label{alg:subprob-heur}
\begin{algorithmic}[1]
\Require{ Dual values of the LP relaxation of $(\mip)$, parameter $p$}
\State Solve LP relaxation of $(\pricing)$ with the dual values from the LP relaxation of $(\mip)$
\State Sort the reduced costs of decision variables $\subEdge_{\nodeU,\nodeV}$ for all $(\nodeU,\nodeV)\in \edgesB$
\State Select the $100p\%$ edges with the highest reduced costs of $\subEdge_{\nodeU,\nodeV}$
\State Solve $(\pricing)$ with only the edge set trimmed to the $p$ fraction\\
\Return Solution of $(\pricing)$ obtained on the trimmed edge set
\end{algorithmic}
\end{algorithm} 
\end{minipage}
\end{center}

\subsection{Full Algorithm}\label{sec:full-alg}
We present the full algorithm for solving the multi-modal transit network planning problem. We employ distinct budget values—$\budget_{LP}$ for solving the LP relaxation of the master problem during the line generation process and $\budget_{MIP}$ for obtaining mixed-integer solutions—to optimize performance.

\algnewcommand{\OR}{\textbf{ or }} 
\algnewcommand{\AND}{\textbf{ and }} 
\begin{center}
\begin{minipage}{0.85\linewidth}
\begin{algorithm}[H]
	\caption{Full Algorithm of Solving the Multi-Modal Transit Network Planning Ridership Maximization Problem}\label{alg:full-alg}
	\begin{algorithmic}[1]
		\Require Initial Candidate Set of Lines $\lineSet$, budget values $\budget_{LP}$ and $\budget_{MIP}$, maximum number of lines $\mathcal{N}$  
		\Repeat
		\State Solve LP relaxation of $(\mipPrime)$ with candidate set of lines $\lineSet$ and budget $\budget_{LP}$, get dual solution
		\State Solve $(\pricingPrime)$ using the pricing problem heuristic
		\State $j_{\pricingPrime}\gets$ objective value of the returned solution of $(\pricingPrime)$ 
		\State $l_{\pricingPrime}\gets$ line(s) extracted from the returned solution of $(\pricingPrime)$ 
		\If{$j_{\pricingPrime} > 0$}
			\State $\lineSet\gets \lineSet\cup\{l_{\pricingPrime}\}$
		\EndIf
		\If{$|\lineSet|<\mathcal{N}$}
			\State Solve LP relaxation of $(\mip)$ with candidate set of lines $\lineSet$ and budget $\budget_{LP}$, get dual solution
			\State Solve $(\pricing)$ using the pricing problem heuristic
			\State $j_{\pricing}\gets$ objective value of the returned solution of $(\pricing)$ 
			\State $l_{\pricing}\gets$ line(s) extracted from the returned solution of $(\pricing)$
			\If{$j_{\pricing} > 0$ }
				\State $\lineSet\gets \lineSet\cup\{l_{\pricing}\}$
			\EndIf
		\EndIf
		\Until{($|\lineSet|=\mathcal{N}$) \OR ($j_{\pricing} \le 0$ \AND $j_{\pricingPrime} \le 0$ )}
		\State Solve $(\mip)$ using the MIP master problem heuristic with candidate set of lines $\lineSet$ and budget $\budget_{MIP}$\\
		\Return Solution of $(\mip)$ obtained using the MIP master problem heuristic with candidate set of lines $\lineSet$ and budget $\budget_{MIP}$
	\end{algorithmic}
\end{algorithm}
\end{minipage}
\end{center}

\section{Computational Results} \label{sec:computatation}
We evaluate the proposed multi-modal transit network design framework through computational experiments. We start by describing the preprocessing of demand and city networks, and then assess the performance of our algorithms in terms of ridership maximization, comparison with the benchmark method and sensitivity analyses. The code is available at: \url{https://github.com/ning74/multimodal_design}.

\subsection{Data Preprocessing}\label{sec:data}

We detail our process for obtaining the demand dataset, the transportation network $\graphT=(\nodesT, \edgesT)$ and the \bus~network $\graphB=(\nodesB, \edgesB)$. We extract city maps and construct travel demand instances for three major cities across the United States: Boston, Chicago, and Atlanta. The city center coordinates are (42.3601, -71.0589), (41.8819, -87.6301), and (33.7681, -84.3806) for Boston, Chicago, and Atlanta, respectively. Our focus is on dense downtown regions to design high-frequency backbone \bus~systems in the multi-modal setting. To generate the passenger trips, the LEHD Origin-Destination Employment Statistics (LODES) data released by the United States Census Bureau \cite{LODES} serves as our foundational dataset. The dataset provides employees' home-work location information at the census block level. We interpret it as representing travel demand during the morning peak period. The process of creating passenger demand in the form of pickup-drop-off locations (in the latitude-longitude format) involves randomly selecting home-work locations within each census block. We then select the generated pickup and drop-off locations within a 3-mile (4,800-meter) radius of the city centers, measured by the geodesic distance. The road networks for the cities are obtained from OpenStreetMap with the distance type parameter set to \texttt{bbox}, which builds a bounding box extending 3 miles north, south, east, and west in approximate geodesic distance. OpenStreetMap \cite{openstreetmap} provides structured geospatial data that simplifies network operations. Given the interdependence between network and demand data processing, we describe the steps together below rather than in separate sections. Boston serves as an illustrative example for our procedure. We apply the same approach to Chicago and Atlanta. 

\begin{figure}[!htb]   
    \begin{subfigure}[b]{0.33\textwidth}
        \textbf{\large Boston}\\[2pt]
        \includegraphics[width=\textwidth]{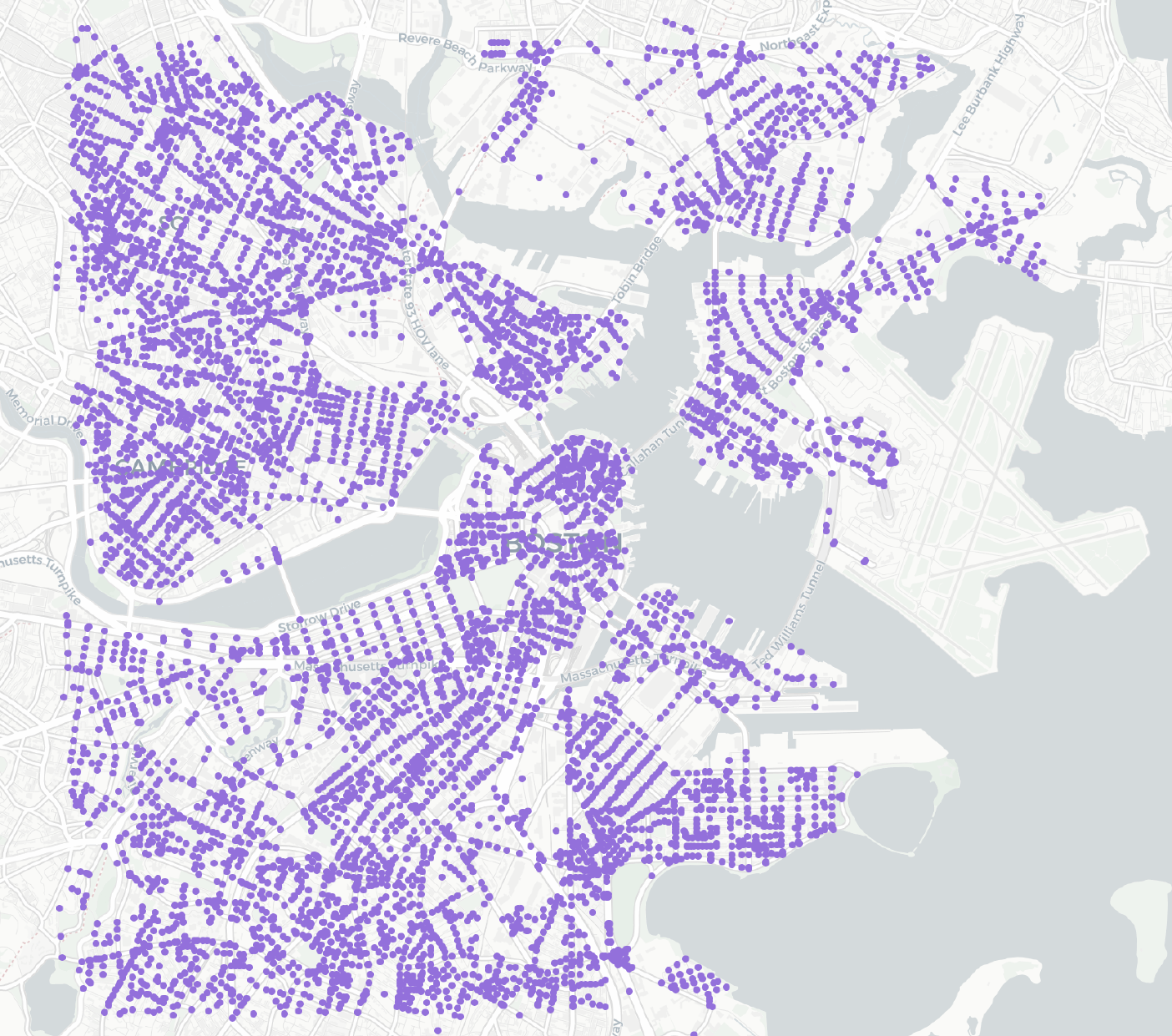}
        \subcaption{OpenStreetMap nodes}
    \end{subfigure}\hfill
    \begin{subfigure}[b]{0.33\textwidth}
        
        \includegraphics[width=\textwidth]{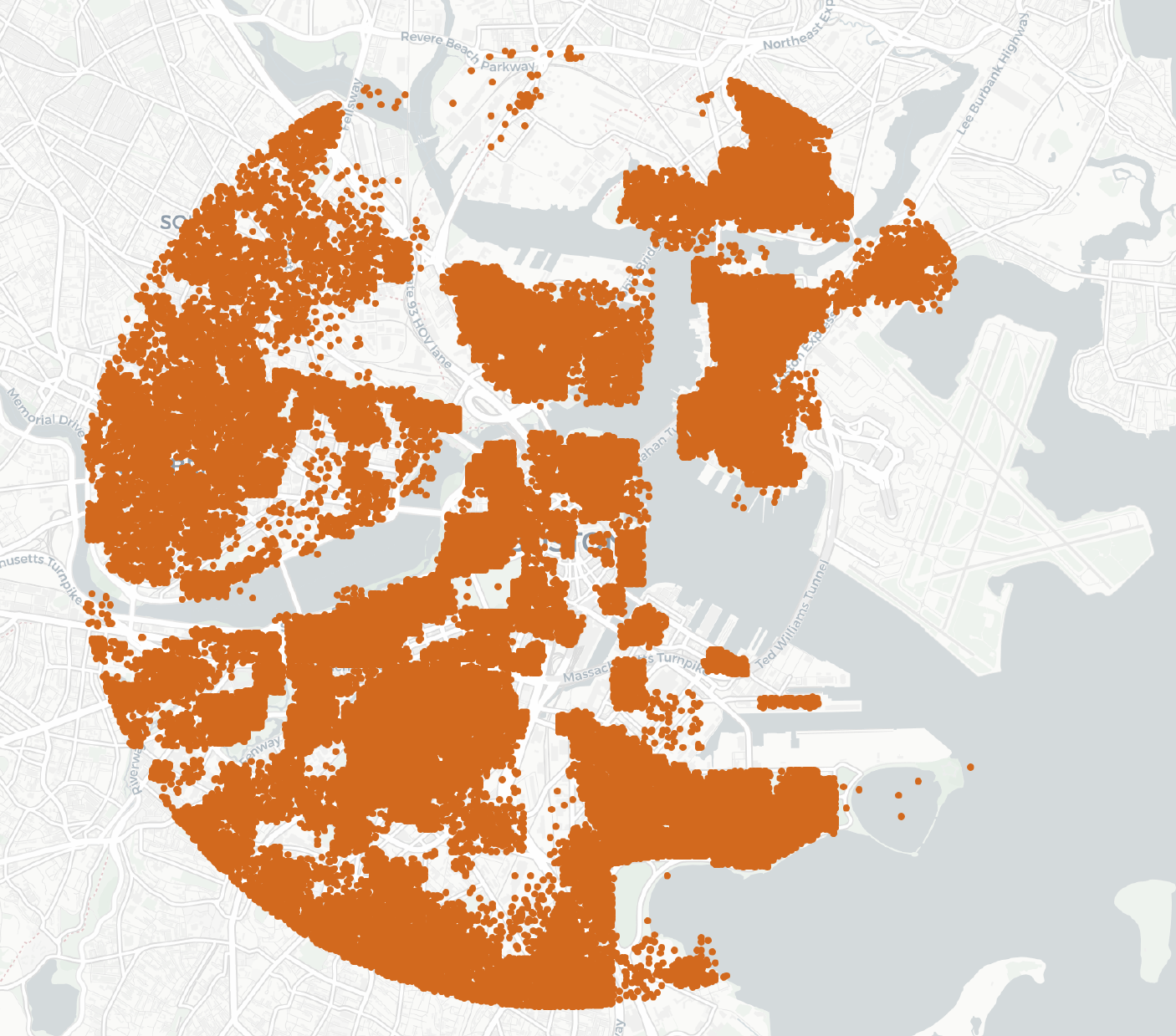}
        \subcaption{Pickup locations}
    \end{subfigure}\hfill
    \begin{subfigure}[b]{0.33\textwidth}
        
        \includegraphics[width=\textwidth]{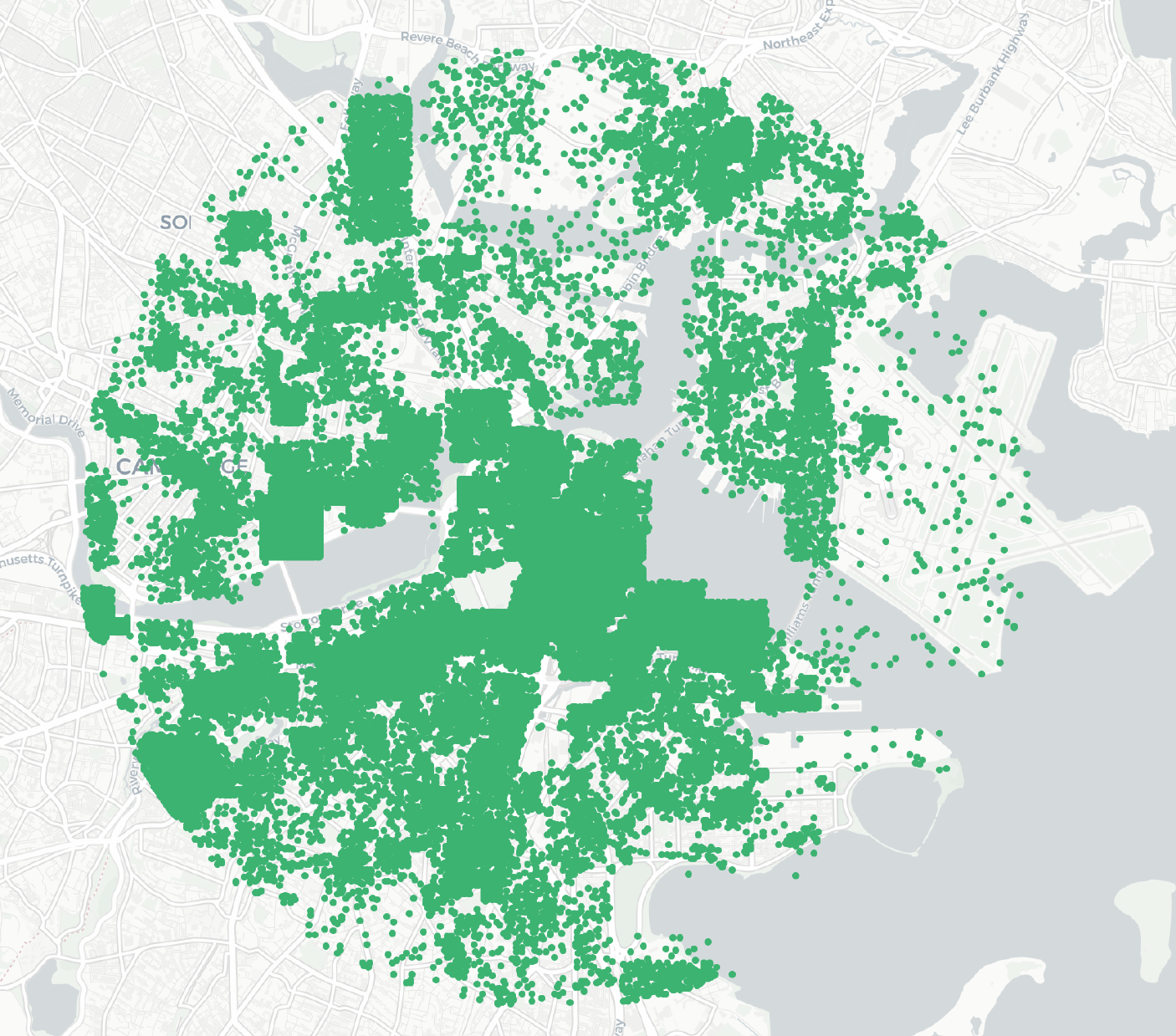}
        \subcaption{Drop-off locations}
    \end{subfigure}

    \begin{subfigure}[b]{0.33\textwidth}
        \textbf{\large Chicago}\\[2pt]
        \includegraphics[width=\textwidth]{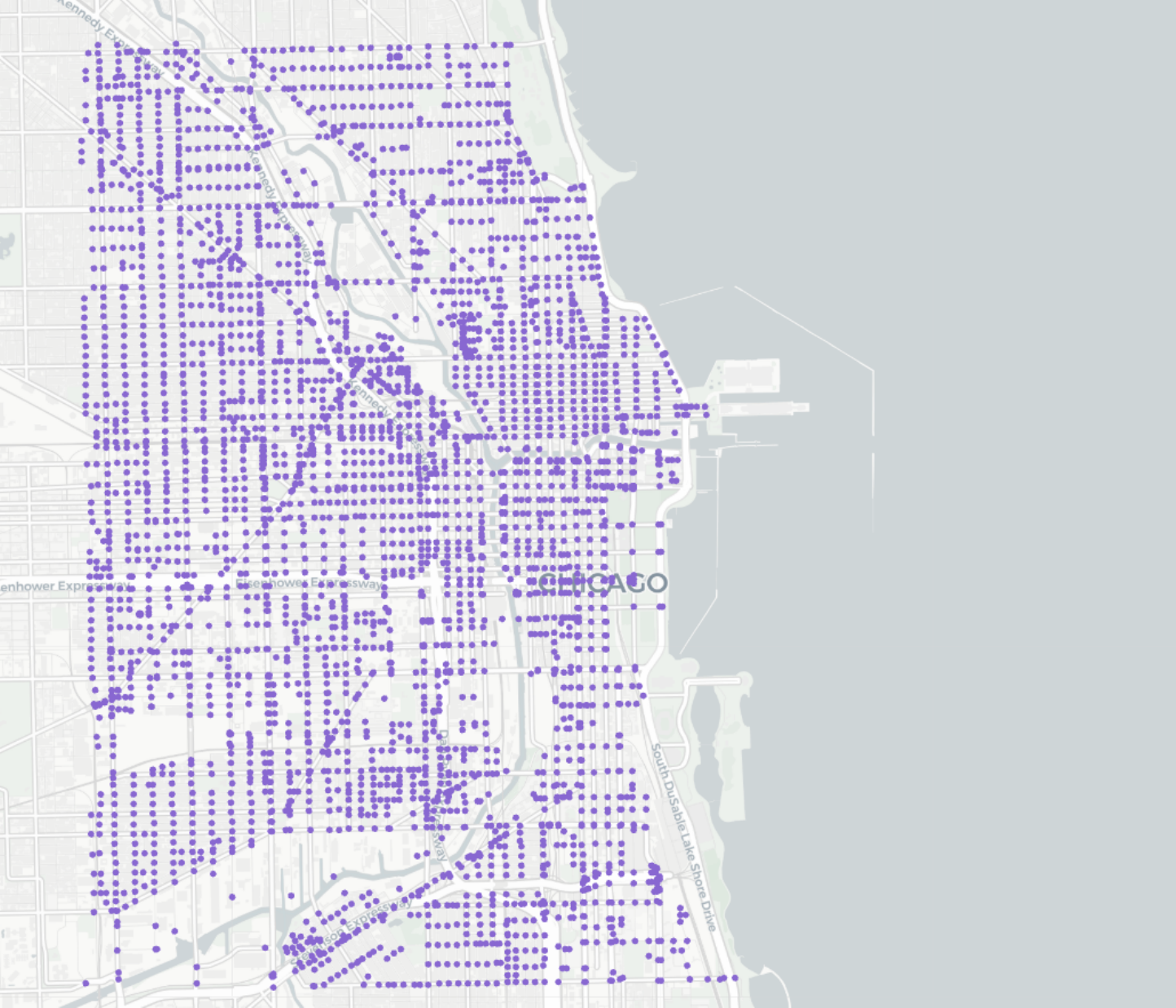}
        \subcaption{OpenStreetMap nodes}
    \end{subfigure}\hfill
    \begin{subfigure}[b]{0.33\textwidth}
        \includegraphics[width=\textwidth]{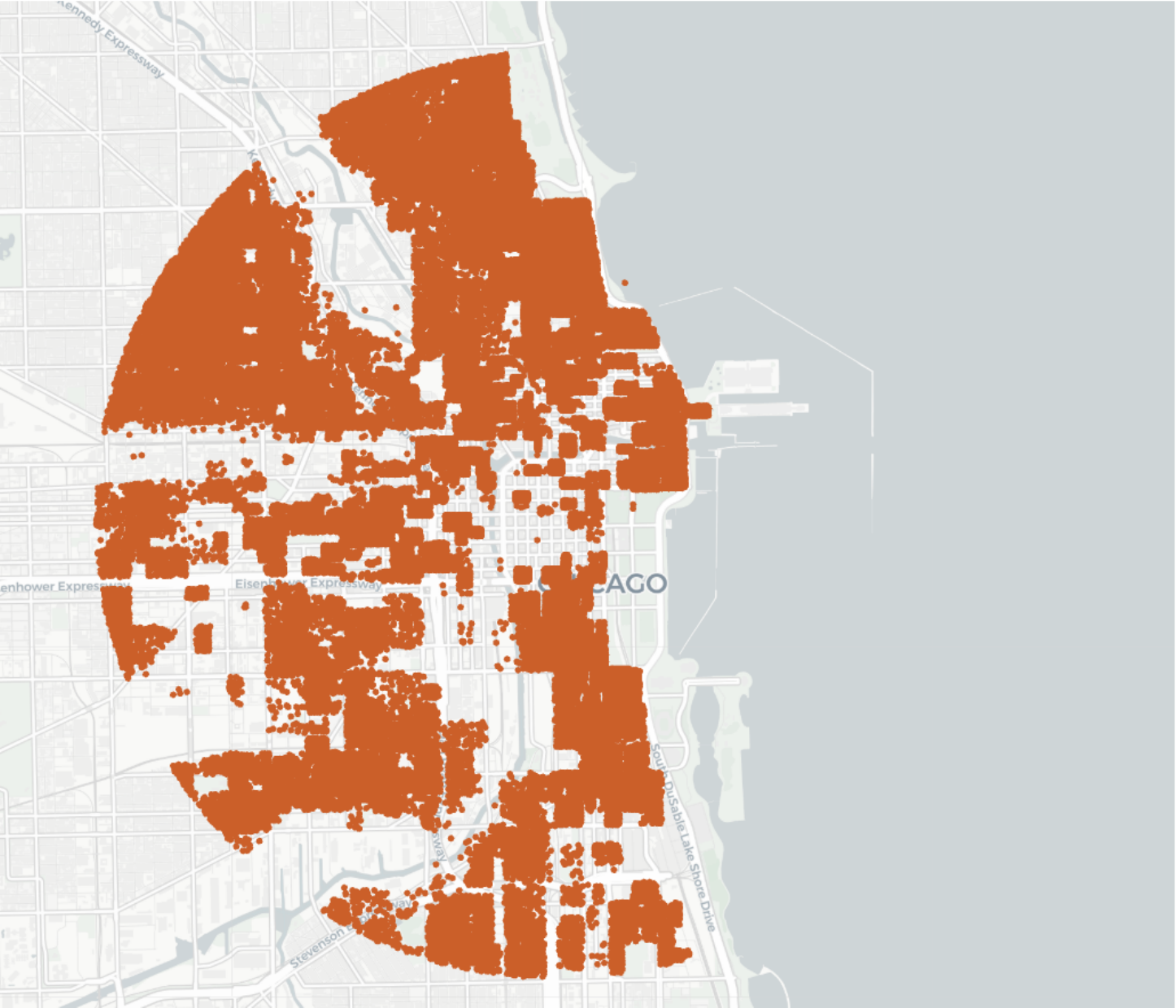}
        \subcaption{Pickup locations}
    \end{subfigure}\hfill
    \begin{subfigure}[b]{0.33\textwidth}
        \includegraphics[width=\textwidth]{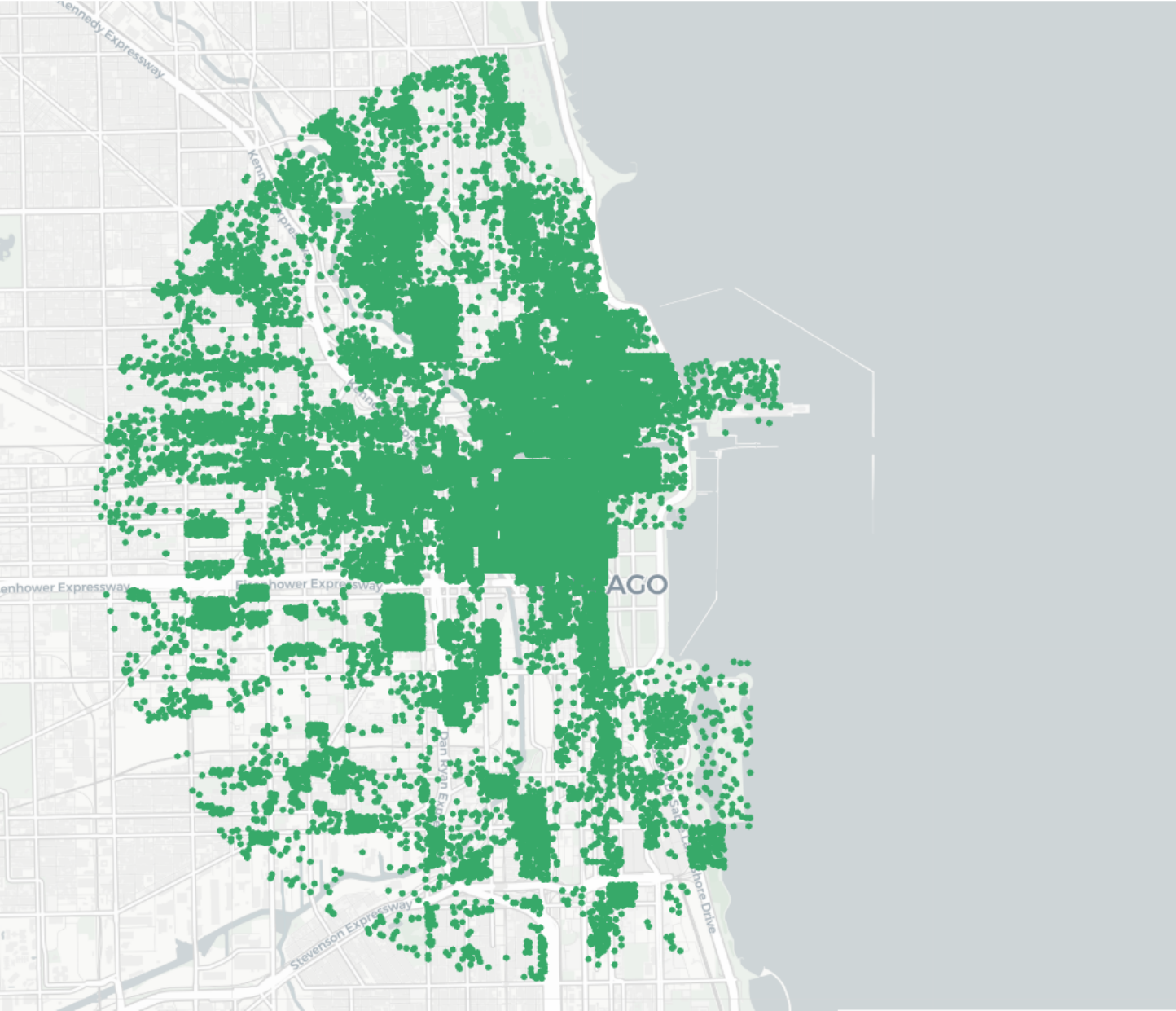}
        \subcaption{Drop-off locations}
    \end{subfigure}
    
    \begin{subfigure}[b]{0.33\textwidth}
        \textbf{\large Atlanta}\\[2pt]
        \includegraphics[width=\textwidth]{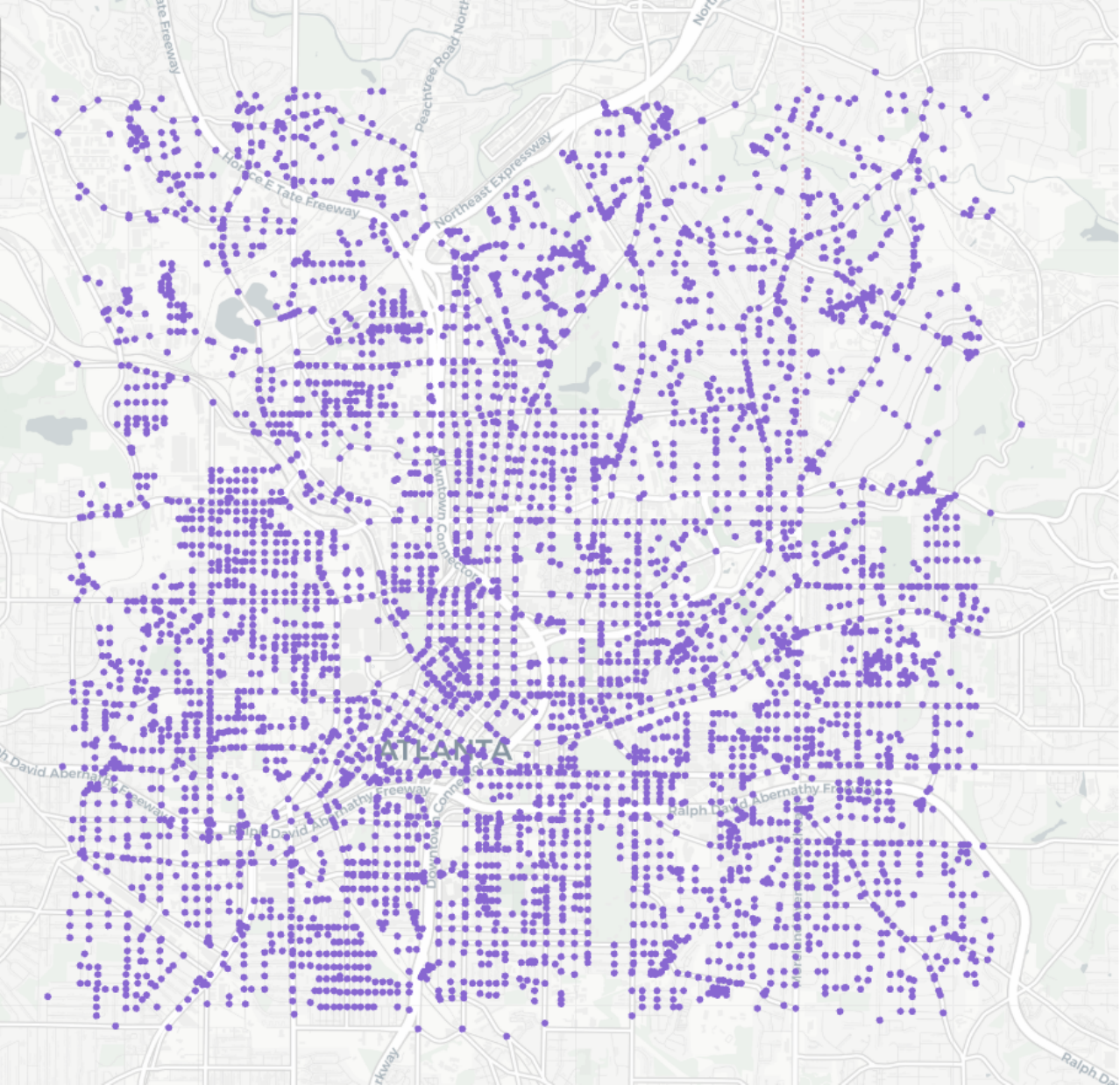}
        \subcaption{OpenStreetMap nodes}
    \end{subfigure}\hfill
    \begin{subfigure}[b]{0.33\textwidth}
        \includegraphics[width=\textwidth]{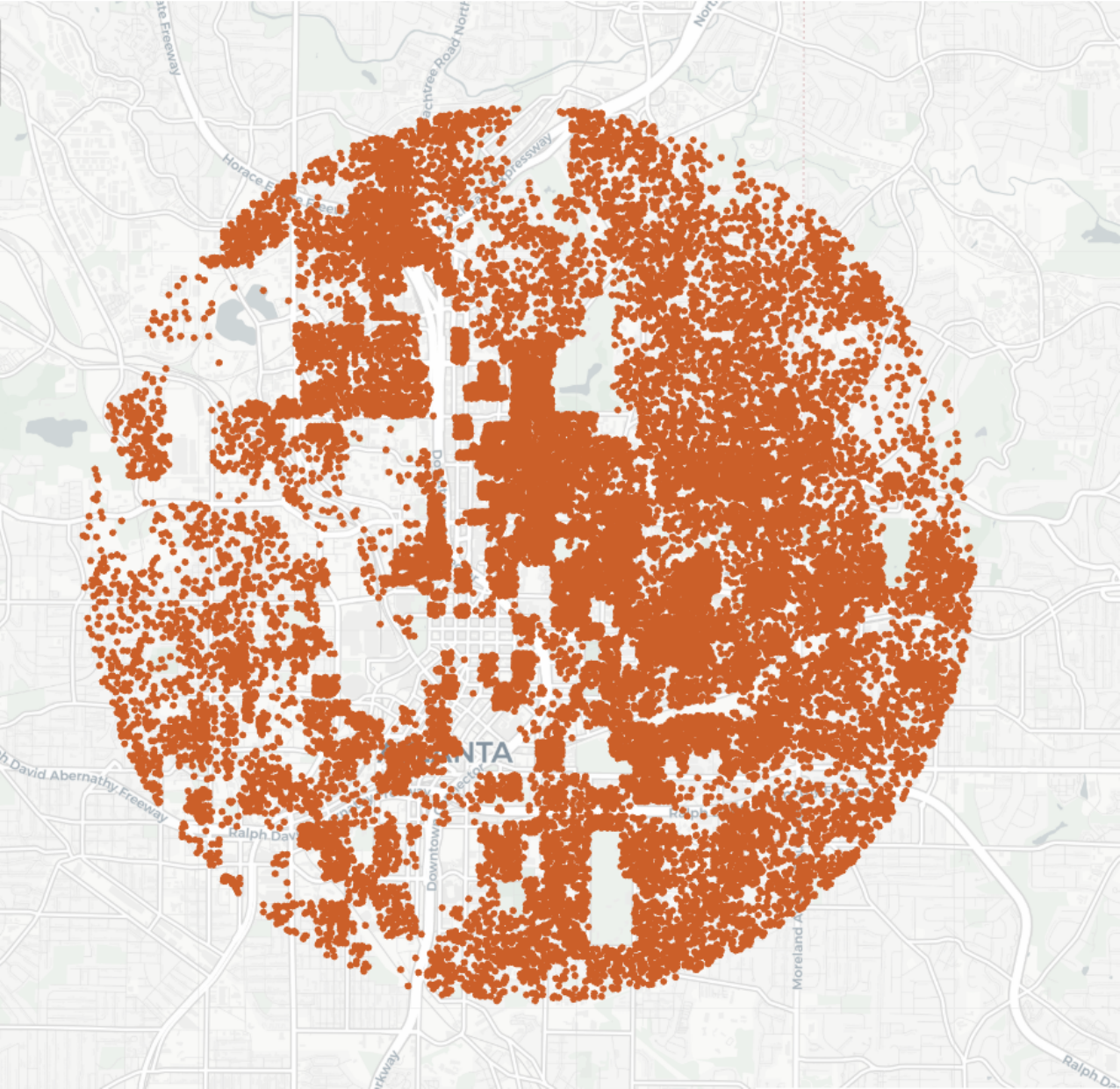}
        \subcaption{Pickup locations}
    \end{subfigure}\hfill
    \begin{subfigure}[b]{0.33\textwidth}
        \includegraphics[width=\textwidth]{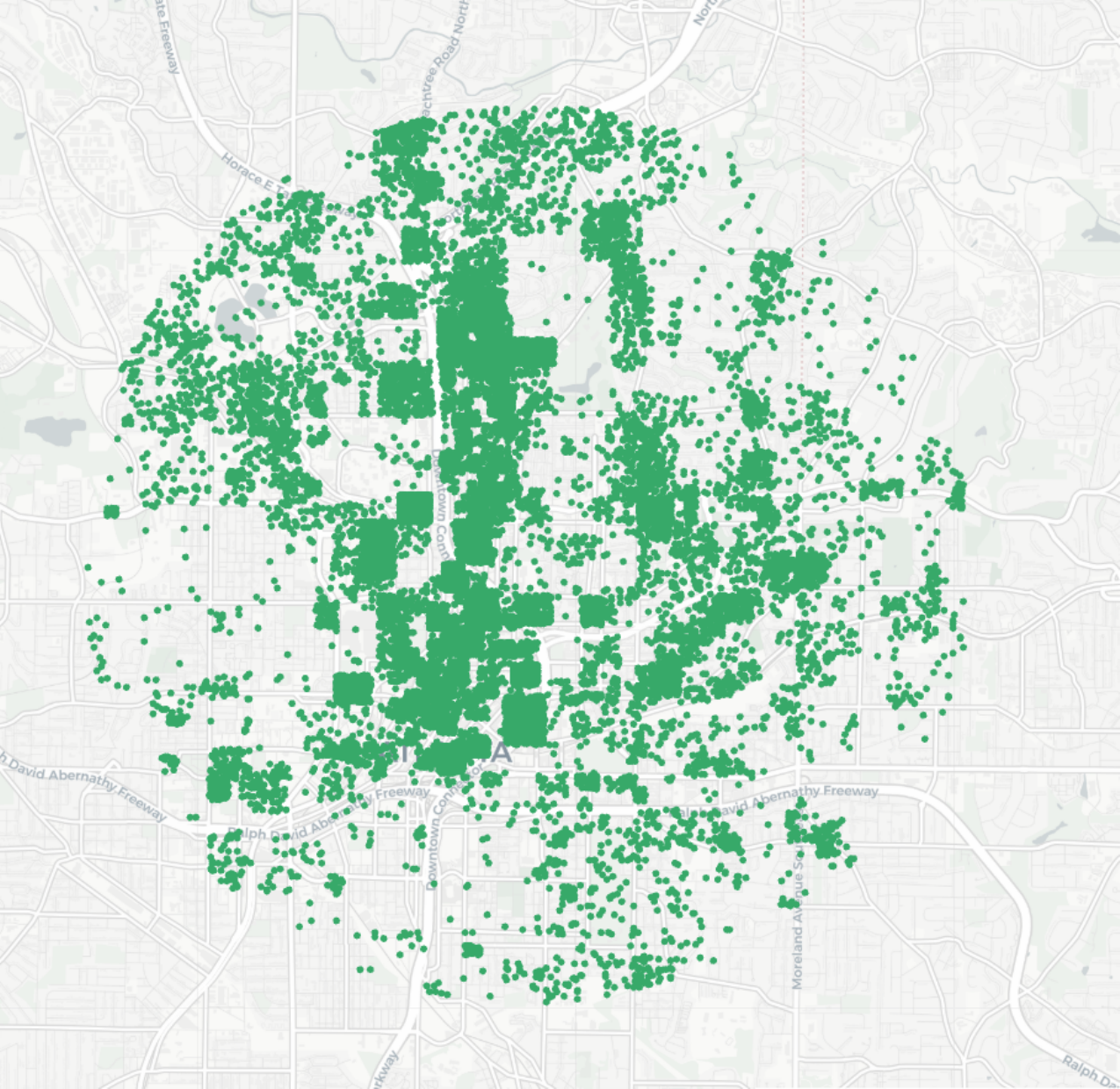}
        \subcaption{Drop-off locations}
    \end{subfigure}
    
    \caption{Initial city networks obtained from OpenStreetMap (within 3 miles of the city centers, omitting edges) and pickup-drop-off demand (in the longitude-latitude format) obtained from LODES. (a), (b), (c): Boston has 7,149 OSM nodes, 15,620 OSM edges, and 99,982 pickup-drop-off demand requests. (d), (e), (f): Chicago has 3,845 OSM nodes, 9,835 OSM edges, and 108,126 pickup-drop-off demand requests. (g), (h), (i): Atlanta has 5,034 OSM nodes, 13,251 OSM edges, and 29,120 pickup-drop-off demand requests.}
    \label{fig:input_data}
\end{figure}

The Boston dataset contains 99,982 passenger records, each indicated by a pickup-drop-off location pair. After obtaining the Boston network from OpenStreetMap, we deliberately exclude road classes of ``motorway", ``motorway link", ``living street" and ``unclassified", as they are unlikely to support \bus~or \taxi~services. Moreover, our focus is on the dense downtown area, where the excluded road classes are rarely present. We select the remaining road classes: ``primary", ``primary link", ``secondary", ``secondary link", ``tertiary", ``tertiary link", ``trunk", ``trunk link", ``residential", ``busway" and ``crossing". We then select the largest strongly connected component after removing self-loops. The resulting network contains 7,149 nodes and 15,620 edges.  \Cref{fig:input_data} illustrates the network nodes and the pickup and drop-off locations of demand requests in Boston, Chicago, and Atlanta. Let $X$ denote the demand set of 99,982 records. Let $N$ denote the set of 7,149 nodes. We select a subset of the 7,149 nodes as virtual stops such that each passenger we serve walks no more than a walking-distance threshold  at both ends of their trip. We set the walking distance to be $400$ meters. The selection has five steps. The resulting set of virtual stops becomes our set of transportation nodes $\nodesT$ on which \taxi~vehicles travel, as defined in the problem definition in \Cref{sec:statement}.

\begin{itemize}
\item Step 1: Obtain $X'$ by filtering $X$ and retaining only trips whose nearest nodes in $N$ to both the pickup and drop-off locations lie within a walking distance of $400$ meters. $X'$ contains 99,756 records, each with a pickup location and a drop-off location.

\item Step 2: Form clusters of all pickup and drop-off locations in $X'$. Each cluster is centered at a location in $X'$, and consists of all locations within 100 meters of that center. Obtain the set of clusters $C$, which has 3,192 elements. 

\item Step 3: Let $\text{distance}(u,v)$ denote the haversine distance between a node $u\in N$ and a location $v$, where $v$ is the center of a cluster in $C$. For each $u\in N$, let $\text{nbh}(u)=\{v \text{ is the center of a cluster in }C: \text{distance}(u,v) \leq 400 $\text{ meters}$\}$ be the neighboring cluster centers for the node. We solve the following hitting set problem \cite{Gainer2017Minimal} to select a subset of nodes $N'\subset N$ covering cluster centers. Let $x_u \in \{0,1\}$ for each $u\in N$ be the decision variable indicating whether node $u$ is included in $N'$. 

\begin{align}
\min \ \sum_{u\in N}&x_{\nodeU}\\
s.t. \ \sum_{\nodeU: \nodeV\in \text{nbh}(\nodeU)} x_{\nodeU} & \geq 1 &&\quad \forall \nodeV\in \{\text{cluster centers in }C\} \\
x_{\nodeU} &\in \{0,1\} &&\quad \forall \nodeU\in N
\end{align}

We set the termination gap to be $10\%$. The resulting solution selects 166 nodes, which form the set of virtual stops $N'$. 

\item Step 4: For each record in $X$, find the nearest node in $N'$ to each pickup and drop-off location. Remove records in $X$ if the nearest node to either the pickup or drop-off is more than $400$ meters away, resulting in 95,921 records. 

\item Step 5: For each of the 95,921 records, match the pickup and drop-off locations to the nearest node in $N'$. They are the origin and destination nodes that each passenger would be willing to walk to. Remove records if pickup and drop-off locations are matched to the same node, resulting in 92,975 records across 11,045 origin-destination (OD) pairs.

\end{itemize}

As described earlier, the resulting set of virtual stops is the set of transportation nodes $\nodesT$, and therefore $\nodesT = N'$ containing 166 nodes. We construct the transportation network $\graphT = (\nodesT, \edgesT)$ to be the directed complete graph on $\nodesT$ where $\edgesT = \{(u,v): u\in \nodesT, v\in \nodesT, u\neq v\}$ containing 27,390 edges (166 $\times$ 165). The distance $\costUV$ for edge $(\nodeU, \nodeV)\in \edgesT$ is the distance of the shortest path from $\nodeU$ to $\nodeV$ in the road network. We regard the demand set of 92,975 records across 11,045 ODs as representing travel demand during a daily 3-hour morning peak period. Bus systems often operate on regular headways, commonly 10, 12, 15, or 20 minutes. We set the interval to be 15 minutes for our design, and we obtain a representative 15-minute demand sample as well as corresponding line-dependent quantities defined in \Cref{sec:statement}. We present two schemes for obtaining the demand sample: probabilistic demand realization and low-demand truncation with rescaling.

\textbf{Probabilistic demand realization:}
For each of the 11,045 OD pairs with 92,975 records, divide the number of demand records by $12$, since a 3-hour period corresponds to twelve 15-minute intervals. Many pairs have demand between 0 and 1. For each OD pair, regard the fractional part of the demand amount as a coin-flip probability. For each OD, we take the floor of the demand and add an additional trip based on the coin-flip probability. This process yields 7,724 demand trips across 4,122 ODs in the 15-minute sample. Scaling back for the full 3-hour peak period, we have 92,688 trips (7,724 × 12) as the input passenger demand we aim to serve.

Since not every node in the transportation node set $\nodesT$ is suitable to serve as a \bus~stop, we select a subset of $\nodesT$ to form the set of \bus~nodes $\nodesB$ and construct the \bus~network $\graphB=(\nodesB, \edgesB)$. We define and solve a facility location \cite{Farahani2010Multiple} problem to select 100 of the 166 nodes in $\nodesT$ to form $\nodesB$, where $100$ is a parameter. (Maintaining the approximate proportionality, we select 70 \bus~nodes out of 117 transportation nodes for Chicago and 117 \bus~nodes out of 183 transportation nodes for Atlanta.) Specifically, let $\text{distance}(i,j)$ be the shortest-path distance between $i\in \nodesT$ and $j\in \nodesT$ in the road network. For each node $j\in \nodesT$, let $\text{volume}(j)$ be the total arrival and departure demand from the 3-hour demand dataset of 92,688 trips. For each $i\in \nodesT$, let $z_i\in \{0,1\}$ be the decision variable indicating whether node $i$ is selected to be one of the 100 \bus~nodes. For each $i,j\in \nodesT$, let $y_{i,j}\in \{0,1\}$ be the decision variable indicating whether demand at node $j$ can potentially be served by constructing a \bus ~stop at node $i$. Let $M$ be a very large constant. We solve the following program to optimality to select $100$ nodes corresponding to the optimal values of $z_i$:

\begin{align}
\min \ \sum_{i,j\in \nodesT}\text{distance}(i,j) &\text{volume}(j) y_{i,j}\\
s.t. \ \sum_{i\in \nodesT} y_{i,j} & = 1 &&\quad \forall j\in \nodesT \\
\sum_{j\in \nodesT} y_{i,j} &\leq M z_i&&\quad \forall i\in \nodesT\\
\sum_{i\in \nodesT} z_i &= 100\\
z_i &\in \{0,1\} &&\quad \forall i\in \nodesT\\
y_{i, j} &\in \{0,1\} &&\quad \forall i\in \nodesT,~\forall j \in \nodesT
\end{align}

After getting the \bus~node set produced by the optimal solution of the above program, we construct a directed complete graph on the 100 nodes. Each edge weight represents the shortest-path distance between the nodes in the original road network. If either of the two directed edges between a node pair has a shortest-path distance greater than or equal to 2,000 meters, we classify the pair as a long-edge pair. We iteratively remove such pairs from the directed complete graph on the 100 bus nodes, starting with the longest edge. Whenever an edge is removed, its opposite directed edge is also removed. We continue removing edges until no long-edge pairs remain and the graph remains strongly connected, resulting in 786 bus edges in $\edgesB$ satisfying connectivity \Cref{assump:connectivity}.

\textbf{Low-demand truncation with rescaling:} For each of the 11,045 OD pairs with 92,975 records, divide the number of demand records by $12$, since a 3-hour period corresponds to twelve 15-minute intervals. Many pairs have demand between 0 and 1. Since LODES does not provide a complete survey of location data, we regard the dataset as an underestimate of actual travel demand. In order to build transportation systems to serve high-corridor demand, we disregard ODs with no more than one demand record. For the remaining ODs, we round up to the next integer. The process results in 6,136 demand trips across 1,663 ODs. Scaling back for the full 3-hour peak period, we have 73,632 trips (6,136 $\times$ 12) as the input demand we aim to serve. We select the bus nodes $\nodesB$ and bus edges $\edgesB$ for the low-demand truncation with rescaling using the same procedure as for the probabilistic demand realization. The selection results in 100 \bus~nodes and 794 \bus~edges.

We conduct our experimental analysis using both demand-processing schemes for each city. The two schemes produce different transportation demand sets and different \bus~networks $\graphB = (\nodesB, \edgesB)$, since the selection of \bus~nodes involves solving a demand-dependent facility location problem. All other inputs are set to be the same for the two schemes for each city, including the transportation network $\graphT=(\nodesT, \edgesT)$. For each city and demand-processing scheme, we generate our \bus~lines from its \bus~network $\graphB=(\nodesB, \edgesB)$. 

\begin{figure}[H]   
    \noindent\makebox[\textwidth][l]{\textbf{\large Boston}}\\
    \begin{subfigure}[t]{0.25\textwidth}
        \includegraphics[width=\textwidth]{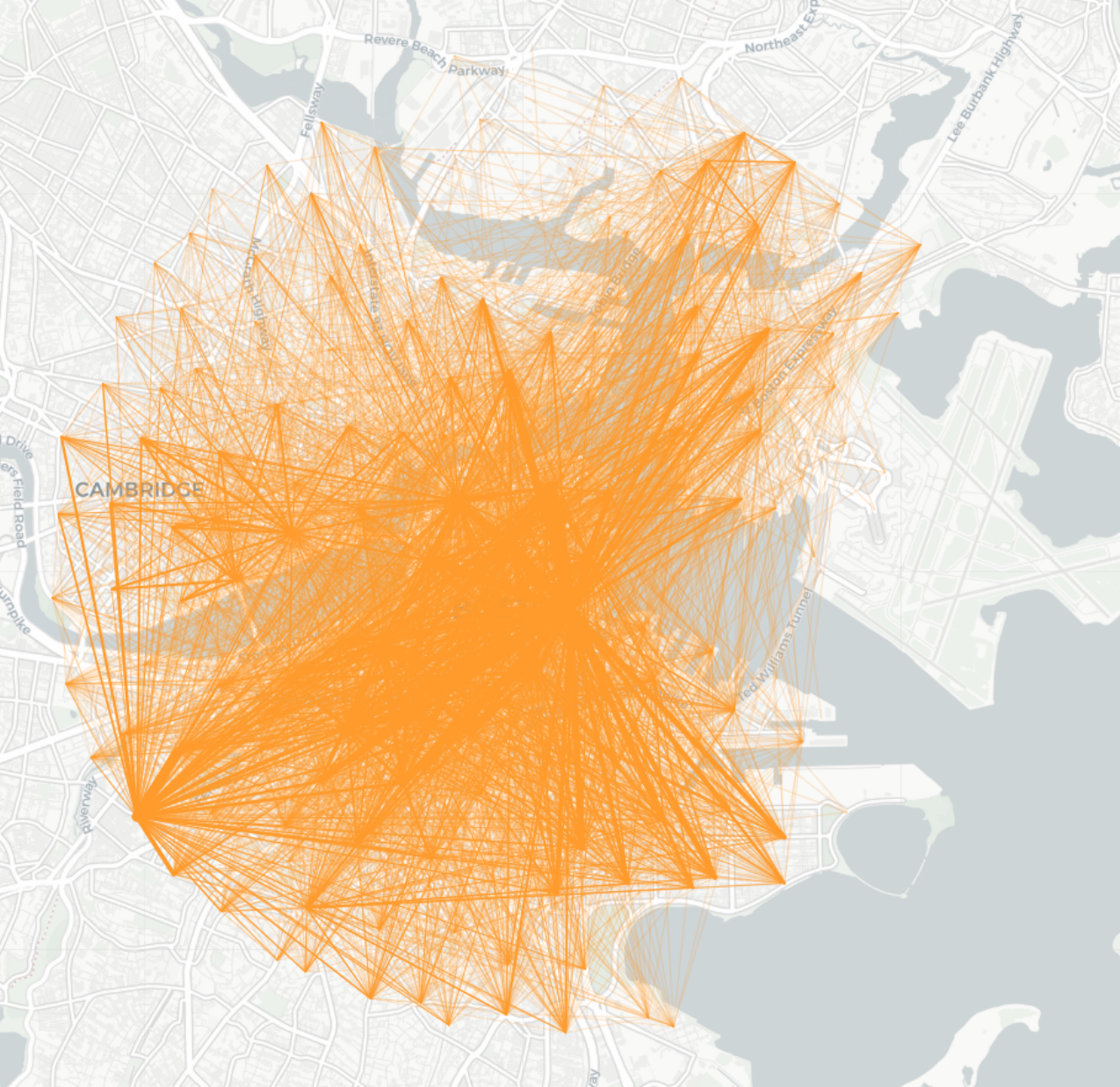}
        \subcaption{Prob. demand realization}
    \end{subfigure}%
    \begin{subfigure}[t]{0.25\textwidth}
        \includegraphics[width=\textwidth]{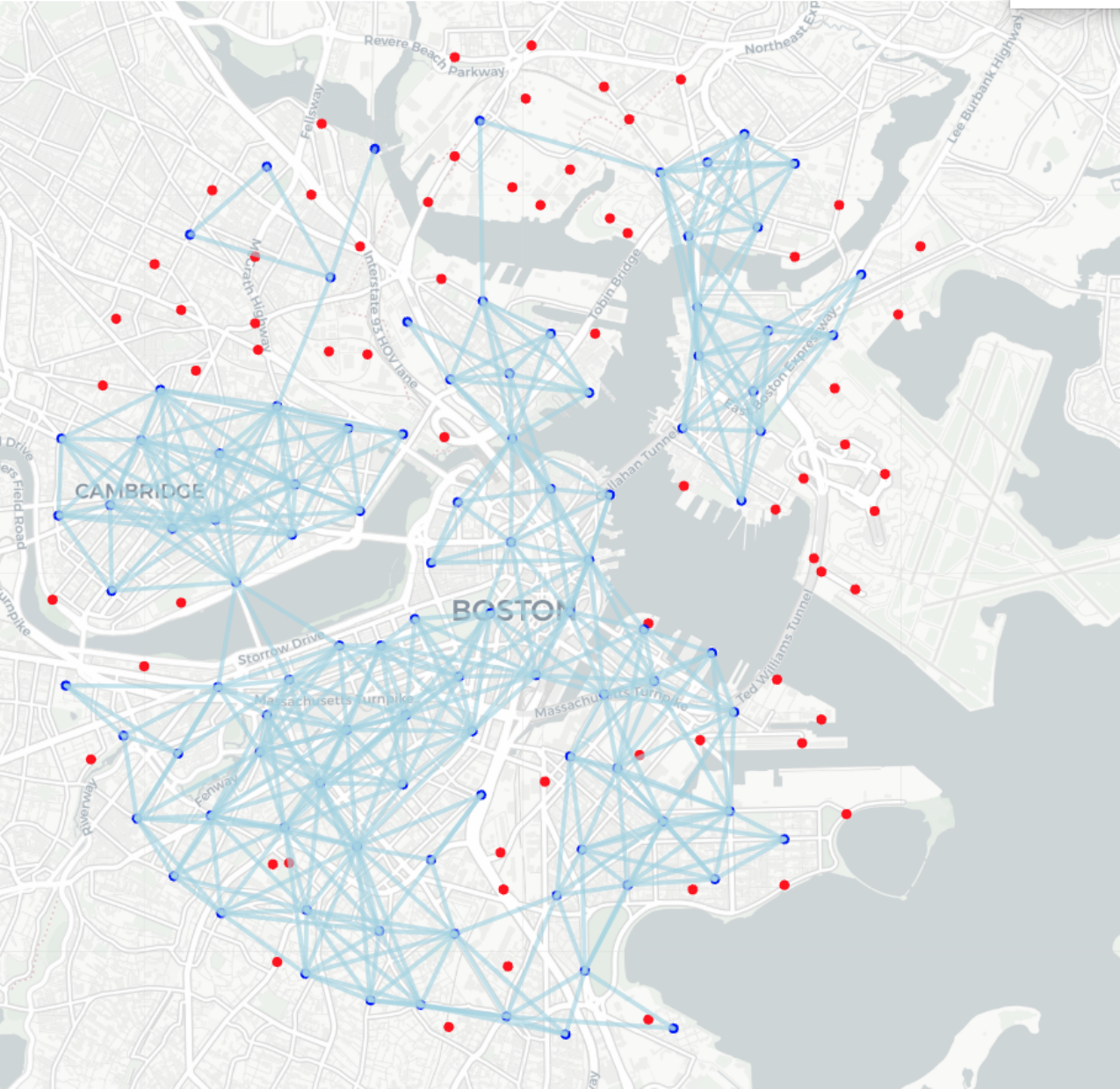}
        \subcaption*{Demand based network}
    \end{subfigure}%
    \hspace{1em}
    \begin{subfigure}[t]{0.25\textwidth}
        \includegraphics[width=\textwidth]{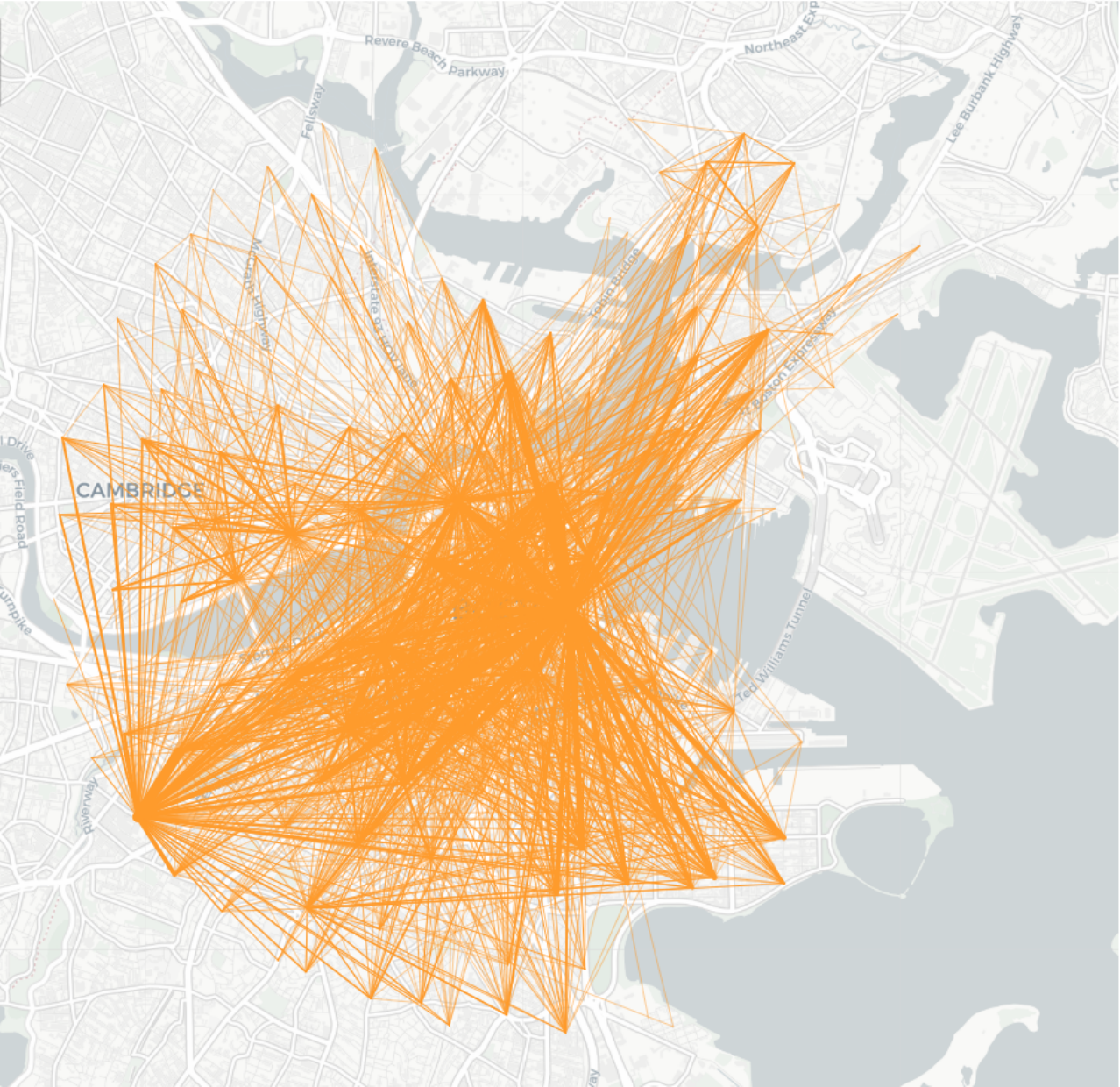}
        \subcaption{Low-demand truncation}
    \end{subfigure}%
    \begin{subfigure}[t]{0.25\textwidth}
        \includegraphics[width=\textwidth]{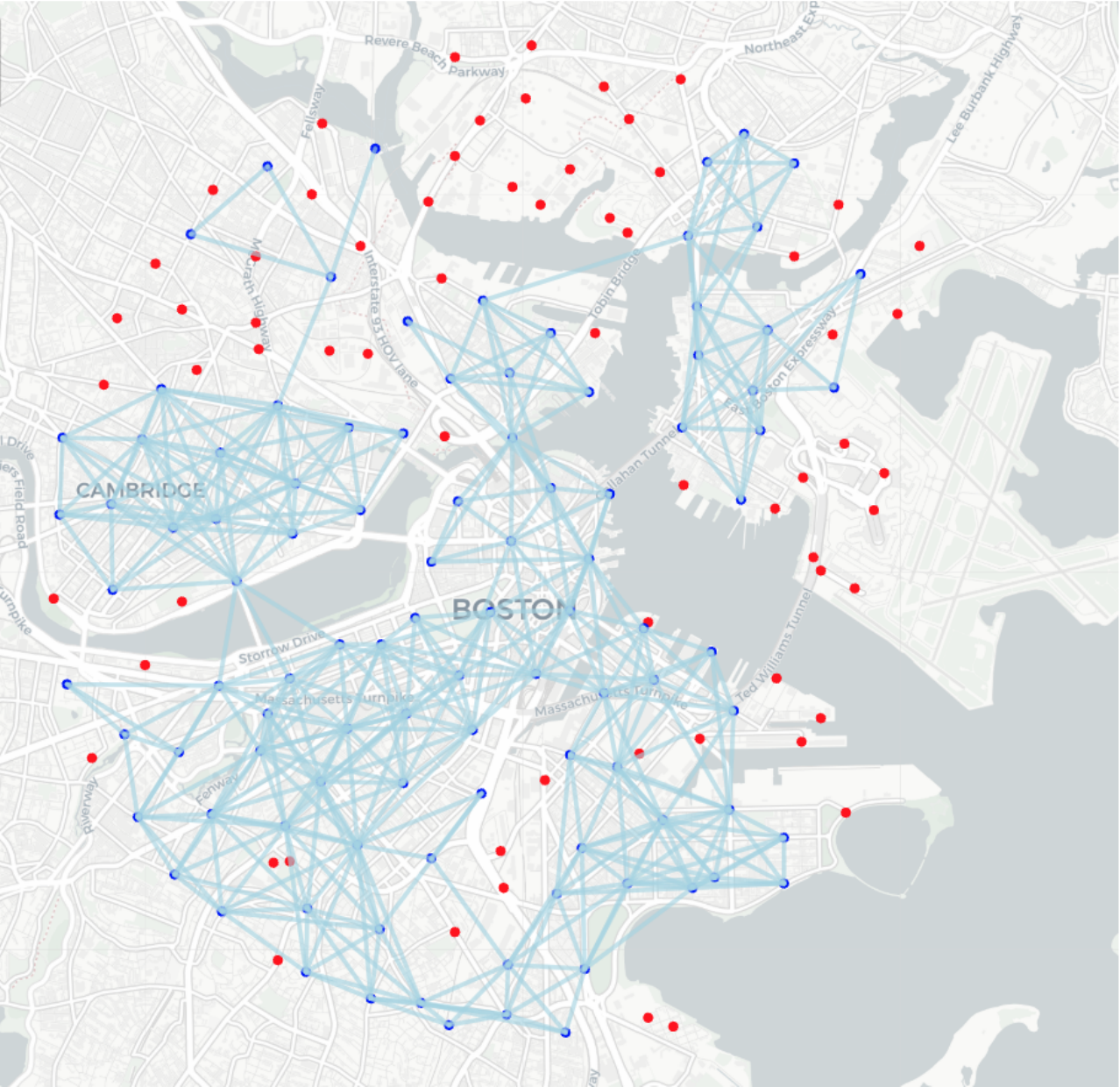}
        \subcaption*{Demand based network}
    \end{subfigure}\\[1em]
    
    \noindent\makebox[\textwidth][l]{\textbf{\large Chicago}}\\
    \begin{subfigure}[t]{0.25\textwidth}
        \includegraphics[width=\textwidth]{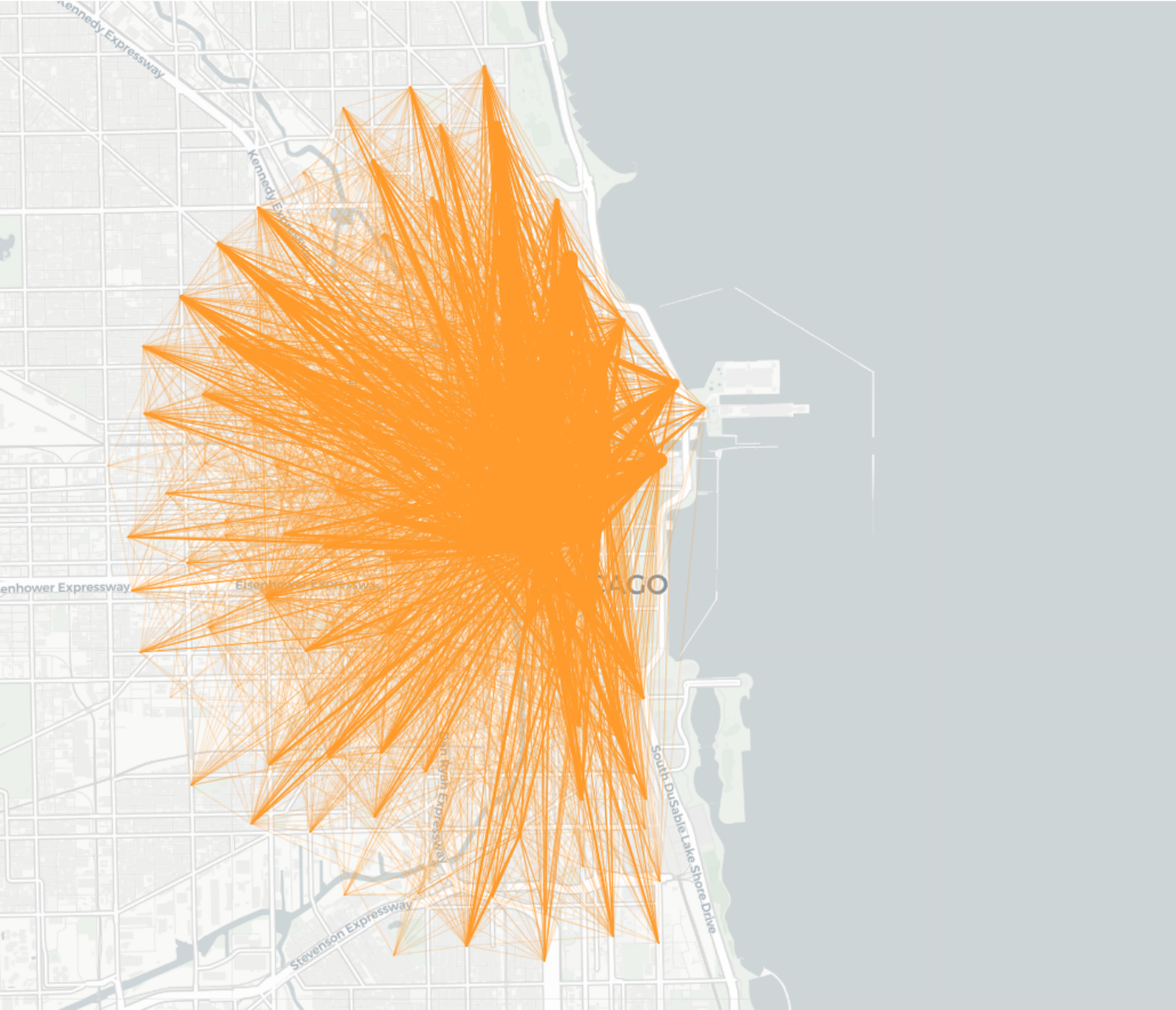}
        \subcaption{Prob. demand realization}
    \end{subfigure}%
    \begin{subfigure}[t]{0.25\textwidth}
        \includegraphics[width=\textwidth]{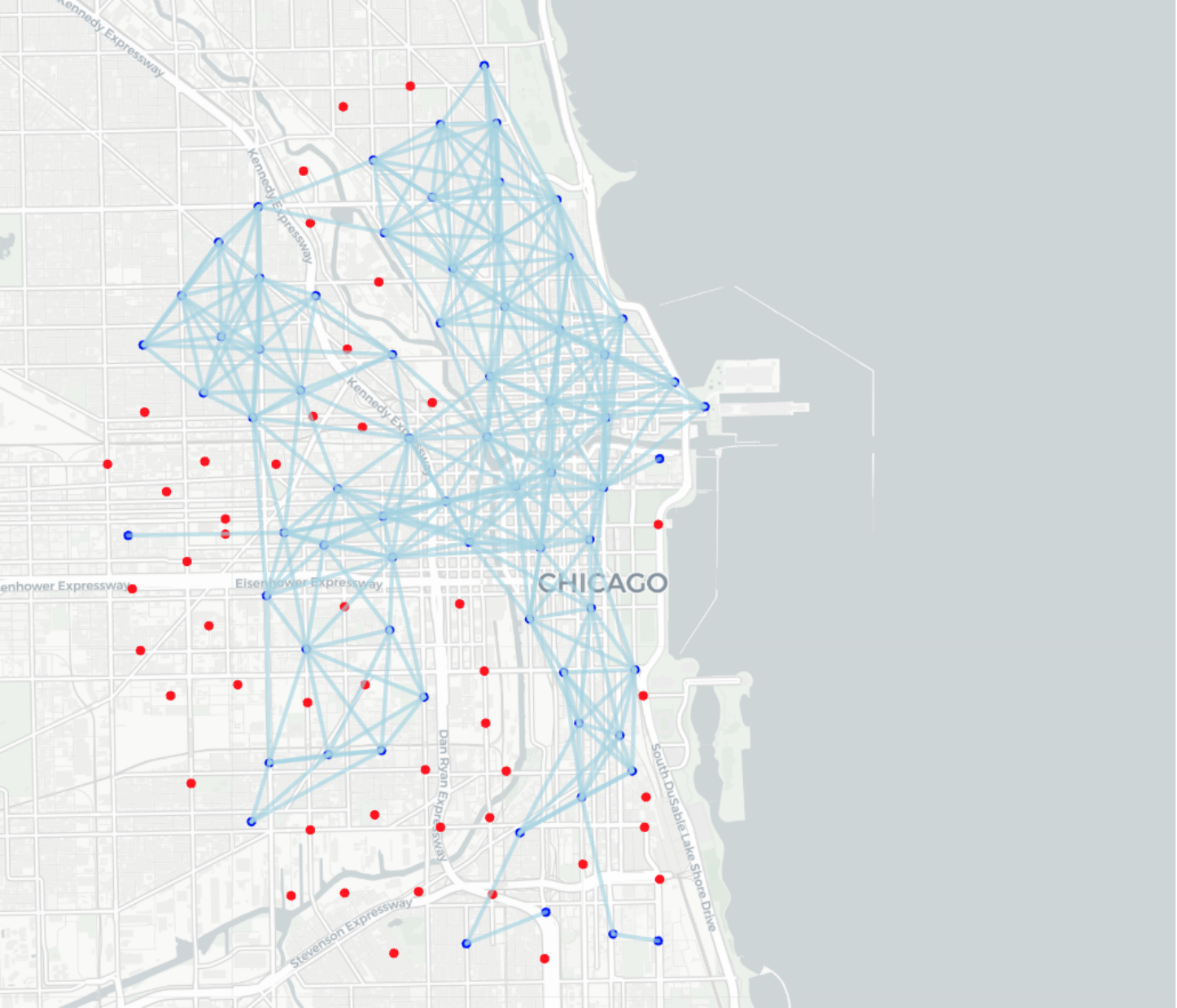}
        \subcaption*{Demand based network}
    \end{subfigure}%
    \hspace{1em}
    \begin{subfigure}[t]{0.25\textwidth}
        \includegraphics[width=\textwidth]{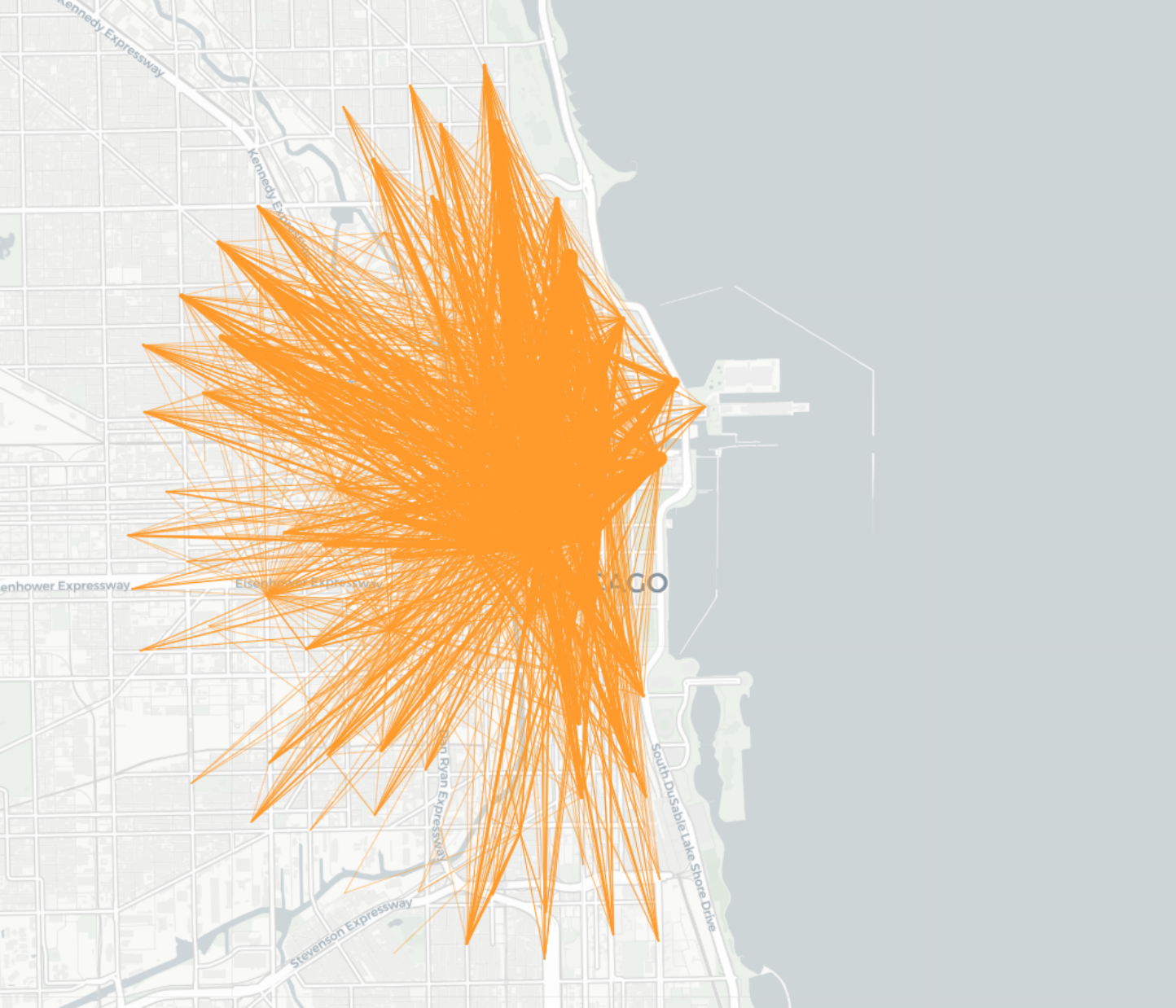}
        \subcaption{Low-demand truncation}
    \end{subfigure}%
    \begin{subfigure}[t]{0.25\textwidth}
        \includegraphics[width=\textwidth]{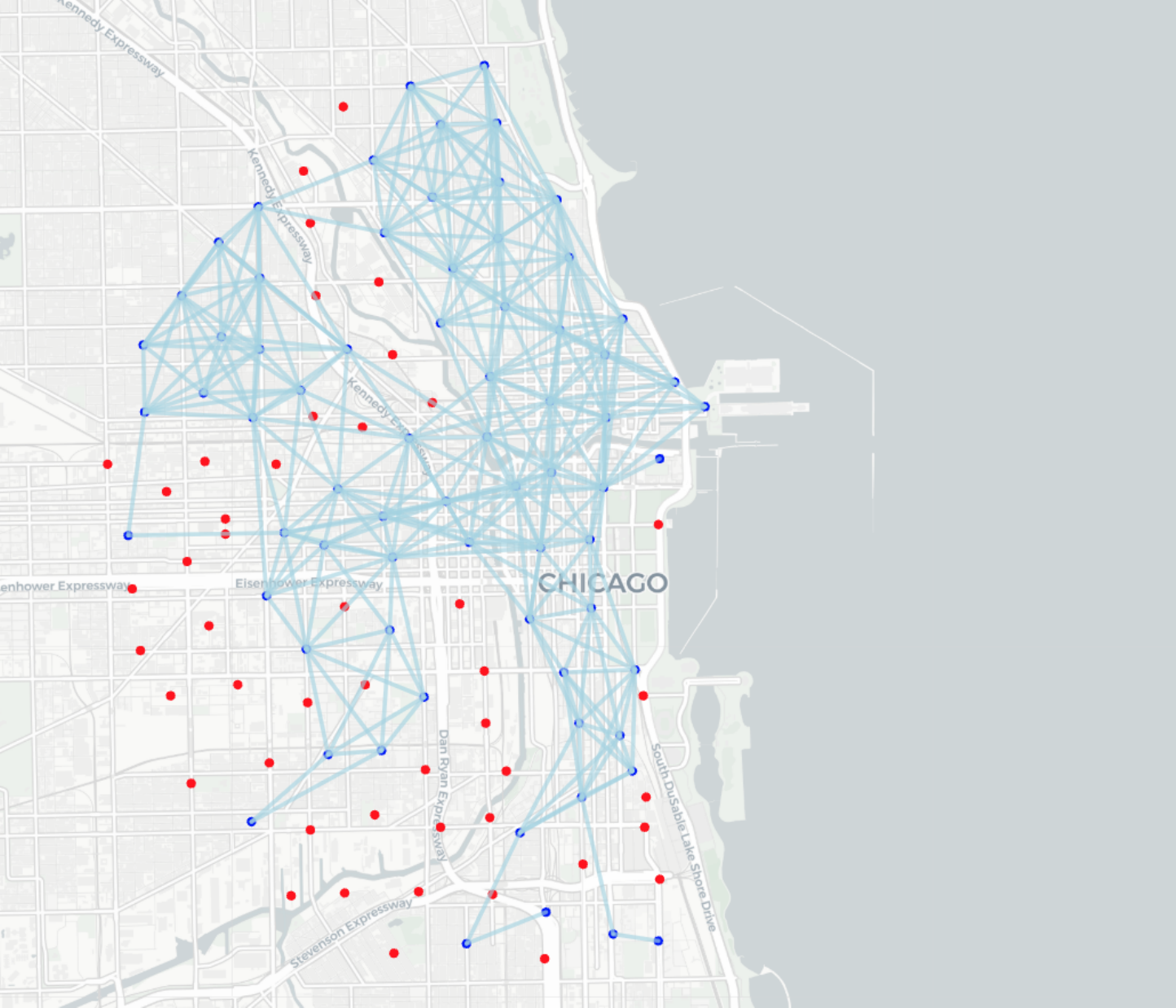}
        \subcaption*{Demand based network}
    \end{subfigure}\\[1em]
    
    \noindent\makebox[\textwidth][l]{\textbf{\large Atlanta}}\\
    \begin{subfigure}[t]{0.25\textwidth}
        \includegraphics[width=\textwidth]{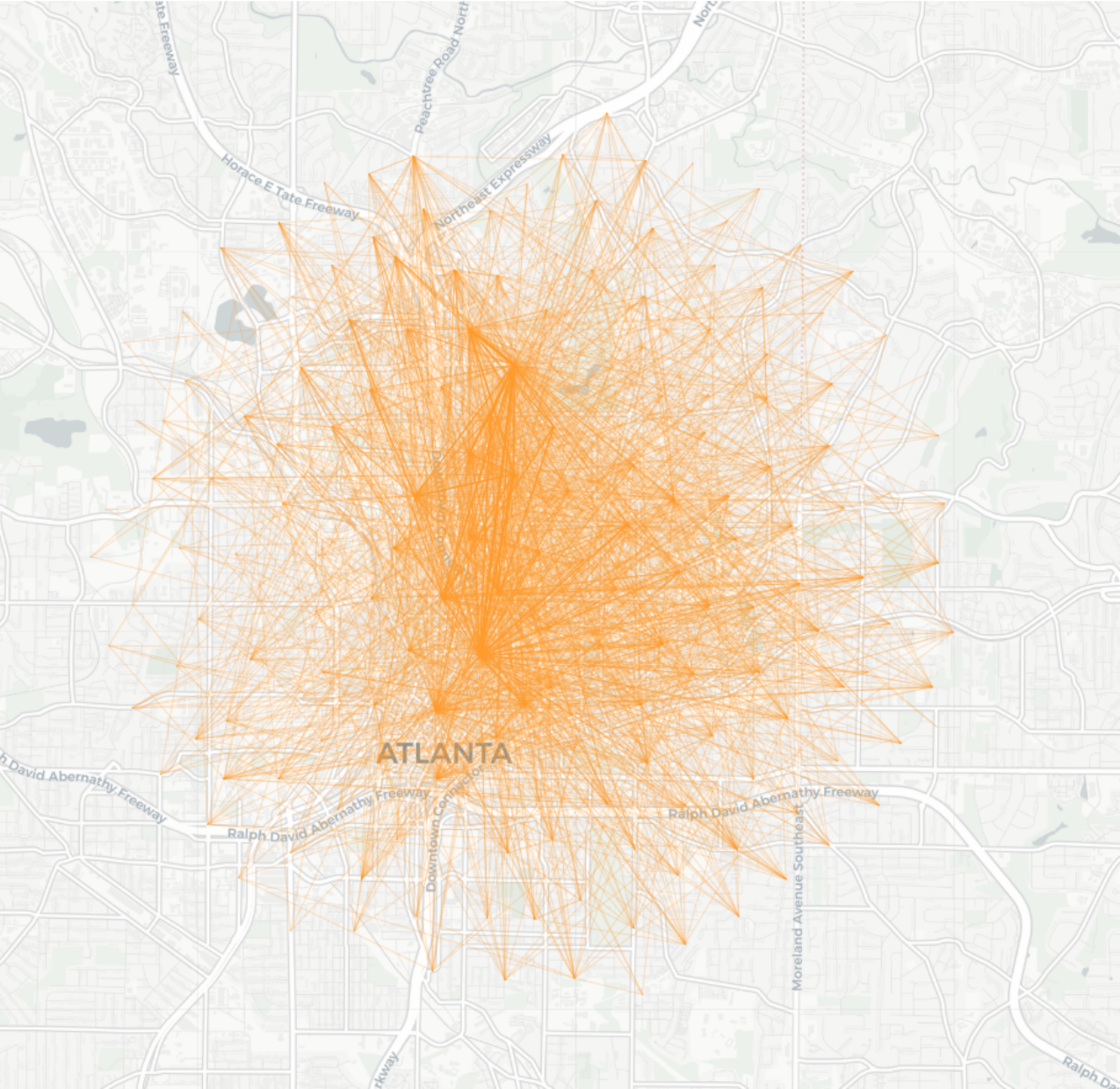}
        \subcaption{Prob. demand realization}
    \end{subfigure}%
    \begin{subfigure}[t]{0.25\textwidth}
        \includegraphics[width=\textwidth]{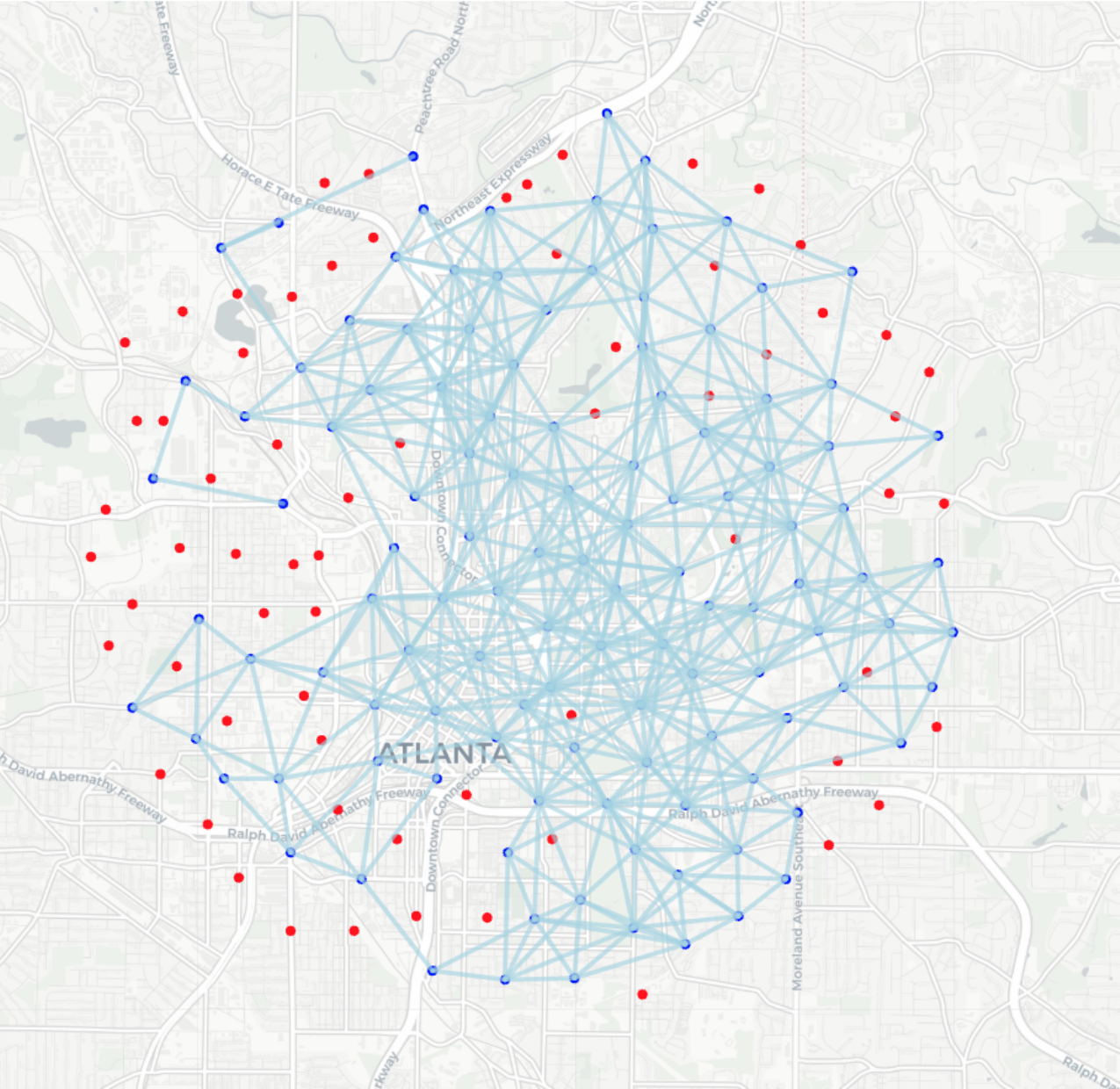}
        \subcaption*{Demand based network}
    \end{subfigure}%
    \hspace{1em}
    \begin{subfigure}[t]{0.25\textwidth}
        \includegraphics[width=\textwidth]{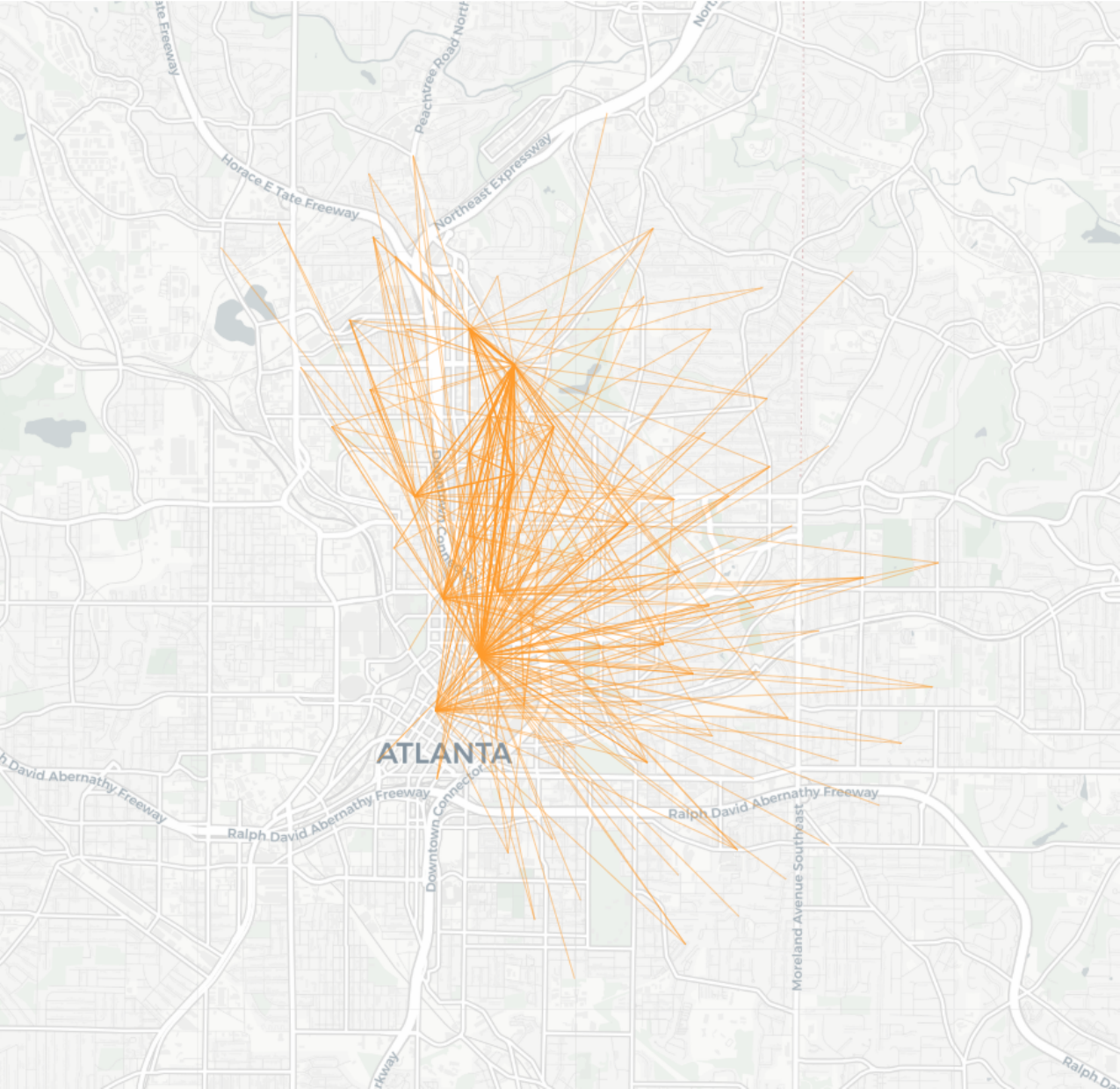}
        \subcaption{Low-demand truncation}
    \end{subfigure}%
    \begin{subfigure}[t]{0.25\textwidth}
        \includegraphics[width=1.015\textwidth]{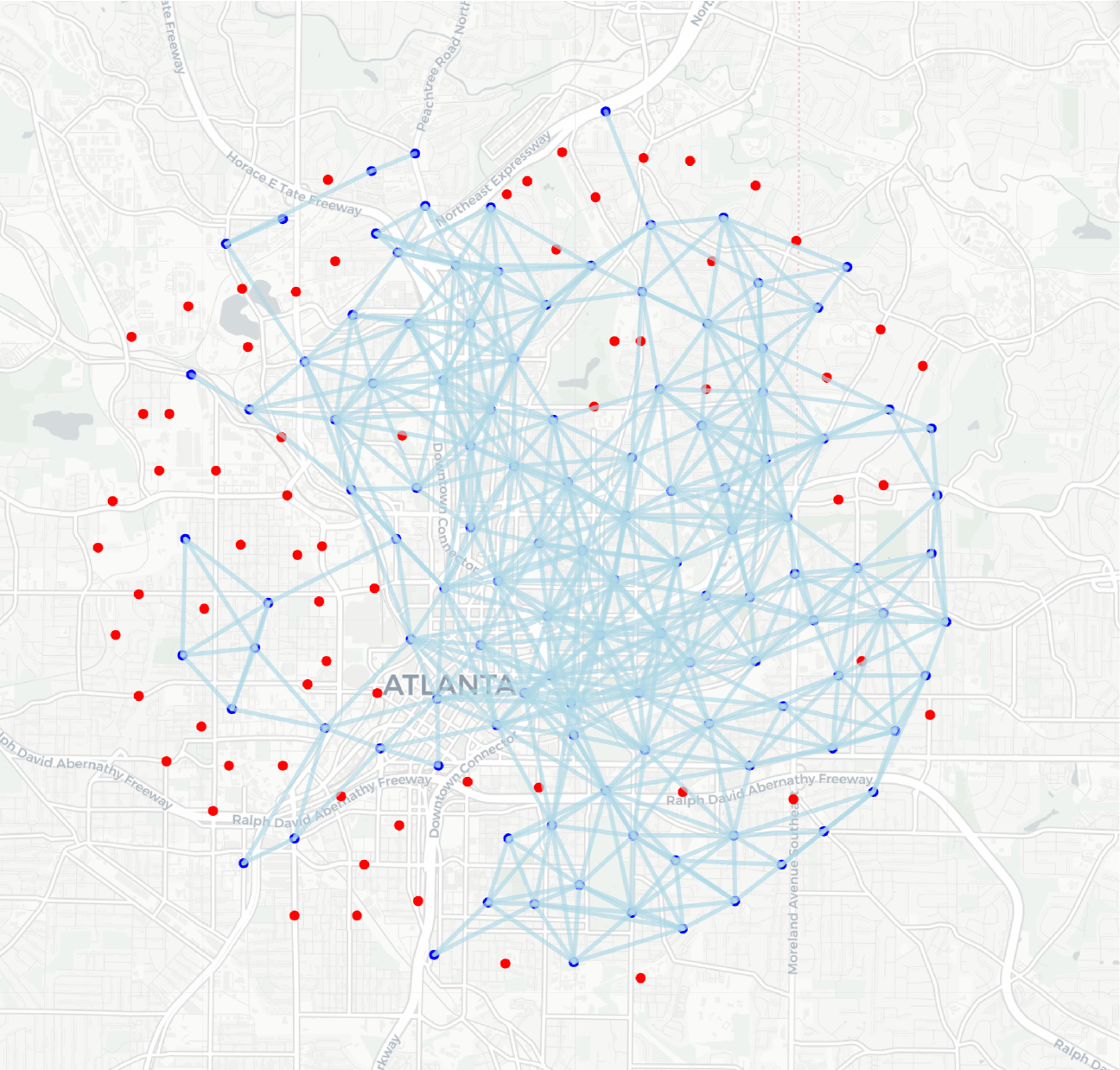}
        \subcaption*{Demand based network}
    \end{subfigure}
    \caption{Visualization of the probabilistic demand realization and the low-demand truncation with rescaling for Boston, Chicago, and Atlanta, as well as the networks $\graphT$ and $\graphB$. Each pair of figures has the demand visualization on the left and the network on the right. In each demand visualization on the left, each yellow edge connects an OD pair (without directional information), where the edge thickness indicates the relative demand magnitude in either direction between the endpoints. The network on the right of each pair of figures consists of transportation network nodes $\nodesT$ (red and blue), \bus~nodes $\nodesB$ (blue), and \bus~edges $\edgesB$ (light blue). We omit the visualization of the transportation network edge set $\edgesT$, which includes all directed edges in the complete graph formed by nodes in $\nodesT$. Each light blue edge represents a pair of directed edges in opposite directions between the endpoints. The underlying data in each pair of figures are the following:\\ (a) 92,688 demand trips across 4,122 ODs, $|\nodesT| = 166$, $|\edgesT| = $27,390, $|\nodesB| = 100$, $|\edgesB| = 786$; (b) 73,632 demand trips across 1,663 ODs, $|\nodesT| = 166$, $|\edgesT| = $27,390, $|\nodesB| = 100$, $|\edgesB| = 794$; (c) 100,296 demand trips across 2,900 ODs, $|\nodesT| = 117$, $|\edgesT| = $13,572, $|\nodesB| = 70$, $|\edgesB| = 602$; (d) 90,348 demand trips across 1,384 ODs, $|\nodesT| = 117$, $|\edgesT| = $13,572, $|\nodesB| = 70$, $|\edgesB| = 620$;\\ (e) 27,084 demand trips across 2,012 ODs, $|\nodesT| = 183$, $|\edgesT| = $33,306, $|\nodesB| = 117$, $|\edgesB| = $1,052; (f) 8,928 demand trips across 306 ODs, $|\nodesT| = 183$, $|\edgesT| = $33,306, $|\nodesB| = 117$, $|\edgesB| = $1,044.}
    \label{fig:demand_realization}
\end{figure}

\Cref{fig:demand_realization} visualizes the resulting input demand of the probabilistic demand realization and the low-demand truncation with rescaling for each city, as well as their input networks. Note that the low-demand truncation with rescaling case yields anomalously low demand for Atlanta due to the data artifact. It only has 8,928 trips across 306 ODs over a $3$-hour peak period, which is insufficient for building transit systems. We therefore omit the subsequent experiments for this case. Across cities, the probabilistic demand realization covers a geographically wider demand area whereas the low-demand truncation with rescaling focuses on high-demand transportation areas. This pattern is consistent with the intended design of the two schemes. Despite the  differences in demand patterns between the two schemes, the resulting \bus~networks show relatively similar \bus-node selections. The similarity between bus networks suggests that both schemes identify comparable high-traffic corridor locations as \bus~nodes based on the facility location selection procedure.

To establish a realistic baseline for our ridership maximization MIP, we first solve a cost-minimization variant of the model. In this setup, we estimate the budget required to satisfy $90\%$ of total demand, keeping all other parameters constant. We run $160$ iterations, each outputting up to five solutions if the solver can find five, for the cost-minimization model yielding a baseline cost for each city under each demand setting. Under the probabilistic demand realization, the baseline costs for Boston, Chicago, and Atlanta are $\$369$K, $\$162$K, and $\$141$K, respectively. Under the low-demand truncation with rescaling setting, the baseline costs for Boston and Chicago are $\$233$K and $\$153$K. Using the baseline budget, we then conduct ridership maximization experiments across five budget levels: the full baseline (100\%) and four reduced budgets representing 90\%, 80\%, 70\%, and 60\% of the baseline, respectively. We outline the procedure for obtaining each baseline budget in \Cref{subsec:experiments-budget}.

\subsection{Optimization Setting}\label{subsec:experiments-solver} 
All computational experiments were conducted on a university high-performance computing (HPC) cluster using CPU resources. Each optimization job was allocated 32 CPU cores and 128 GB of memory. All optimization models were solved using Gurobi Optimizer 12.0.1 (Gurobi Optimization, LLC) through the Python API. During each column generation iteration, the search is focused on a reduced network subgraph as explained, with edge-trimming parameter $p=0.3$. Each pricing problem is solved with a 10-minute wall-clock time limit or until the solver reaches a 5\% gap termination criterion. In each iteration, the solver returns the five best feasible solutions (i.e., five lines) it identifies, if the solver can find five such solutions. The ridership-maximization line-generation process terminates after $120$ iterations. Using this set as input, we solve $(\mip)$ under each budget value for each system. For the multi-modal system, we apply the heuristic for the MIP master problem to first select 200 lines (using parameter $m=200$ in the algorithm) and then remove 10 lines at a time until only $50$ lines remain. Then $(\mip)$ is solved with a 0\% gap termination criterion on this reduced selection. The master problem heuristic completes within $1$ hour for all instances. For the resulting multi-modal systems, the final designs typically use between 30 and 40 active bus lines. To obtain MIP solutions for the bus-only system, we solve the model with its set of all candidate lines with a 0\% gap termination criterion, except for Chicago, where a 2\% gap termination criterion is used, due to high computational cost. All linear programs in our approach are solved to optimality. Additionally, we provide further experimental details in the appendix, including model input parameters (\Cref{subsec:experiments-model-parameter}), budget values (\Cref{subsec:experiments-budget}) and  algorithmic settings (\Cref{subsec:experiments-algo-setting}).

\subsection{Ridership Maximization through Multi-Modal Integration}
\begin{figure*}[!htb]
\centering
\includegraphics[width=0.93\textwidth]{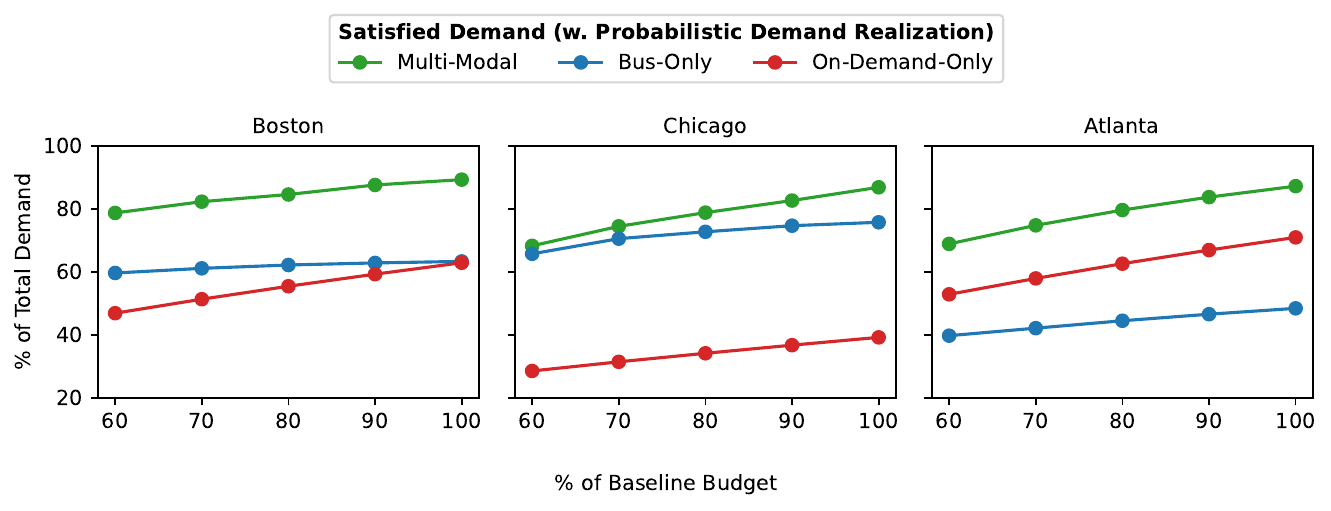}
\caption{Satisfied demand of the multi-modal system, \bus-only system, and \taxi-only system in each city using the probabilistic demand realization. The multi-modal system (green) increases from 78.76\% to 89.34\% in Boston, 68.32\% to 86.91\% in Chicago and 68.93\% to 87.27\% in Atlanta as the budget value changes from the lowest to the highest in each city. The \bus-only system (blue) increases from 59.69\% to 63.37\% in Boston, 65.78\% to 75.79\% in Chicago and 39.78\% to 48.47\% in Atlanta as the budget value changes from the lowest to the highest in each city. The \taxi-only system (red) increases from 46.92\% to 62.94\% in Boston, 28.58\% to 39.25\% in Chicago and 52.92\% to 71.00\% in Atlanta as the budget value changes from the lowest to the highest in each city. For a fixed budget, our multi-modal system achieves ridership percentages up to 40.97\% and 67.84\% higher than those of the bus-only and on-demand-only systems, respectively, in Boston; 14.67\% and 139.01\% higher in Chicago; and 80.06\% and 30.25\% higher in Atlanta.} 
\label{fig:sparse-results}
\end{figure*}

\begin{figure*}[!h]
\centering
\includegraphics[width=0.6\textwidth]{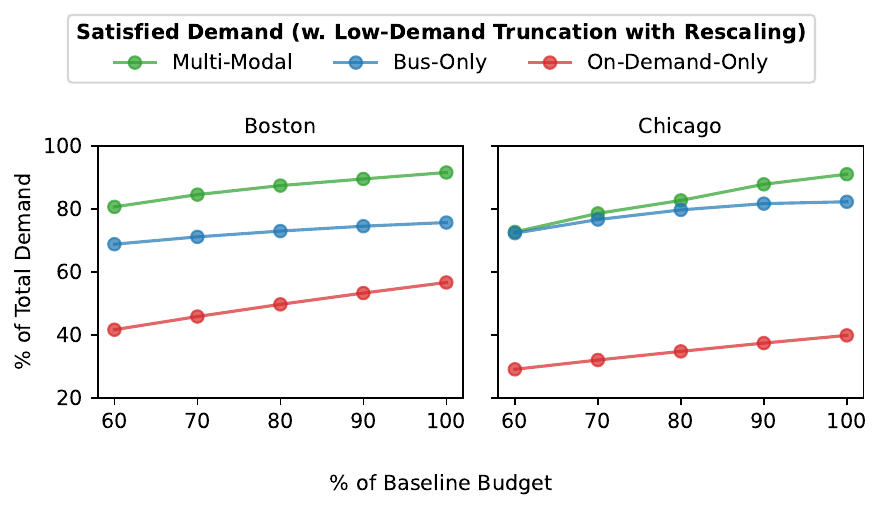}
\caption{Satisfied demand of the multi-modal system, \bus-only system, and \taxi-only system in each city using the low-demand truncation with rescaling. The multi-modal system (green) increases from 80.69\% to 91.61\% in Boston and 72.73\% to 91.06\% in Chicago as the budget value changes from the lowest to the highest in each city. The \bus-only system (blue) increases from 68.83\% to 75.71\% in Boston and 72.34\% to 82.31\% in Chicago as the budget value changes from the lowest to the highest in each city. The \taxi-only system (red) increases from 41.68\% to 56.68\% in Boston and 29.11\% to 39.87\% in Chicago as the budget value changes from the lowest to the highest in each city. For a fixed budget, our multi-modal system achieves ridership percentages up to 20.99\% and 93.58\% higher than those of the bus-only and on-demand-only systems, respectively, in Boston, and 10.63\% and 149.84\% higher in Chicago.} 
\label{fig:dense-results}
\end{figure*}

\Cref{fig:sparse-results} and \Cref{fig:dense-results} present the ridership results obtained using our design approach for each city, with the probabilistic demand realization and the low-demand truncation with rescaling, respectively. For all cities shown in both figures, the multi-modal systems achieve the highest level of satisfied demand, outperforming the \bus-only and the \taxi-only systems. The relative performance of the \bus-only system and the \taxi-only system varies across cities due to different demand distributions in each figure. For example, based on panels (a), (c), and (e) of \Cref{fig:demand_realization}, the average per-OD demand for Boston, Chicago, and Atlanta is $22.48$, $34.58$, and $13.46$, respectively. Higher average per-OD demand is reflected in the larger increase in satisfied demand for the \bus-only systems relative to the \taxi-only systems across all cities shown in \Cref{fig:sparse-results}. Greater per-OD demand potentially enables more efficient passenger aggregation and thus more effective \bus~utilization. A similar analysis based on panels (b) and (d) of \Cref{fig:demand_realization} reveals the same pattern in the comparison between the \bus-only and \taxi-only systems shown in \Cref{fig:dense-results}. The average per-OD demands for Boston and Chicago are 44.27 and 65.28, respectively. Since both schemes yield qualitatively consistent conclusions, we focus on the low-demand truncation with rescaling for the remaining experiments and omit Atlanta due to its low demand, as explained in \Cref{sec:data}.

\Cref{fig:dense-results} shows the percentage of demand satisfied by the multi-modal, \bus-only, and \taxi-only systems across five budget levels. While the multi-modal system leverages \taxi~for first- and last-mile connectivity, the \bus-only system is restricted to passengers whose origins and destinations are at \bus~stops. We compare the systems to demonstrate the advantages of integrating the two modes within a unified design.

The experimental results across both cities demonstrate the consistent superiority of the multi-modal system, highlighting the value of integrating \taxi~services with fixed-route \bus~networks. Across all budget values in Boston, the multi-modal system achieves ridership increases of up to 20.99\% over the \bus-only baseline and 93.58\% over the \taxi-only baseline. These gains quantify the cost of optimizing the modes independently. While efficient at serving high-density corridors, \bus-only systems leave broader demand unmet due to routing inflexibility. Conversely, \taxi-only systems suffer from high per-passenger costs that severely limit their scalability. By effectively bridging these gaps, the multi-modal system satisfies between 80.69\% and 91.61\% of total demand in Boston, compared with just 68.83\% to 75.71\% for the \bus-only system.

In \Cref{fig:dense-results}, Chicago exhibits similar trends, with multi-modal ridership gains of up to 10.63\% and 149.84\% over the \bus-only and \taxi-only baselines, respectively. The comparatively smaller gain over the \bus-only system in Chicago relative to Boston reflects distinct structural differences in the underlying transit networks and demand distributions of the two cities. Nevertheless, the significant gain over the \taxi-only baseline in Chicago reaffirms that on-demand services alone are ill-suited as a standalone transportation solution at scale for serving aggregate downtown peak demand.

\begin{figure}[!htb]
\centering
\includegraphics[width=0.55\textwidth]{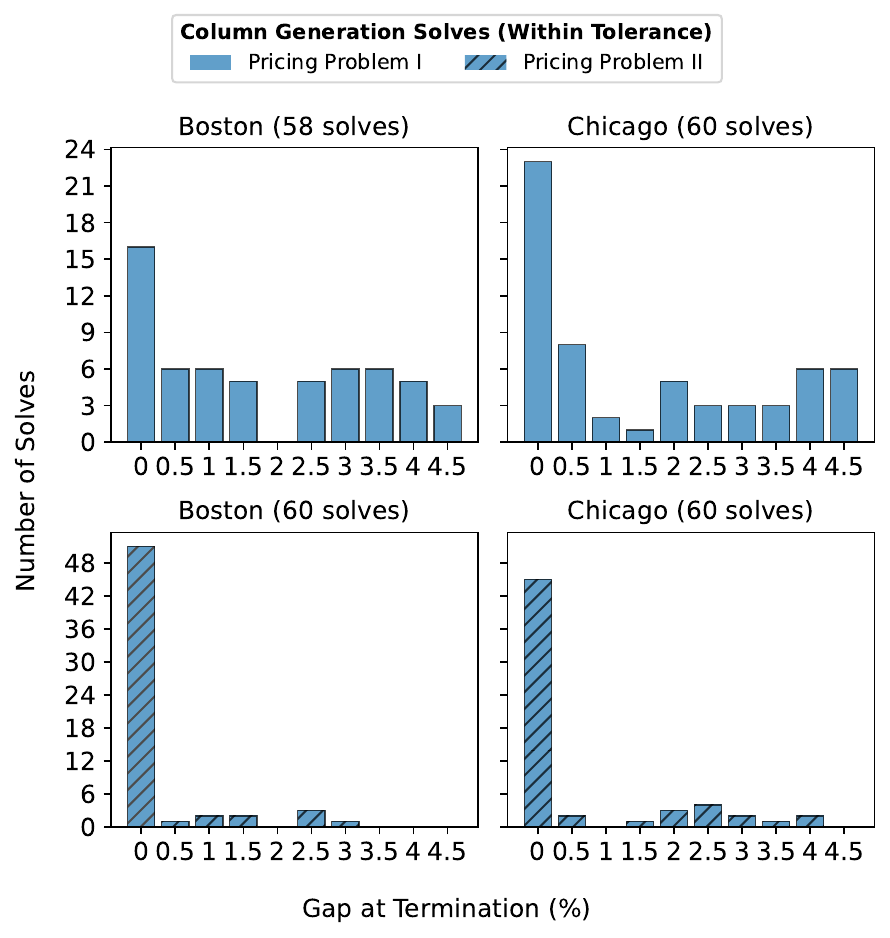}
\captionsetup{skip=1pt}
\caption{Histogram of pricing problem termination gaps. Most pricing problem instances terminate within the 5\% tolerance specified by the Gurobi MIPGap parameter before reaching the 10-minute time limit.}
\label{fig:subprob}
\end{figure}

To quantify solution quality within the final candidate set of 50 lines, we compare the MIP solution with the LP relaxation over the same set at the highest budget level. For Boston, the LP relaxation attains 92.18\% ridership versus 91.61\% for the MIP; for Chicago, 91.53\% versus 91.06\%. The remaining within-set integrality gap is therefore small. Because the pricing problems are solved heuristically over reduced subgraphs, the procedure does not certify optimality over the full set of feasible lines. Histograms of pricing problem termination gaps (\Cref{fig:subprob}) show that nearly all solves terminate within the 5\% tolerance before the 10-minute limit; Pricing Problem II terminates with smaller gaps, reflecting its reduced complexity, potentially at the cost of solution quality.

\subsection{Sensitivity Analysis on On-Demand Cost}
\begin{figure*}[!h]
\centering
\includegraphics[width=0.6\textwidth]{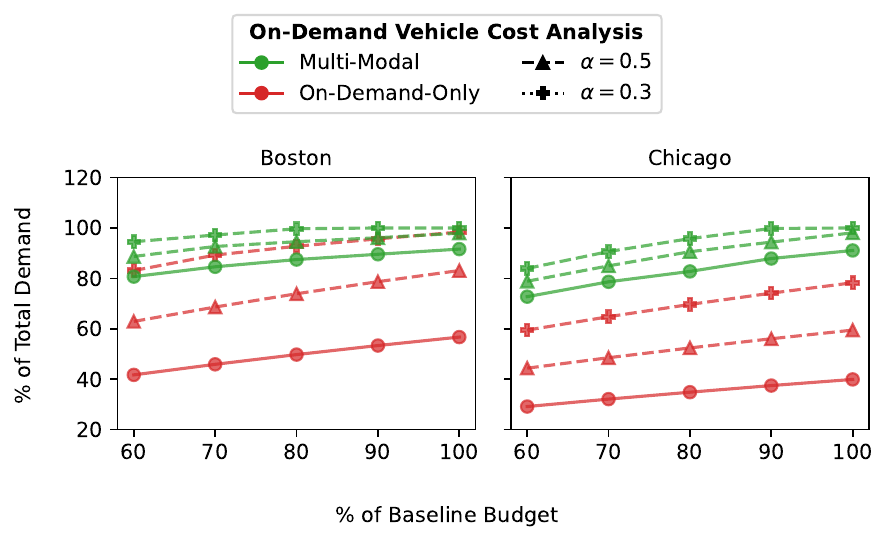}
\caption{Sensitivity analysis: \taxi~cost. Lower \taxi~costs yield substantial ridership improvements in both cities. At $\scaleTaxi=0.3$, satisfied ridership in the multi-modal system increases from 94.54\% to 100\% in Boston and from 83.94\% to 100\% in Chicago; while at $\scaleTaxi=0.5$, it rises from 88.68\% to 97.93\% in Boston and from 78.85\% to 98.07\% in Chicago. Relative to the default baseline ($\scaleTaxi=1$) in the multi-modal system, the settings $\scaleTaxi=0.3$ and $\scaleTaxi=0.5$ achieve maximum ridership gains of 17.16\% and 9.91\% in Boston, and 15.75\% and 9.46\% in Chicago, respectively. The $\scaleTaxi=0.3$ gains do not appear proportionally larger than the $\scaleTaxi=0.5$ gains because satisfied ridership reaches the 100\% saturation limit in both cities at $\scaleTaxi=0.3$.}
\label{fig:dense-cost}
\end{figure*}

\begin{figure*}[!h]
\centering
\includegraphics[width=0.6\textwidth]{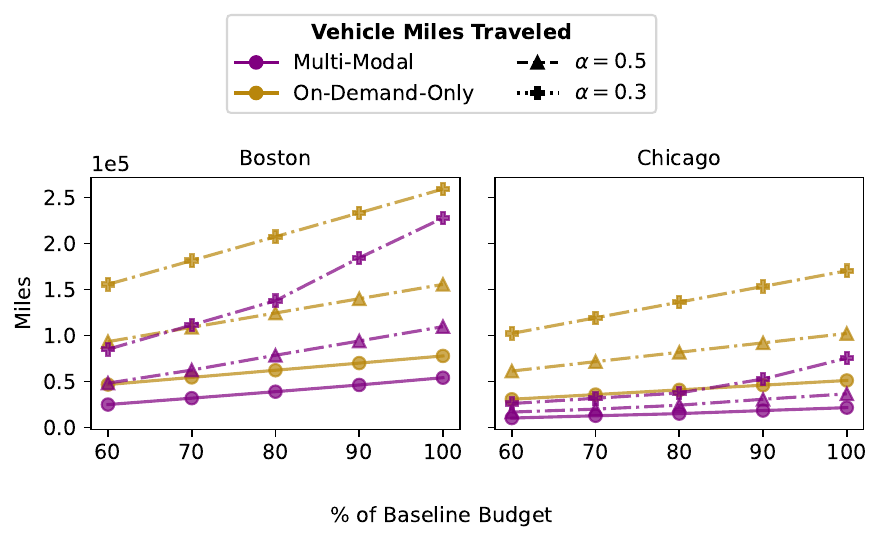}
\caption{Vehicle miles traveled for
$\scaleTaxi$ = 0.3, 0.5, and 1. The total vehicle miles traveled is the sum of \taxi~vehicle miles traveled and \bus~vehicle miles traveled. The total VMT increases as the budget increases and as $\scaleTaxi$ decreases. For Boston, as the budget increases, the multi-modal system VMT increases from $84$K to $227$K miles, from $48$K to $109$K miles, and from $24$K to $54$K miles for $\scaleTaxi$ = 0.3, 0.5, and 1; the \taxi-only system VMT increases from $155$K to $259$K miles, from $93$K to $155$K miles, and from $46$K to $77$K miles for $\scaleTaxi$ = 0.3, 0.5, and 1. For Chicago, as the budget increases, the multi-modal system VMT increases from $25$K to $75$K miles, from $16$K to $36$K miles, and from $10$K to $21$K miles for $\scaleTaxi$ = 0.3, 0.5, and 1; the \taxi-only system VMT increases from $102$K to $170$K miles, from $61$K to $102$K miles, and from $30$K to $51$K miles for $\scaleTaxi$ = 0.3, 0.5, and 1. }
\label{fig:dense-cost-VMT}
\end{figure*}

As commercial autonomous vehicles and ride-pooling services grow in prevalence, their integration with \bus~lines becomes an inevitable dimension of urban mobility planning. Transit agencies are actively exploring how \bus~systems could benefit from integrating with a variety of \taxi~services including shared and autonomous services~\cite{ruter}. \Cref{fig:dense-cost} presents sensitivity analyses on the \taxi~cost parameter ($\alpha$). We use distinct cost parameters to represent the varying cost profiles for traditional \taxi~services, ride-pooling services, and autonomous vehicles---with lower per-passenger costs reflecting savings from shared rides and the elimination of driver wages (in AVs). Our results indicate that as autonomous vehicles and ride-pooling drive down operational costs, the multi-modal system achieves significantly higher demand satisfaction within the same budget. However, the multi-modal system with higher \taxi~costs still generally outperforms a cheaper \taxi-only system. \Cref{fig:dense-cost-VMT} demonstrates that the multi-modal design realizes these higher ridership levels while simultaneously reducing total vehicle miles traveled (VMT) relative to the \taxi-only system: moving more people with fewer vehicle-miles than an \taxi-only system. The reduction in VMT stems from the different roles of \taxi~vehicles in the \taxi-only and multi-modal systems. In the \taxi-only system, all trips are served by \taxi~vehicles. In contrast, the multi-modal system uses \taxi~services for direct trips and flexible first- and last-mile connections, while high-capacity \bus~corridors carry concentrated passenger flows. This functional separation enables more efficient allocation of travel demand to vehicles and therefore reduces VMT.

\subsection{Comparison with the Benchmark Approach}

\begin{figure}[!htb]
\centering
\includegraphics[width=0.55\textwidth]{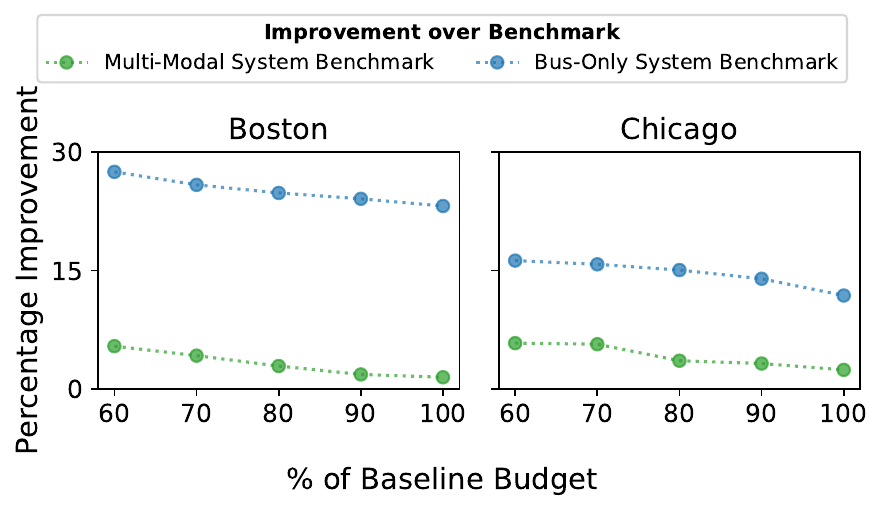}
\captionsetup{skip=1pt}
\caption{Percentage improvement over benchmark method: cases for the multi-modal system and the bus-only system.}
\label{fig:benchmark}
\end{figure}

To rigorously evaluate our line-generation procedure, \Cref{fig:benchmark} compares our method against a benchmark multi-modal line-generation heuristic based on randomized shortest paths, as used in \cite{Perivier2021}. We use this heuristic as a benchmark because it is the most recent line-generation method for the same problem setting. The heuristic in \cite{Perivier2021} constructs candidate lines by randomly selecting four nodes from the \bus~network and connecting them with shortest paths. We additionally require the generated lines to satisfy our detour, maximum-length, and subtour-elimination constraints.

Our optimization-based approach improves ridership by 5.42\% at the lowest budget to 1.48\% at the highest budget in Boston, and by 5.80\% at the lowest budget to 2.42\% at the highest budget in Chicago, relative to the benchmark. This improvement is most pronounced under tight budget constraints where efficient resource allocation is critical---precisely the conditions in which our method provides the greatest practical value. Notably, when evaluating the \bus-only case, our method improves ridership by 27.49\% at the lowest budget to 23.17\% at the highest budget in Boston and by 16.26\% at the lowest budget to 11.83\% at the highest budget in Chicago, relative to the benchmark. The bus-only comparison is the cleanest test of line-generation quality, as it isolates the generated lines from the on-demand component. The substantial gain, achieved without any \taxi~support, underscores the superior efficiency of our underlying line-generation algorithm.

\section{Conclusion}
Motivated by the inherent limitations of modern mass transit in the face of emerging ride-hailing services, we have developed a framework for designing multi-modal systems that systematically integrate these two modes. Our approach introduces a scalable flow-based mixed-integer programming formulation coupled with a tailored column generation heuristic scheme. Unlike previous methods, this framework explicitly accounts for realistic operational features, such as passenger waiting times and a structural guarantee that passengers do not transfer between \bus~lines, while remaining computationally tractable for urban-scale networks.

Our computational results demonstrate three main findings. First, the framework significantly increases ridership compared to single-mode systems within the same fixed budget, achieving gains of up to 20.99\% in Boston and 10.63\% in Chicago over traditional bus-only systems. Second, the multi-modal design improves the efficiency of existing transit backbones, effectively using \taxi~services to bridge the first- and last-mile gaps that often render traditional transit inaccessible in sprawling urban environments. Third, our optimization-based approach consistently outperforms the benchmark line-generation method, particularly in designing the high-capacity bus routes that form the core of the system.

As cities face growing pressure to reduce congestion and vehicle miles traveled (VMT), our results show that integrated, data-driven planning offers a practical pathway to more efficient and resilient urban mobility. By providing a scalable tool for central planners, this research helps devise transportation systems that are not only cost-effective but also capable of providing ubiquitous access to personal mobility.

\newpage
\bibliographystyle{naturemag}
\bibliography{reference}

\newpage
\appendix
\section{Additional Experimental Details} \label{app:appendix}

\subsection{Experimental Setup Details}

\subsubsection{Model Parameters}\label{subsec:experiments-model-parameter}We present the default parameter values employed in the experiments for the main MIP model and the subproblems. An \taxi~vehicle traveling on edge $(\nodeU, \nodeV) \in \edgesT$ incurs a transportation cost $\scaleTaxi\costUV$, where $\costUV$ is the shortest-path distance between $\nodeU\in \nodesT$ and $\nodeV\in \nodesT$ in the original road network and the default value of parameter $\scaleTaxi$ is 1. We assume the \taxi~mode is served by conventional taxis in the baseline experiments. Taxi vehicles have a capacity of 1, whereas bus vehicles have a capacity of $\capacity=50$. We assume a maximum passenger waiting time of $T=15$ minutes, during which a bus travels $R = $ 4,000 meters \cite{bus_travel_per_period}.  Each line $\oneline$ (a bidirectional simple path) has total length $M_l=\sum_{(\nodeU,\nodeV)\in \edgesline}\costUV$. Based on these values, we construct the feasible set of the number of vehicles $\freqL$ for each line,  the per-vehicle maximum passenger load $\capacity_{\oneline}$ and per-vehicle operating cost $\scaleBus\costL$. According to \cite{bus_per_vehicle_mile}, the national total operating expenses and vehicle revenue miles for \bus es in 2023 are \$28,016.9 million and 1,863.3 million. Therefore, we estimate the \bus~operating cost per mile per vehicle to be \$28,016.9/1,863.3 $\approx \$15$. The median hourly wage for taxi drivers nationwide is $\$16.67$ \cite{taxi_driver_wage} and the driving speed is estimated to be $12$ miles per hour \cite{Huang2018}. Hence, the estimated wage for drivers is $\$16.67 / 12 = \$1.39$ per mile. We estimate that the average cost to own and operate a new vehicle is $\$0.82$ per mile \cite{own_vehicle}. Therefore, for a company managing taxi services, we estimate the per-mile, per-vehicle operating cost to be the sum of $\$1.39$, $\$0.82$, and other regulatory costs. We round up to $\$3$ per mile to account for operator margin and regulatory costs. In summary, we estimate the default per-mile cost of a \bus~vehicle to be $\$15$ and the default per-mile cost of a taxi vehicle to be $\$3$. These two figures are measured at the level at which a transit agency incurs them: on-demand service is typically procured from an operator at a contracted per-mile rate, for which per-mile wage and vehicle costs are the relevant components, while bus service is operated in-house, for which operating expense per revenue mile is the budget-relevant quantity. The budget constraint therefore prices each mode as an agency would purchase it. Accordingly, the physical cost estimates are $\$3$ per vehicle-mile for on-demand service and $\$15$ per vehicle-mile for bus service. We normalize these values by the on-demand cost, yielding $\scaleTaxi = 1$ and $\scaleBus = 5$ in the experiments.

For the pricing problem, the detour constraint parameter $\detourConst$ is set to $2$, and the maximum path length parameter $\maxLength$ is set to 20,000 meters. At least one of the first-mile and last-mile taxi trips for a demand request in the multi-modal system must be within a maximum distance of 1,000 meters, reducing the number of decision variables in $(\pricing)$, $(\pricingPrime)$, $(\mip)$ and $(\mipPrime)$. This is a deliberate modeling choice: we require only that one end of the trip have a short connection to the \bus~network, rather than bounding both legs, so as not to exclude multi-modal trips in which one leg is short and the other longer. The budget constraint \eqref{eqn:Budget} also weighs against unnecessarily long taxi legs, since each taxi vehicle-mile draws against the same operating budget. The next section details the procedure for obtaining the input operating budget value $B$. 

\subsubsection{Operating Budget}\label{subsec:experiments-budget}
To obtain the total budget value, we run our algorithm on the cost-minimization version of the model for the multi-modal transit network planning problem to estimate the budget required to satisfy at least 90\% of total demand over the 3-hour peak period. To convert $(\mip)$ from the ridership-maximization version to the cost-minimization version, we modify it as follows.

\begin{enumerate}
\item[(i)] Replace \objcref{eqn:Ridership} with the cost-minimization objective $\min \sum_{(\nodeU,\nodeV)\in \edgesT} \scaleTaxi\costUV \numTaxi_{\nodeU,\nodeV} + \sum_{\oneline\in \lineSet} \scaleBus \costL \numBus_{\oneline}$.
\item[(ii)] Substitute \constrscref{eqn:Request} with $\sum_{(\ori,\dest)\in \requestSet}\directTrip_{\ori,\dest} + \sum_{(\ori,\dest)\in \requestSet}\sum_{\nodeV\in \nodesB} \toDest^{\ori}_{\nodeV,\dest} \geq 0.9\sum_{(\ori,\dest)\in \requestSet} \requestST$, ensuring at least 90\% of total demand is satisfied.
\item[(iii)] Remove \constrcref{eqn:Budget}.
\item[(iv)] Restate all inequality constraints in $\geq$ form, the conventional form for a minimization program.
\end{enumerate}

Applying the analogous set of modifications (i)--(iv) to $(\mipPrime)$ yields its cost-minimization counterpart. Based on the cost-minimization version of $(\mip)$ and $(\mipPrime)$, we replace the objective functions of $(\pricing)$ and $(\pricingPrime)$, respectively, with
 \begin{align}\min~~& \frac{\scaleBus}{\capacity}\sum_{(\nodeU,\nodeV)\in \edgesB} \left(\costUV+ c_{v,u}\right)\subEdge_{\nodeU,\nodeV} +\sum_{\ori\in \nodesOri} \sum_{\nodeU\in \nodesB} \tilde \taxiDual_{\ori,\nodeU}\getOn^{\ori}_{\nodeU} + \sum_{\ori\in \nodesOri}\sum_{\nodeU\in \nodesB}  \tilde \busDual_{\ori,\nodeU}\getOff^{\ori}_{\nodeU}\label{eqn:Sub1-mincost}\\
 \intertext{and}
 \min~~& \frac{\scaleBus}{\capacity}\sum_{(\nodeU,\nodeV)\in \edgesB}\left(\costUV+  c_{\nodeV, \nodeU}\right)\subEdge_{\nodeU, \nodeV} - \sum_{(\nodeU,\nodeV)\in \edgesB} (\tilde \capDual_{\nodeU,\nodeV} +\tilde \capDual_{\nodeV,\nodeU})\subEdge_{\nodeU, \nodeV}.\label{eqn:Sub2-mincost}
 \end{align}
 
We detail below the differences between the ridership-maximization and cost-minimization pricing objective functions (i.e., \eqref{eqn:Sub1-maxride} and \eqref{eqn:Sub2-maxride} vs. \eqref{eqn:Sub1-mincost} and \eqref{eqn:Sub2-mincost}) in terms of the dual variables of their respective master LPs. For any primal decision variable, its contribution to the reduced-cost objective is $c_j-\sum_i a_{ij}\tilde\pi_i$, where $c_j$ is its coefficient in the master problem’s original objective, $a_{ij}$ is its coefficient in constraint $i$, and $\tilde\pi_i$ is the corresponding dual value from the master LP; the sum is taken over the constraints i in which the variable appears. This reduced-cost rule is independent of whether the master is a maximization or minimization problem.

\textbf{For Pricing Problem I $(\pricing)$:}
\begin{itemize}
\item \Taxi~dual $\tilde\taxiDual$, from \eqref{eqn:Taxi}. Since \eqref{eqn:Taxi} is restated in (iv), the coefficient of $\getOn^{\ori}_{\oneline,\nodeV}$ changes from $+1$ (maximization form) to $-1$ (minimization form). The rule above gives a contribution of $-\tilde\taxiDual\getOn$ in the ridership-maximization pricing objective and $+\tilde\taxiDual\getOn$ in the cost-minimization pricing objective.

\item Alighting-balance dual $\tilde\busDual$, from \eqref{eqn:Agg}. Since \eqref{eqn:Agg} is not restated, the coefficient of $\getOff^{\ori}_{\oneline,\nodeV}$ remains $-1$ in both master formulations. The rule above gives a contribution of $+\tilde\busDual\getOff$ in both the ridership-maximization and cost-minimization pricing objectives.
\end{itemize}

\textbf{For Pricing Problem II $(\pricingPrime)$:}
\begin{itemize}
\item Aggregated capacity dual $\tilde\capDual$, from \eqref{eqn:AggCapa}.
The edge-selection variables $\subEdge_{\nodeU,\nodeV}$ in $(\pricingPrime)$ are a proxy for $\numBus_{\oneline}$, the variable that actually appears in \eqref{eqn:AggCapa}. In the ridership-maximization master, $\numBus_{\oneline}$ has zero direct objective coefficient and enters \eqref{eqn:AggCapa} (in $\leq$ form) with a negative coefficient per line edge, giving a reduced-cost contribution of $+\tilde\capDual$ per edge. A positive reduced cost is what an improving column requires in a maximization master, so this term is added in \eqref{eqn:Sub2-maxride}. In the cost-minimization master, $\numBus_{\oneline}$ carries a direct objective cost of $\scaleBus\costL$, and \eqref{eqn:AggCapa} is restated in $\geq$ form by (iv), changing its constraint coefficient to positive and its reduced-cost contribution to $-\tilde\capDual$ per edge. Since an improving column requires a negative reduced cost in a minimization master, and the pricing searches for the most negative value by minimizing, this term is subtracted in the cost-minimization pricing objective. 
\end{itemize}

\textbf{For both $(\pricing)$ and $(\pricingPrime)$:}
\begin{itemize}
\item Budget dual $\tilde\budgetDual$, from \eqref{eqn:Budget}. The corresponding constraint is removed by (iii), so no corresponding dual-weighted term appears in either cost-minimization pricing objective. The line construction cost $\frac{\scaleBus}{\capacity}(\costUV+c_{\nodeV,\nodeU})$ enters directly into the cost-minimization pricing objectives.
\end{itemize}

We apply our approach to the cost-minimization version of $(\mip)$ to obtain an MIP solution. Specifically, to generate the candidate set of lines for the cost-minimization problem, we create up to 800 lines (160 column generation iterations, each identifying up to 5 lines if the solver can find five) for each of Boston, Chicago, and Atlanta. The LP objective values plateau in the later stages of line generation. Under the probabilistic demand realization, the LP objective values for Boston, Chicago, and Atlanta plateau at \$251,585, \$117,661, and \$45,585, respectively. We then solve the cost-minimization MIP model with all the generated lines as the candidate set of lines, with a termination gap of $12\%$. The MIP objective values for Boston, Chicago, and Atlanta are \$369,046, \$162,427, and \$141,902, respectively. Under the low-demand truncation with rescaling, the LP objective values for Boston and Chicago plateau at \$121,518 and \$91,131, respectively. To obtain the final MIP solution for each case, we apply \Cref{alg:mip-heur} to iteratively remove lines. We first select 200 lines (using parameter $m=200$ in the algorithm) and then remove 10 lines at a time until only $50$ lines remain (i.e., $m = 190, 180,..., 50$).  The MIP objective values for Boston and Chicago are \$233,145 and \$153,207, respectively. We note that distinct approaches were employed to obtain the MIP values for the probabilistic demand realization and for the low-demand truncation with rescaling. Both approaches were evaluated with a solver time limit of 5 hours, and the best results are reported.

\subsubsection{Algorithmic Setting}\label{subsec:experiments-algo-setting}
We apply our approach to design three types of systems to maximize ridership given the input demand using our model $(\mip)$: the multi-modal system, the \bus-only system, and the \taxi-only system. In both the multi-modal system and the \bus-only system, \bus~lines are constructed using \Cref{alg:full-alg}. The multi-modal system aims to serve all passengers by combining different modes, whereas the \bus-only system serves only those whose origin node and destination node are in $\nodesB$. In the \taxi-only system, no \bus~lines are present, and all passengers are served by \taxi~vehicles.

In \Cref{alg:full-alg}, some steps involve solving the LP relaxation of $(\mipPrime)$ to obtain the dual solutions. The dual variables for \constrscref{eqn:AggCapa} depend on edges in $\edgesB$ covered by the constructed lines. To ensure valid dual values in the process, we initialize the algorithm with an initial candidate set of lines that cover all edges in $\edgesB$. The initial candidate sets contain $56$, $33$, and $49$ lines for Boston, Chicago, and Atlanta, respectively, under the probabilistic demand realization, and $55$ and $50$ for Boston and Chicago, respectively, under the low-demand truncation with rescaling. Each line is a bidirectional simple path, constructed greedily by iteratively selecting adjacent edges from set $\edgesB$. Specifically, each line is formed by starting from a randomly chosen edge and repeatedly adding adjacent edges until no further extension is possible. Since the path is bidirectional, each selected edge is added as a pair of directed edges in opposite directions. Covered edges are removed from $\edgesB$ and we continue to construct lines in this way to cover the remaining ones in $\edgesB$. Each line in this initial candidate set is assigned a per-vehicle operating cost of $100\scaleBus\costL$ in the multi-modal system and $1000\scaleBus\costL$ in the \bus-only system, compared to the standard cost $\scaleBus\costL$. The higher cost serves to penalize infeasibilities introduced by the greedy construction process, which may violate the detour constraints in the subproblems. As a result, these lines typically do not appear in the optimal solutions in the LP relaxation of $(\mip)$ and $(\mipPrime)$ during the line generation process, until 100 additional lines have been generated beyond the initial candidate set.

Our approach uses pairs of budgets $\budget_{LP}$ and $\budget_{MIP}$. For both the multi-modal and \bus-only systems in each city, we generate lines for $40$ iterations using $\budget_{LP}$ from Budget Pair 1, followed by $20$ iterations each from Budget Pairs 2 through 5, also using $\budget_{LP}$. The resulting lines from the total $120$ iterations form the candidate set for each multi-modal and \bus-only system. Using this set as input, we solve $(\mip)$ under 
$\budget_{MIP}$ from each of the five budget pairs, yielding five MIP solutions per system. 


\subsection{Additional Computational Results}

\subsubsection{Sensitivity Analysis on Bus Cost and Capacity}

\begin{figure}[H]
\centering
\includegraphics[width=0.6\textwidth]{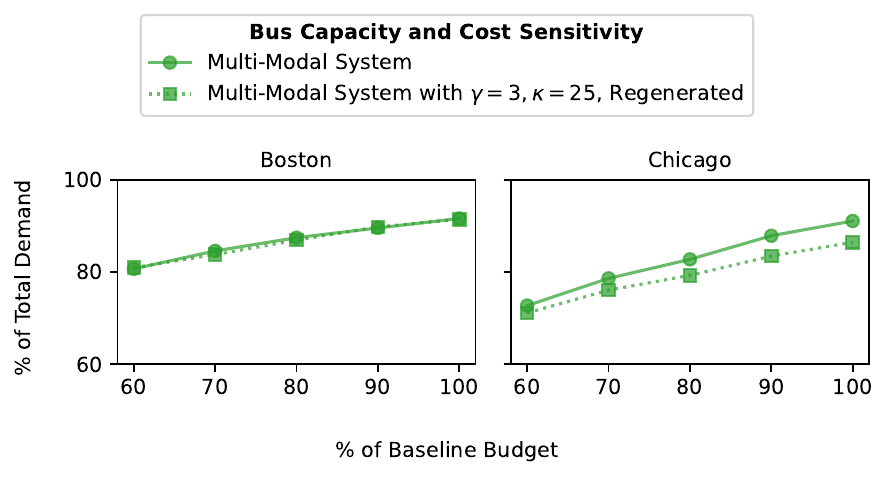}
\caption{Sensitivity analysis on \bus~capacity and cost in the multi-modal system. While the two parameter settings achieve similar results in Boston, the default system ($\scaleBus=5$ and $\capacity=50$) outperforms the alternative ($\scaleBus=3$ and $\capacity=25$) by up to 5.33\% in Chicago.}
\label{fig:dense-capacity-g3}
\end{figure}

 We conduct a sensitivity analysis of \bus~capacity and cost to understand the resource utilization. We use the default \bus~cost parameter $\scaleBus=5$ and per-vehicle \bus~capacity $\capacity=50$. In addition to the default setting, we investigate the effect of smaller bus vehicles with $\capacity=25$. Since operating costs do not scale linearly with vehicle size due to fixed costs (e.g., labor), we set $\scaleBus=3$ for $\capacity=25$. We regenerate the multi-modal system under these parameters.  With smaller capacity and a higher cost-to-seat ratio, the new multi-modal system offers greater flexibility in \bus~seat allocation, trading off higher per-seat costs. In \Cref{fig:dense-capacity-g3}, we compare the two multi-modal systems and find that neither consistently outperforms the other in Boston. For Chicago, the default setting with $\scaleBus=5$ and $\capacity=50$ achieves higher satisfied demand than the combination of $\scaleBus=3$ and $\capacity=25$, suggesting that higher vehicle capacity with lower per-seat costs is preferred to increased line flexibility. One explanation for this phenomenon is that Chicago has a higher per-OD demand (65.28) than Boston (44.27), indicating a more aggregated demand pattern that is well suited to high-capacity vehicles.

\begin{figure}[H]
\centering
\includegraphics[width=0.6\textwidth]{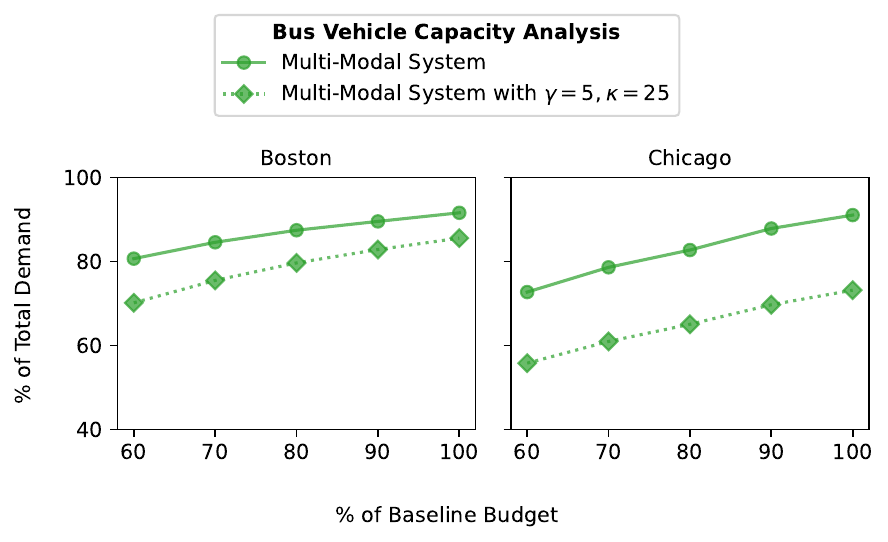}
\caption{Sensitivity analysis on decreasing \bus~capacity in the multi-modal system. Compared with the lower-capacity setting, the default setting achieves up to 14.99\% higher ridership in Boston and up to 30.34\% higher ridership in Chicago with the same candidate set of \bus~lines as input for the optimization program. }
\label{fig:dense-capacity-low}
\end{figure}

In another set of experiments, we investigate the capacity utilization in the default multi-modal system. We take the candidate set of all lines from the default multi-modal system as input for $(\mip)$. We solve the program with per-vehicle capacity set to $\capacity=25$ instead of $\capacity=50$, while maintaining the cost parameter $\scaleBus=5$. \Cref{fig:dense-capacity-low} shows the decrease in  ridership as the capacity drops from $50$ to $25$. We expect an increase in line utilization with the smaller capacity, with line utilization calculated as

\begin{equation}
\text{line-utilization} = \frac{\text{passenger-miles traveled}}{\text{total seat-miles available}}\times 100\%. 
\end{equation}

Indeed, the increases in line utilization when capacity decreases from $\capacity=50$ to $\capacity=25$ are 14.73\%, 21.28\%, 20.97\%, 21.90\%, 20.35\% for Boston, and 8.26\%, 8.46\%, 11.01\%, 10.56\%, 9.02\% for Chicago, from the lowest to the highest budget. With the smaller \bus~capacity, Boston experiences a smaller drop in total ridership and a greater increase in line utilization compared to Chicago. The higher per-OD demand in Chicago makes higher-capacity vehicles more advantageous. The higher the per-OD demand, the greater the advantage of investing in \bus~systems.

\subsubsection{Types of Trips in the Multi-Modal System}

We provide a description of our procedure for disaggregating passenger trips (hybrid, pure \bus, and pure \taxi) from a solution of model $(\mip)$ in \Cref{alg:disaggregate}. Note that the resulting disaggregation is not unique. 
 Specifically, the procedure extracts the number of pure \bus~trips from the solution. The number of hybrid trips is then readily obtained by subtracting both pure \bus~and pure \taxi~trips from the total. Let $\ttilde \numTaxi_{\nodeU,\nodeV}, \ttilde \numBus_{\oneline}, \ttilde \directTrip_{\ori,\dest}, \ttilde \getOn^{\ori}_{\oneline,\nodeV}, \ttilde \busFlow^{\ori}_{\oneline,\nodeU,\nodeV}, \ttilde \getOff^{\ori}_{\oneline,\nodeV}, \ttilde \toDest^{\ori}_{\nodeV,\dest}$ denote a feasible solution of $(\mip)$. Since each \bus~line forms a bidirectional simple path, we describe positional relations along the path as left and right.

 In \Cref{alg:disaggregate}, Lines 2--6 eliminate opposing passenger flows on anti-parallel edges from the solution. Lines 7--9 compute $F^{\ori}_{\dest}$, the total number of passengers traveling from $\ori$ to destination $\dest$ who alight from a \bus~at node $\dest$, with or without a first-mile \taxi~segment. For each passenger origin node $\ori$ and line $\oneline$ containing $\ori$, Line 11 computes $D^{\ori, \oneline}$, the passenger flow originating at $\ori$, boarding $\oneline$ at $\ori$, and traveling to the right of $\ori$ on the line. Assume $\busFlow^{\ori}_{\oneline,\ori-1,\ori}$ is zero if there is no node to the left of $s$. Similarly, assume $\busFlow^{\ori}_{\oneline,\ori,\ori+1}$ is zero if there is no node to the right of $s$.  For each node $b$ to the right of $\ori$ on line $\oneline$, Line 13 obtains the passenger flow originating at $\ori$ and arriving at $b$ on line $\oneline$, denoted by $E^{\ori, \oneline}_{b}$. Line 14 updates $D^{\ori, \oneline}$ to be the passenger flow with origin $\ori$ that boards $\oneline$ at $\ori$ and arrives at node $b$. Next, Line 15 computes the passenger flow on line $\oneline$ using pure \bus~trips, denoted by $Bus$, that board at $\ori$ and alight at b. We deduct $Bus$ from both $D^{\ori,\oneline}$ and $F^{\ori}_b$. We update the total number of pure \bus~trips extracted so far, denoted by $TotalBus$, which concludes one iteration of the extraction process of Lines 10--20. We then proceed to each subsequent node to the right along the line $\oneline$ to extract pure \bus~trips. Lines 21--31 compute pure \bus~trips corresponding to passenger flow on the line that boards at $\ori$ and alights at nodes to the left of $\ori$. In the end, we return the total pure \bus~trips extracted. 

\begin{algorithm}[!htb]
\caption{Disaggregating \Bus~Ridership from the Total Satisfied Demand in the Multi-Modal System}\label{alg:disaggregate}
\begin{algorithmic}[1]
\setlength{\baselineskip}{1.25\baselineskip}
\Require Solution of $(\mip)$ denoted by $\ttilde \numTaxi_{\nodeU,\nodeV}, \ttilde \numBus_{\oneline}, \ttilde \directTrip_{\ori,\dest}, \ttilde \getOn^{\ori}_{\oneline,\nodeV}, \ttilde \busFlow^{\ori}_{\oneline,\nodeU,\nodeV}, \ttilde \getOff^{\ori}_{\oneline,\nodeV}, \ttilde \toDest^{\ori}_{\nodeV,\dest}$

\State $TotalBus\gets 0$
\For{each origin $\ori\in \nodesOri$, each line $\oneline\in \lineSet$, and edge $(\nodeU, \nodeV)\in \edgesline$}
\State $x\gets \min \{\ttilde \busFlow^{\ori}_{\oneline, \nodeU, \nodeV}, \ttilde \busFlow^{\ori}_{\oneline, \nodeV,\nodeU}\}$
\State $\ttilde \busFlow^{\ori}_{\oneline,\nodeU,\nodeV}\gets \ttilde \busFlow^{\ori}_{\oneline,\nodeU,\nodeV}-x$
\vspace{-0.1em}
\State $\ttilde \busFlow^{\ori}_{\oneline,\nodeV,\nodeU}\gets \ttilde \busFlow^{\ori}_{\oneline,\nodeV,\nodeU}-x$ 
\EndFor
\vspace{-0.4em}
\For{each origin node $\ori\in \nodesOri$ and bus node $\dest\in \nodesB$}
\State $F^{\ori}_t\gets \min \{{\scriptstyle\sum\nolimits}_{\oneline\in \lineSet: \dest\in \nodesline} \ttilde \getOff^{\ori}_{\oneline,t}, \ttilde \toDest^{\ori}_{t,t}\}$, if $(\ori,\dest) \in \requestSet$; $F^{\ori}_t\gets 0$, otherwise 
\EndFor
\vspace{-0.4em}
\For{each origin node $\ori\in \nodesOri$ and line $\oneline \in \lineSet$ containing $\ori$}
\State $D^{\ori, \oneline}\gets \ttilde \busFlow^{\ori}_{\oneline,\ori,\ori+1} - \ttilde \busFlow^{\ori}_{\oneline,\ori-1,\ori}$, where $\ori+1$ and $\ori-1$ denote the nodes to the right and left of $\ori$
\For{each node $b$ to the right of $\ori$ on line $\oneline$}
\State $E^{\ori, \oneline}_{b}\gets \ttilde \busFlow^{\ori}_{\oneline,b-1,b}$, where $b-1$ denotes the node to the left of $b$
\State $D^{\ori, \oneline}\gets\min\{D^{\ori, \oneline}, E^{\ori, \oneline}_{b}\}$ 
\State $Bus \gets \min \left\{D^{\ori, \oneline}, F^{\ori}_b\right\}$ 
\State $D^{\ori, \oneline}\gets D^{\ori, \oneline}-Bus$
\State $F^{\ori}_b\gets F^{\ori}_b-Bus$
\State $TotalBus\gets TotalBus + Bus$
\EndFor
\EndFor

\For{each origin node $\ori\in \nodesOri$ and line $\oneline \in \lineSet$ containing $\ori$}
\State $D^{\ori, \oneline}\gets \ttilde \busFlow^{\ori}_{\oneline,\ori,\ori-1} - \ttilde \busFlow^{\ori}_{\oneline,\ori+1,\ori}$, where $\ori-1$ and $\ori+1$ denote the nodes to the left and right of $\ori$
\For{each node $b$ to the left of $\ori$ on line $\oneline$}
\State $E^{\ori, \oneline}_{b}\gets \ttilde \busFlow^{\ori}_{\oneline,b,b+1}$, where $b+1$ denotes the node to the right of $b$  
\State $D^{\ori, \oneline}\gets\min\{D^{\ori, \oneline}, E^{\ori, \oneline}_{b}\}$ 
\State $Bus \gets \min \{D^{\ori, \oneline}, F^{\ori}_b\}$ 
\State $D^{\ori, \oneline}\gets D^{\ori, \oneline}-Bus$
\State $F^{\ori}_b\gets F^{\ori}_b-Bus$
\State $TotalBus\gets TotalBus + Bus$
\EndFor
\EndFor\\
\Return $TotalBus$
\end{algorithmic}
\end{algorithm}

\end{document}